\documentclass[11pt]{article}
\usepackage{amsmath}
\usepackage{amssymb}
\usepackage{amsthm}
\numberwithin{equation}{section}
\newtheorem{theorem}{Theorem}[section]

\newtheorem{corollary}[theorem]{Corollary}

\newtheorem{definition}[theorem]{Definition}
\newtheorem{example}[theorem]{Example}
\newtheorem{exercise}[theorem]{Exercise}
\newtheorem{lemma}[theorem]{Lemma}

\newtheorem{proposition}[theorem]{Proposition}
\newtheorem{remark}[theorem]{Remark}

\def\sd{\partial}
\def\xyb{(\bar x,\bar y)}

\def\dis{\displaystyle}

\def\ran{\rangle}
\def\lan{\langle}

\def\Pi{{\cal P}}
\def\B{{\cal B}}

\def\N{\mathbb{N}}
\def\R{I\!\!R}
\def\la{\lambda}
\def\al{\alpha}
\def\del{\delta}
\def\ep{\varepsilon}
\def\sig
\def\sd{\partial}
\def\xb{\overline x}

\def\pb{\bar p}
\def\yb{\overline y}
\def\ub{\overline u}
\def\vb{\overline v}
\def\zb{\overline z}

\def\wb{\overline w}
\def\tb{\bar t}

\def\xzb{(\bar x,\bar z)}
\def\bpsi{\overline{\psi}}

\def\dom{{\rm dom}~}
\def\gr{{\rm Graph}~}

\def\rra{\rightrightarrows}

\def\xn{x_n}

\def\ep{\varepsilon}
\def\la{\lambda}

\def\al{\alpha}

\def\vf{\varphi}
\def\del{\delta}

\def\ga{\gamma}

\def\sig{\sigma}
\def\nab{\nabla}

\def\cf{{\mathcal F}}
\def\cll{{\mathcal L}}
\def\ca{{\mathcal A}}

\def\cd{{\mathcal D}}

\def\ce{{\mathcal E}}

\def\cp{{\mathcal P}}

\def\cq{{\mathcal Q}}

\def\cu{{\mathcal U}}

\def\ch{{\mathcal H}}
\def\cs{{\mathcal S}}

\def\ct{{\mathcal T}}

\def\cl{{\rm cl}}

\def\dom{{\rm dom}~}
\def\epi{{\rm epi}~}
\def\intr{{\rm int}~}

\def\im{{\rm Im}~}
\def\ker{{\rm Ker}~}
\def\conv{{\rm conv}~}

\def\ran{\rangle}
\def\lan{\langle}

\def\sur{{\rm sur}}

\def\cont{{\rm contr}}
\def\calm{{\rm calm}}
\def\reg{{\rm reg}}
\def\subreg{{\rm subreg}}
\def\rad{{\rm rad}}
\def\ress{{\rm ress}}

\def\Lip{{\rm lip}}
\def\lip{{\rm lip}}

\def\dis{\displaystyle}

\def\bo{\overset{\circ}B}

\title{ METRIC  REGULARITY -- A SURVEY}
\author{A.D. IOFFE \thanks{Department of Mathematics, Technion, Haifa 32000, Israel}}

\begin{document}
	
	
	
\maketitle

	\begin{flushright} \begin{tabular}{r}
			{\it In science things
				should be made}\\ {\it as simple as possible.}
			\\ {\rm Albert Einstein} \\	
			
			\\		
			
			{\it All the great things are simple.}\\ Winston Churchill
		\end{tabular}
	\end{flushright}
	
	\begin{abstract} Metric regularity theory lies in the very heart of
	variational analysis, a relatively new discipline whose appearance was to a large extent determined by needs of modern optimization theory in which
	such phenomena as non-differentiability and  set-valued mappings naturally appear.
	The roots of the theory go back to such fundamental results of the classical analysis
	as the implicit function theorem, Sard theorem and some others. The paper offers
	a survey of the state-of-the-art of some principal parts of the theory along with a variety of its applications in analysis and optimization.
	\end{abstract}
	
	\vskip 4mm

	\noindent{\bf \large Contents}
	
	\vskip 3mm
	
		\noindent{\bf Introduction}

	
	\vskip 2mm
	
    	\noindent{\bf Part 1. Theory}
	
	\vskip 2mm
	
	\noindent{\bf 1.  Classical theory: five great theorems}
	
	1.1 \ Banach-Shauder open mapping theorem
	
	1.2 \ Regular points of smooth maps
	
	1.3 \ Inverse and implicit function theorems
	
	1.4 \ Sard theorem. Transversality.
	
	\vskip 2mm
	
	\noindent{\bf 2. Metric theory. Definitions and equivalences}

	2.1 \ Local regularity

	2.2 \  Non-local regularity
	
	\vskip 2mm
	
	\noindent{\bf 3. Metric theory. Regularity criteria}
		
		3.1\  General criteria
		
		3.2 \ An application: density theorem
		
		3.3 \ Infinitesimal criteria
		
		3.4 \ Related concepts: metric subregularity, calmness, controllability, linear
		recession .
	
   	\vskip 2mm
   	
   	\noindent{\bf 4. Metric theory. Perturbations and stability}
   	
   		4.1 \ Stability under Lipschitz perturbation
   		
   		4.2 \ Strong regularity and metric implicit function theorem
   		
   	\vskip 2mm
   	
   	\noindent{\bf 5.  Banach space theory}
   	
   		5.1 \  Techniques of variational analysis in Banach spaces
   		
   		\qquad 5.1.1 \ Homogeneous set-valued mappings
   		
   	\qquad	5.1.2 \ Tangent cones and contingent derivatives
   	
   		\qquad 5.1.3 \ Subdfferentials, normal cones and coderivatives
   		
   		5.2 \  Separable reduction
   		
   		5.3 \ Contingent derivatives and primal regularity estimates
   		
   		5.4 \ Dual regularity estimates
   		
   		\qquad 5.4.1 \ Neighborhood estimates
   		
   		\qquad 5.4.2 \ Perfect regularity and linear perturbations
   		
   			\vskip 2mm
   			
   			\noindent{\bf 6. Finite dimensional theory}
   			
   			
   			7.1 \ Regularity
   			
   			7.2 \ Subregularity and error bounds
   			
   			7.3 \ Transversality

   		\vskip 4mm
   		
   		\noindent{\bf Part 2. Applications}
	
	\vskip 2mm
	
	\noindent{\bf 7. Special classes of mappings}
	
		7.1 \ Error bounds
		
		\qquad 7.1.1 \ Error bounds for convex functions
		
		\qquad 7.1.2 \ Some general results on global error bounds
		
		7.2  \ Mappings with convex graphs
		
		\qquad 7.2.1 \ Convex processes
		
		\qquad 7.2.2 \ Theorem of Robinson-Ursescu
		
		\qquad 7.2.3 \ Mappings with convex graphs. Regularity rates
		
		7.3 \ Single-valued Lipschitz maps
		
		7.4 \ Polyhedral and semi-linear sets and mapping
		
		7.5 \ Semialgebraic mappings, stratifications and the Sard theorem
		
		\qquad 7.5.1 \ Basic properties
		
		\qquad 7.5.2 \ Transversality

	\vskip 2mm
	
	\noindent{\bf 8. Some applications to analysis and optimization}
	
		8.1 \ Subdfferential calculus
		
		8.2 \ Necessary conditions in constrained optimization
		
		\qquad 8.2.1 \ Noncovering principle
		
		\qquad 8.2.2 \ Exact penalty
		
		\qquad 8.2.3 \ Optimality alternative
		
		\qquad 8.2.4 \ Optimal control of differential inclusions
		
		\qquad 8.2.5 \ Constraint qualifications
		
		8.3 \ An abstract relaxed optimal control problem
		
		8.4 \ Genericity in tame optimization
		
		8.5 \ Method of alternating projection
		
		8.6 \ Generalized equations
		
		8.7 \ Variational inequalities over polyhedral sets
		
		8.8 \ Differential inclusions (existence of solutions)

	\section*{Introduction}
	
	Metric regularity has emerged during last 2-3 decades as one of the central
	concepts of a young discipline now
	often called {\it variational analysis}. The roots of this concept go
	back to a circle of
	fundamental regularity ideas of classical analysis embodied in such results as
	the implicit
	function theorem, Banach open mapping theorem, theorems of Lyusternik and
	Graves, on
	the one hand, and the Sard theorem and the Thom-Smale transversality theory, on the other.
	
	Smoothness is the key property of the objects to which the classical results are
	applied.
	Variational analysis, on the other hand, appeals to objects that may lack this
	property:
	functions and maps that are non-differentiable at points of interest, set-valued
	mappings
	etc.. Such phenomena naturally appear in optimization theory and not only
	there\footnote{Grothendick mentions "ubiquity of stratifed structures in
		practically all domains of geometry" in his1984 {\it Esquisse d'un Programme},
		see \cite{AG}}.

	In the traditional nonlinear analysis, regularity of a mapping (e.g. from a
	normed space
	or a manifold to another) at a certain point means that its derivative at the
	point is onto (the target space or the tangent space of the target manifold).
	This property, translated through available analytic or topological means to
	corresponding
	local properties of the mapping, plays a crucial role in studying some basic
	problems
	of analysis such as existence and behavior of solutions of a nonlinear equation
	$F(x)=y$
	(with $F$ and $y$ viewed as data and $x$ as unknown) under small perturbations
	of the data.
	Similar problems appear if, instead of equation, we consider inclusion
	\begin{equation}\label{pr1}
	y\in F(x)
	\end{equation}
	(with $F$ a set-valued mapping this time) which, in essence, is the main object
	to study
	in variational analysis. The challenge here is evident: no clear way to
	approximate the
	mapping by simple objects like linear operators in the classical case.

	The key step in the answer to the challenge was connected with the understanding
	of the metric nature of some basic phenomena that appear in the classical
	theory. This
	eventually led to the choice of the class of metric spaces as the main
	playground and subsequently
	to abandoning approximation as the primary tool of analysis in favor of a direct
	study of the phenomena as such. The "metric theory" offers a rich collection of
	results
	that, being fairly general and stated in purely metric language, are nonetheless
	easily adaptable to
	Banach and finite dimensional settings (still among the most important in
	applications) and to
	various classes of mappings with special structure. Moreover, however surprising
	this may
	sound, the techniques coming from the metric theory  sometimes appear more efficient,
	flexible and easy to use than the available Banach space techniques (associated with
	subdifferentials and coderivatives, especially in infinite dimensional Banach
	spaces). We shall not once see that proper use of metric criteria may lead to dramatic simplification of proofs and clarification of the ideas behind them.
	This occurs at all  levels of generality, from  results valid in arbitrary metric spaces to specific facts about even fairly simple classes of finite dimensional mappings.
	
	It should be added furthermore that the central role played by distance
	estimates has determined a quantitative character of the theory (contrary to the
	predominantly
	qualitative character of the classical theory). Altogether, this opens gates to
	a
	number of new applications, such as say metric fixed point theory, differential
	inclusions,
	all chapters of optimization theory, numerical methods.
	
	This paper has appeared as a result of two short courses I gave in the
	University of
	Newcastle and the University of Chile in 2013-2014.
	The goal was to give a brief  account of some major principles of the theory of
	metric regularity along with the impression of how they work in various areas of
	analysis and optimization.  The three principal
	themes that will be in the focus of  attention are:
	
	(a) regularity criteria (containing quantitative estimates for rates of
	regularity)
	including formal comparisons of their relative power and precision;
	
	(b) stability problems relating to
	the effect of perturbations of the mapping on its regularity properties, on the
	one hand,
	and  to solutions of  equations, inclusions etc. on the other;
	
	(c) role of metric regularity in analysis and optimization.

	The existing regularity theory of variational analysis may look very technical.
	Many available proofs take a lot of space and use heavy
	techniques.
	But the ideas behind most  basic results,
	especially in the metric theory, are rather simple and in many cases
	proper application of the ideas leads to noticeable (occasionally even dramatic)
	simplification and clarification of the proofs.
	This is a survey paper, so many results are quoted and discussed,
	often  without proofs. As a rule, a proof is given if
	(a) the result is of a primary importance and the proof is sufficiently simple,
	(b) the result is new, (c) the access to the original publication containing the result is not very easy and especially (d) the proof is
	simpler (shorter, or looking more transparent) than available in the literature
	known to me.
	
	And of course there are topics (some important) not touched upon in the paper,
	especially those that can be found in monographic literature. I mean first of
	all
	the books by Dontchev and Rockafellar \cite{DR} and Klatte and Kummer
	\cite{KK02}
	in which metric regularity, in particular its finite dimensional chapter, is
	prominently	presented. Among more specialized topics not touched upon in the survey,
	I would mention nonlinear regularity models, point subdifferential regularity criteria  with associated compactness properties of subdifferentials
	and directional regularity.
	
	The survey consists of two parts. The first part called `Theory' contains an account of the	basic ideas and principles
	 of the metric regularity theory, first in traditional settings of the classical
	 analysis and then for arbitrary set-valued mappings between various classes of spaces.
	In  the second part `Applications' we show how the theory works for some
	specific classes of maps that typically  appear in variational analysis and
	and for a variety of fundamental existence, stability and optimization  problems. 	
	 In preparing this part of the survey the main efforts were focused on finding
	a productive balance between general principles and specific results and/or
	methods associated with the problem. This declaration may look as a sort
	of truism but the point is that publications in which  over-attachment to  certain particular	techniques of variational analysis  (e.g. associated with generalized differentiation) leads to long and poorly
	digestible proofs  of sufficiently simple and otherwise easily provable results
	is not an exceptional phenomenon.

	To conclude the introduction I wish to express my thanks to J. Borwein and A.
	Joffre
	for inviting me to give the lectures that were the basis for this paper and
	to J. Borwein especially for his suggestion to write the survey. I also wish to thank D. Drusvyatskij and A. Lewis for the years of cooperation and many fruitful discussions and to A. Kruger and D. Klatte for many helpful remarks.
	
	\vskip 1mm
	
	\noindent{\bf Dedication}. 2015 and late 2014 have witnessed remarkable jubilees of six my good old friends. I dedicate this paper, with gratitude for the past and warm wishes for the future to
	
	\vskip 2mm
	
	\begin{tabular}{ll}
	Prof. Vladimir Lin\qquad\qquad\qquad& Prof. Terry Rockafellar\\ Prof. Louis Nirenberg & Prof. Vladimir Tikhomirov\\Prof. Boris Polyak& Prof.  Nikita Vvedenskaya	
	\end{tabular}\

	\vskip 5mm

	\noindent{\bf Notation}.
	
	$d(x,Q)$ -- distance from $x$  to $Q$;
	
	$d(Q,P)=\inf \{ \|x-u\|:\; x\in Q,\; u\in P\}$ -- distance between $Q$ and $P$;
	
	${\rm ex}(Q,P)=\sup\{d(x,P):\; x\in Q)\}$ -- excess of $Q$ over P;
	
	$h(Q,P)=\max\{{\rm ex}(Q,P),{\rm ex}(P,Q)\}$ -- Hausdorff distance between $Q$ and $P$;
	
	$B(x,r)$ --  closed ball of radius $r$ and center at $x$;
	
	$\bo (x,r)$ -- open ball of radius $r$ and center at $x$;

	$F|_Q$ -- the restriction of a mapping $F$ to the set $Q$;
	
	$F: X\rra Y$ -- set-valued mapping;
	
	$\gr F=\{(x,y): \; y\in F(x)\}$ -- graph of $F$;

	$I$ -- the identity mapping (subscript, if present, indicates the space, e.g.
	$I_X$);
	
	$\epi f=\{ (x,\al):\; \al\ge f(x)\}$ --  epigraph of $f$;
	
	$\dom f= \{x:\; f(x)<\infty\}$ -- domain of $f$;
	
	$i_Q(x)$ -- indicator of $Q$ (function equal to $0$ on $Q$ and $+\infty$
	outside);
	
	$[f\le \al]=\{ x:\; f(x)\le \al\}$ etc.;
	
	$X\times Y$ -- Cartesian product of spaces;

	$X^*$ -- adjoint of $X$;
	
	$\lan x^*,x\ran$ -- the value of $x^*$ on $x$ (canonical bilinear form on
	$X^*\times X$);
	
	$\R^n$ -- the $n$-dimensional Euclidean space;
	
	$B$ -- the closed unit ball in a Banach space (sometimes  indicated by
	a subscript,
	e.g. $B_X$ is the unit ball in $X$);
	
	$S_X$  --  the unit sphere in $X$;

	$\ker A$ -- kernel of the (linear) operator $A$;
	
	$L^{\perp}=\{x^*\in X^*:\; \lan x^*,x\ran=0,\;\forall\; x\in L\}$ -- annihilator
	of a subspace $L\subset X$;
	
	$K^{\circ}=\{ x^*\in X^*:\; \lan x^*,x\ran\le 0,\;\forall\; x\in K\}$ -- the polar of a cone $K\subset X$
	
	$\im A$ -- image of the operator $A$;
	
	$\cs(X)$ -- collection of closed separable subspaces of $X$;

	$\cll(X,Y)$ -- the space of linear bounded operators $X\to Y$ with
	the {\it operator norm}:
	$$
	\| A\|=\sup_{\| x\|\le 1}\| Ax\|.
	$$
	
	$L\oplus M$ -- direct sum of subspaces;
	
	$T_xM$, $N_xM$  -- tangent and normal space to a manifold $M$ at $x\in M$;
	
	$T(Q,x)$ -- contingent cone to a set $Q$ at $x\in Q$;
	
	$N(Q,x)$ -- normal cone to $Q$ at $x\in Q$, often with a subscript (e.g $N_F$ is
	a Fr\'echet normal cone etc.)

	\vskip 1mm
	
	\noindent We use the standard conventions \
	$d(x,\emptyset)=\infty;\; \inf\emptyset=\infty;\; \sup\emptyset=-\infty$ with one exception: when we deal with non-negative quantities we set $\sup \emptyset=0$.

	\vskip 1cm
	
	\centerline{\bf \large Part 1. Theory}
	
	

	\section{Classical theory: five great theorems.}
	In this section all spaces are Banach.
	\subsection{Banach-Shauder open mapping theorem}
	
	\begin{theorem}[\cite{SB,JS30}]\label{bansha} Let $A:\ X\to
		Y$ be a linear bounded operator onto $Y$, that is $A(X)=Y$. Then
		$0\in\intr A(B)$.
	\end{theorem}
	
	The theorem means that there is a $K>0$ such that for any $y\in Y$
	there is an $x\in X$ such that $A(x)=y$ and $\| x\|\le K\| y\|$
	(take as $K$ the reciprocal of the radius of a ball in $Y$
	contained in the image of the unit ball in $X$ under $A$).

	\begin{definition}[Banach constant]  {\rm Let $A: X\to Y$ be a bounded linear
			operator.
			The quantity
			$$
			C(A) = \sup\{ r\ge 0:\; rB_Y\subset A(B_X)\} = \inf\{\| y\|:\; y\not\in
			A(B_X)\}
			$$
			will be called the {\it Banach constant} of $A$}.
	\end{definition}
	
	The following simple proposition offers two more expressions for the Banach
	constant. Given a linear operator $A: X\to Y$, we set
	$$
	\| A^{-1}\|=\sup_{\| y\|\le 1}d(0,A^{-1}(y))=\sup_{\| y\|=1}\inf\{\| x\|:\;
	Ax=y\}.
	$$
	Of course, if $A$ is a linear homeomorphism, this coincides with the usual norm
	of the inverse operator.
	\begin{proposition}[calculation of $C(A)$]\label{calca}
		For a bounded linear operator $A: X\to Y$
		$$
		C(A) = \inf_{\| y^*\|=1}\| A^* y^*\| = \| A^{-1}\|^{-1}.
		$$
	\end{proposition}

	\subsection{Regular points of smooth maps. Theorems of
		Lyusternik and Graves.}
	
	Let $F: X\to Y$ be Fr\'echet differentiable at $\xb\in X$. It is said that
	$F$ is {\it regular} at $\xb$ if its derivative $F'(\xb)$ is a linear operator
	{\it onto} $Y$. Let $M\subset X$ be a smooth manifold.   The {\it tangent space}   $T_xM$ to $M$ at
	$x\in M$ is the collection of $h\in X$ such that $d(x+th,S)= o(t)$
	when $t\to+0$.
	
	\begin{theorem}[Lyusternik \cite{LAL}] Suppose that $F$ is continuously
		differentiable and regular at $\xb$. Then the tangent space to the
		level set $M=\{ x: \; F(x)=F(\xb)\}$ at $\xb$ coincides with $\ker
		F'(\xb)$.
	\end{theorem}

	\begin{theorem}[Graves \cite{LMG50}]\label{graves} Let $F$ be a continuous
		mapping from a
		neighborhood of $\xb\in X$ into $Y$. Suppose that there are a
		linear bounded operator $A:\ X\to Y$   and positive numbers
	    $\del>0$, $\ga>0$, $\ep>0$  such that $C(A)>\del+\ga$ and
		$$
		\| F(x')-F(x)-A(x'-x)\|< \del \| x'-x\|,
		$$
		whenever $x$ and $x'$ belong to the open $\ep$-ball around $\xb$. Then
		$$
		B(F(\xb),\ga t)\subset F(B(\xb,t))
		$$
		for all $t\in (0,\ep)$.
	\end{theorem}
	
	Here is a slight modification (quantities explicitly added) of the original proof by
	Graves.
	\proof  We may harmlessly assume that $F(\xb)=0$.
	Take $K>0$ such that $KC(A) > 1>K(\del+\ga)$, and let $\| y\|<\ga t$ for some
	$t<\ep$.
	Set $x_0=\xb$, $y_0=y$ and define recursively $x_n$, $y_n$ as follows:
	$$
	y_{n-1}=A(x_n-x_{n-1}),\; \| x_n-x_{n-1}\|\le K\| y_{n-1}\|;\quad
	y_n=A(x_n-x_{n-1}) -(F(x_n)-F(x_{n-1})).
	$$
	It is an easy matter to verify that
	$$
	\| x_n-x_{n-1}\|\le (K\del)^{n-1}K\| y\|,\quad \| y_n\|\le (K\del)^{n}\| y\|
	$$
	and $y_{n-1}-y_n= F(x_n)-F(x_{n-1})$, so that $(x_n)$ converges to some $x$ such
	that $ F(x)=y$ and
	$$
	\| x-\xb\|\le \frac{K}{1-K\del}\| y\|\le \ga^{-1}\|y\|<t
	$$
	as claimed.\endproof
	
	The theorem of Lyusternik was proved in 1934 and the theorem of
	Graves in 1950. Graves was apparently unaware of Lyusternik's result
	and Lyusternik, in turn, of the open mapping theorem by Banach-Shauder.
	Nonetheless
	the methods they used in their proves were very
	similar. For that reason the following statement which is somewhat
	weaker than the theorem of Graves and somewhat  stronger than the
	theorem of Lyusternik is usually called the Lyusternik-Graves
	theorem.
	
	\begin{theorem}[Lyusternik-Graves theorem] Assume that $F:
		X\to Y$ is continuously differentiable and regular at $\xb$. Then
		for any positive $r< C(F'(\xb))$, there is an $\ep>0$ such that
		$$
		B(F(\xb),r t)\subset F(B(\xb,t)),
		$$
		whenever $\| x-\xb\|<\ep,\; 0\le t<\ep$.
	\end{theorem}
	
	It should be also emphasized that no differentiability assumption is made in the
	theorem of Graves. In this respect Graves was much ahead of time. Observe that
	the mapping $F$ in the theorem of Graves can be viewed as a perturbation of $A$
	by a $\del$-Lipschitz mapping. With this interpretation the theorem of Graves
	can be also viewed as a direct predecessor of
	Milyutin's perturbation theorem (Theorem \ref{milt1} in the fourth section),
	which is one of the central results in the regularity theory of variational
	analysis.

	\subsection{Inverse and implicit function theorem}
	
	\begin{theorem}[Inverse function theorem]\label{lugr} Suppose
		that $F$ is continuously differentiable at $\xb$ and the derivative
		$F'(\xb)$ is an invertible operator onto $Y$. Then there is a
		mapping $G$ into $X$ defined in a neighborhood of $\yb=F(\xb)$,
		strictly differentiable at $\yb$ and such that
		$$
		G'(\yb)= \big(F'(\xb)\big)^{-1}\quad {\rm and}\quad F\circ G= I_Y
		$$
		in the neighborhood.
	\end{theorem}
	
	The shortest among standard proofs of the theorem is based on the contraction
	mapping principle (see e.g. the second proof of the theorem in \cite{DR}). But
	equally short proof follows from the theorem of Lyusternik-Graves.
	
	\proof Set $A=F'(\xb)$. Then
	$F(x')-F(x)- A(x'-x)= r(x',x)\| x'-x\|$,
	where $\| r(x',x)\|\to 0$ when $x,x'\to \xb$. As $A$ is
	invertible, there is a $K>0$ such that $\| Ah\|\ge K\| h\|$.
	Hence $\| F(x')-F(x)\|\ge (K-r(x,x'))\| x'-x\|>0$ if $x,\  x'$ are close to
	$\xb$.
	This means that $F$ is one-to-one in a neighborhood of $\xb$.
	But by the Lyusternik-Graves theorem,
	$F(U)$ covers a certain open neighborhood  of $\yb$. Hence $G=F^{-1}$
	is defined in a neighborhood of $F(\xb)$.
	So given $y$ and $y'$ close to $\yb=F(\xb)$
	and let $x', \ x $ be such that $F(x')=y',\ F(x)=y$. Then as
	we have seen $\| y-y'\|\ge K\| x-x'\|$.  We have
	
	$$
	A^{-1}\big(F(x')-F(x)-A(x'-x)\big)=A^{-1}(y'-y) - G(y')-G(y),
	$$
	so  that
	$$
	\begin{array}{lcl}
	\|G(y')-G(y) - A^{-1}(y'-y)\|&\le&\| A\|^{-1}\|F(x')-F(x)- A(x'-x)\|\\
	&=&  \| A^{-1}\|\| r(x',x)\|\| x'-x\|\| \le q(y,y')\| y'-y\|,
	\end{array}
	$$
	where $q(y,y')=Kr(G(y),G(y'))$ obviously goes to zero when $y,y'\to
	\yb$.
	\endproof

	\begin{theorem}[implicit function theorem]\label{implic} Let $X,\ Y,\ Z$ be
		Banach spaces, and let $F$ be a mapping into $Z$ which is defined
		in a neighborhood of $(\xb,\yb)\in X\times Y$ and strictly
		differentiable at $(\xb,\yb)$. Suppose further that the partial
		derivative $F_y(\xb,\yb)$   is an invertible operator. Then there
		are neighborhoods $U\subset X$ \ of $\xb$ and $W\subset Z$ of \
		$\zb=F(\xb,\yb)$ and a mapping $S:\; U\times W\to Y$ such that
		$(x,z)\mapsto (x,S(x,z))$ is a homeomorphism of $U\times W$ onto a neighborhood
		of $(\xb,\yb)$ in $X\times Y$ and
		$$
		F(x,S(x,z))=z,\quad \forall  \ x\in U,\; \forall \ z\in W
		$$
		The mapping $S$ is strictly differentiable at $\xzb$ with
		\begin{equation}\label{1.3.1}
		S_z\xzb=\big(F_y(\xb,\yb)\big)^{-1},\quad
		S_x\xzb=\big(F_y(\xb,\yb)\big)^{-1}F_x(\xb,\yb).
		\end{equation}
	\end{theorem}
	
	The simplest proof of the theorem is obtained by application  of the inverse
	mapping theorem to the following map $X\times Y\to X\times Z$ (see e.g.
	\cite{DR}):
	$$
	\Phi (x,y) = \left(\begin{array}{c}x \\ F(x,y)\end{array}\right).
	$$

	\subsection{Sard theorem. Transversality.}
	
	\begin{definition}[critical and regular value]\label{defcrsm}
		{\rm Let $X$ and $Y$ be
			Banach spaces, and let $F$ be a mapping into $Y$ defined and
			continuously differentiable on an open set of $U\subset X$.   A
			vector $y\in Y$ is called a {\it critical value} of $F$ if there
			is an $x\in U$ such that $F(x)=y$ and $x$ is a singular point of
			$F$. Any point in the range space which is not a critical value is
			called a {\it regular value}, even if it does not belong to $\im
			F$. Thus $y$ is a regular value if either $y\neq F(x)$ for any $x$
			of the domain of $F$ or $\im F'(x)=Y$ for every $x$ such that	 $F(x)=y$.}
	\end{definition}

	\begin{theorem}[Sard \cite{Sard42}]\label{sard} Let $\Omega$ be an open set
		in $\R^n$ and $F$ a $C^k$-mapping from $\Omega$ into $\R^m$.  Then the Lebesgue
		measure of the set of critical values of $F$ is equal to zero,
		provided $k\ge n-m+1$.
	\end{theorem}
	\noindent For a proof of a "full" Sard theorem see \cite{AR}; a much shorter
	proof for $C^{\infty}$ functions can be found in \cite{LN}.
	
	\begin{definition}[transversality]\label{deftran}
		{\rm
			Let $F:X\to Y$ be a $C^1$-mapping,
			and let $M\subset Y$ be a $C^1$-submanifold. Let finally $x$ be in
			the domain of $F$. We say that $F$ is {\it transversal to $M$ at}
			$x$ if either $y=F(x)\not\in M$ or $y\in M$ and $\im F'(x)+T_yM=Y$.
			It is said that $F$ is {\it transversal to} $M$:\ $F\pitchfork M$,
			if it is transversal to $M$ at every $x$ of the domain of $F$.}
	\end{definition}
	
	We can also speak about transversality of two manifolds $M_1$ in $M_2$ in $X$:
	$M_1\pitchfork M_2$ at $x\in M_1\cap M_2$ if $T_xM_1+T_xM_2=X$. For our future
	discussions, it is useful to have in mind that the latter property can be
	equivalently expressed in dual terms: $N_x M_1\cap N_xM_2=\{ 0\}$, where
	$N_xM\subset X^*$ is the {\it normal space} to $M$ at $x$, that is the
	annihilator of $T_xM$.
	
	A connection with regularity is immediate from the definition: if
	$(L,\vf)$ is a local parametrization for $M$
	at $y$ and $y=F(x)$, then
	transversality of $F$ to $M$ at $x$ is equivalent to regularity at
	$(x,0,0)$ of the mapping $\Phi: X\times L\to Y$ given by $\Phi(u,v)=
	F(u)-\vf (v)$.
	
	The connection of  transversality and regularity is actually  much
	deeper. Let $P$ be also a Banach space  and let $F: X\times P\to
	Y$. We can view $F$ as a family of mappings from $X$ into $Y$
	parameterized by elements of $P$. Let us denote ``individual''
	mappings $x\to F(x,p)$ by $F(\cdot,p)$. Let further $M\subset Y$
	be a submanifold, and let $\pi: X\times P\to P$ be the standard
	Cartesian projection $(x,p)\to p$.
	
	\begin{proposition}\label{prethom}Suppose $F$ is transversal to $M$ and
		$Q=F^{-1}(M)$ is a manifold. Let finally $\pi|_Q$ stands for the restriction
		of $\pi$ to $Q$. Then $F(\cdot,p)$ is transversal to
		$M$, provided $p$ is a regular value of $\pi|_Q$.
	\end{proposition}

	Combining the  proposition  with the Sard theorem, we get the following (simple
	version of)  transversality theorem of Thom
	
	\begin{theorem}[see e.g. \cite{GP}]\label{thomsm} Let $X$, $Y$ and $P$ be
		finite dimensional Banach spaces
		Let $M\subset Y$ be a $C^r$-manifold, and let $F: X\times P\to Y$ be a
		$C^k$-mapping $(k\le r)$. Assume that $F\pitchfork M$ and $k> \dim
		X-{\rm codim} M$. Then
		$F(\cdot,p)\pitchfork M$ for each $p\in P$ outside of a subset of $P$
		with $\dim P$-Lebesgue measure zero.
	\end{theorem}

	\section{Metric theory. Definitions and equivalences.}
	
	Here $X$ and $Y$ are metric space. We use the same notation for the metrics in both
	hoping this would not lead to any difficulties.
	
	\subsection{Local regularity}
	We start with the simplest and the most popular case of local regularity near
	a
	certain point of the graph. So let an $F: X\rra Y$ be given as well as a
	$\xyb\in\gr F$.
	\begin{definition}[local regularity properties]\label{defloc}{\rm
			We say that $F$ is
			
			\vskip 1mm
			
			$\bullet$ {\it open} or {\it covering at a linear rate near} $\xyb$ if there
			are $r>0$, $\ep >0$ such that
			$$
			B(y,rt)\cap B(\yb,\ep)\subset F(B(x,t)), \quad \forall\; (x,y)\in\gr F,\;
			d(x,\xb)<\ep, \; t\ge 0.
			$$
			The upper bound $\sur F(\xb|\yb)$ of such $r$ is the {\it modulus} or {\it
				rate of surjection} of
			$F$ near $\xyb$. If no such $r$,  $\ep$ exist, we set $\sur F(\xb|\yb)=0$;
			
			\vskip 1mm
			
			$\bullet$ {\it metrically regular} near $\xyb\in\gr F$ if there are $K>0$,
			$\ep >0$ such that
			$$
			d(x,F^{-1}(y))\le Kd(y,F(x)),\quad {\rm if}\; d(x,\xb)<\ep,\; d(y,\yb)<\ep.
			$$
			The lower bound $\reg F(\xb|\yb)$ of such $K$ is the {\it modulus} or {\it
				rate of metric regularity}
			of $F$ near $\xyb$. If no such $K$, $\ep$ exist, we set $\reg
			F(\xb|\yb)=\infty$.
			
			\vskip 1mm
			
			$\bullet$ {\it pseudo-Lipschitz} or has the {\it Aubin property} near $\xyb$
			if there are $K>0$ and $\ep>0$
			such that
			$$
			d(y,F(x))\le Kd(x,u), \quad {\rm if}\; d(x,\xb)<\ep,\; d(y,\yb)<\ep,\; y\in
			F(u).
			$$
			The lower bound $\lip F(\xb|\yb)$ is the {\it Lipschitz modulus} or {\it
				rate} of $F$ near
			$\xyb$. If no such $K$, $\ep$ exist, we set $\lip F(\xb|\yb)=\infty$.
		}
	\end{definition}
	
	Note a difference between the covering property and the conclusions of
	theorems of Lyusternik and Graves: the
	theorems deal only with the given argument $\xb$ while in the definition we
	speak about all $x\in\dom F$ close to $\xb$. This difference that
	was once a subject
	of  heated discussions is in fact illusory as under the assumptions of the
	theorems of Lyusternik and Graves
	the covering property in the sense of the just introduced definitions is
	automatically  satisfied.

	The key and truly remarkable fact for the theory is that the three parts of the definition actually	speak about the same phenomenon.
	Namely the following holds true unconditionally for any set-valued mapping
	between two metric spaces.
	
	\begin{proposition}[local equivalence]\label{loceq}
		$F$ is open at a linear rate near $\xyb\in\gr F$ if and only if it is
		metrically regular
		near $\xyb$ and if and only if $F^{-1}$ has the Aubin property near
		$(\yb,\xb)$. Moreover,
		under the convention that $0\cdot\infty=1$,
		$$
		\sur F(\xb|\yb)\cdot\reg F(\xb|\yb)= 1;\quad \reg F(\xb|\yb)=\lip
		F^{-1}(\yb|\xb).
		$$
	\end{proposition}
	
	\begin{remark}
		{\rm In view of the proposition it  makes sense to use the word {\it regular}
			to characterize the three properties. This terminology would also emphasize the
			ties with the classical regularity concept. We observe further that while the rates of
regularity are connected with specific distances in $X$ and $Y$, the very fact that $F$ is regular near  certain point is independent of the choice of specific metrics. Thus, although the definitions
explicitly use metrics} the regularity is a topological property.
	\end{remark}
	
	The proof of the proposition is fairly simple (we shall get it as a consequence of
	a more general equivalence theorem later in this section). But the way to it
	was surprisingly long (see brief bibliographic comments at the end of the
	section).
	
	There are other equivalent formulations of the properties.
	For instance, the definition of linear openness/ covering  can be modified by
	adding the constraint $0\le t<\ep$ (see \cite{AI11a}); a well known modification of the
	definition of metric regularity includes the condition that $d(y,F(x))<\ep$.
	The only difference is that the $\ep$'s in the original and modified
	definitions may be different.

	\begin{definition}[graph regularity \cite{LT94}] {\rm $F$ is said to be {\it
				graph-regular at (or near) $\xyb\in\gr F$} if there are
			$K>0,\ \ep >0$ such that the inequality
			\begin{equation}\label{2.2.1}
			d(x,F^{-1}(y))\le Kd((x,y),\gr F),
			\end{equation}
			holds, provided $d(x,\xb)<\ep,\; d(y,\yb)<\ep$}.
	\end{definition}

	\begin{proposition}[metric regularity vs graph regularity
		\cite{LT94}]\label{graphreg}
		
		Let $F: X\rra Y$,  and
		$\xyb\in\gr F$. Then $F$ is metrically regular at $\xyb$
		if and only if it is graph-regular at $\xyb$.
	\end{proposition}
	
	Note that, unlike the equivalence theorem, the last proposition is purely
	local:
	the straightforward non-local extension of this result (e.g. along the lines
	of the
	subsection below) is wrong.

	\subsection{Non-local regularity.}
	
	As we have already mentioned, most of current researches focus on local
	regularity.\
	(although the first abstract definition of the covering property given in
	\cite{DMO}
	was absolutely non-local).  To a large extent this is because of the close
	connection of modern variational analysis studies with optimization theory which
	is basically interested in local results: optimality conditions, stability of
	solutions under small perturbations, etc. Another less visible reason is that
	non-local regularity is a more delicate concept: in the non-local case we cannot
	freely change the regularity domain
	that is an integral part of the definition. Meanwhile non-local regularity is,
	a powerful instrument for proving e.g. various existence
	theorems (see e.g.  subsection 8.7).
	
	Let  $U\subset X$ and $V\subset Y$ (we usually assume $U$ and $V$ open), let
	$F: X\rra Y$, and
	let $\ga(\cdot)$ and $\del(\cdot)$ be  extended-real-valued  functions on
	$X$ and $Y$ assuming  positive values (possibly infinite) respectively on $U$
	and $V$.
	\begin{definition}[non-local regularity properties
		\cite{AI11a}]\label{nonlocdef} {\rm
			We say that $F$ is
			
			$\bullet$ \ {\it $\ga$-open (or $\ga$-covering) at a linear rate} on  $U\times
			V$ if there is an $r>0$ such that
			$$
			B(F(x),rt)\bigcap V\subset F(B(x,t)),
			$$
			if $x\in U$  and $t<\ga(x)$.
			Denote by $\sur_{\ga}F(U|V)$
			the upper  bound of such $r$. If no such $r$ exists, set
			$\sur_{\ga}F(U|V)=0$.  We shall call
			$\sur_{\ga}F(U|V)$ the {\it modulus} (or {\it rate}) {\it of
				$\ga$-openness} of $F$ on $U\times V$;
			
			$\bullet$ \ {\it  $\ga$-metrically  regular on   $U\times V$} if there is a
			$K>0$ such that
			$$
			d(x,F^{-1}(y))\le Kd(y,F(x)),
			$$
			provided $x\in U,\; y\in V$ and $Kd(y,F(x))<\ga(x)$. Denote by
			$\reg_{\ga}F(U|V)$
			the lower bound of such $K$.
			If no such $K$ exists, set $\reg_{\ga}F=\infty$.  We shall call
			$\reg_{\ga}F(U|V)$ the {\it modulus}
			(or {\it rate}) {\it of $\ga$-metric regularity} of $F$ on $U\times V$;
			
			$\bullet$ \ {\it  $\del$-pseudo-Lipschitz} on  $U\times V$ if  there
			is a $K>0$
			such that
			$$
			d(y,F(x))\le Kd(x,u)
			$$
			if   $x\in U,\; y\in V,\; K d(x,u)< \del (y)$ and $y\in F(u)$.
			Denote by $\Lip_{\del}F(U|V)$   the lower bound of such $K$.
			If no such $K$ exists, set $\Lip_{\del}F =\infty$ . We shall call
			$\Lip_{\del}F(U|V)$
			the {\it $\del$-Lipschitz modulus} of $F$ on $U\times V$.}
	\end{definition}
	
	If $U=X$ and
	$V=Y$, let us agree to write
	$  \sur_{\ga} F,\  \reg_{\ga}F,\ \lip_{\del} F$
	instead of $\sur_{\ga}F(X|Y)$, etc.
	The role of the functions $\ga$ and $\del$ is clear from the definitions. They
	determine how
	far we shall reach from any given point in verification of the defined
	properties.
	It is therefore natural to call them {\it regularity horizon} functions.  Such
	functions are inessential for local regularity (see e.g. Exercise 2.8 below). But for fixed
	$U$ and $V$ regularity horizon function is an essential element of the
	definition. Regularity properties corresponding to different $\ga$
	may not be equivalent (see Example 2.2 in \cite{AI14} and also Exercise 2.8 below).
	
	\begin{theorem}[equivalence theorem]\label{equiv3}
		The following three properties are equivalent for any pair of metric spaces
		$X,\ Y$,
		any $F: X\rra Y$, any  $U\subset X$ and $V\subset Y$ and any
		(extended-real-valued) function $\ga(x)$
		which is positive on $U$:
		
		\vskip 1mm
		
		(a) \ $F$ is $\ga$-open at a linear rate on $U\times V$;
		
		\vskip 1mm
		
		(b) \  $F$ is  $\ga$-metrically regular on $U\times V$;
		
		\vskip 1mm
		
		(c) \  $F^{-1}$ is $\ga$-pseudo-Lipschitz  on $V\times U$.
		
		\vskip 1mm
		
		\noindent   Moreover (under the convention that $0\cdot\infty=1$)
		$$
		\sur_{\ga}F(U|V)\cdot\reg_{\ga}F(U|V)=1,\quad
		\reg_{\ga}F(U|V)=\Lip_{\ga}F^{-1}(V|U).
		$$
	\end{theorem}
	
	\proof
	The implication (b) $\Rightarrow$ (c) is  trivial. Hence $\Lip_{\ga}F^{-1}(V|U)\le
	\reg_{\ga}F(U|V)$.
	To prove that  (c) $\Rightarrow$ (a),  take a $K>\lip_{\ga}F^{-1}$ and
	an $r< K^{-1}$ ,  let $t<\ga(x)$, and let $x\in U$,  \ $y\in V$, \ $v\in F(x)$
	and  $y\in B(v,tr)$.
	Then $d(y,v)<r\ga(x)$ and by (c) $d(x,F^{-1}(y))\le K d(y,v)< r^{-1}d(y,v)\le
	t$. It follows that
	there is a $u$ such that $y\in F(u)$ and $d(x,u)<t$. Hence $y\in F(B(x,t))$.
	It follows that
	$r\le\sur_{\ga}F$, or equivalently $1\le K\sur_{\ga}F$. But $r$ can be chosen
	arbitrarily close to $K^{-1}$ and
	and $K$ can be chosen arbitrarily close to $\Lip_{\ga}F^{-1}$. So we conclude
	that $\sur_{\ga}F\cdot\Lip_{\ga}F^{-1}\ge  1.$
	
	Let finally (a) hold with some $r>0$, let $x\in U,\  y\in V$, and let
	$d(y,F(x))<\ga(x)$. Choose a $v\in F(x)$ such that  $d(y,v)<r\ga(x)$
	and set $t=d(y,v)/r$. By (a) there is a $u\in F^{-1}(y)$ such that $d(x,u)\le
	t$.
	Thus $d(x,F^{-1}(y))\le t=d(y,v)/r$. But $d(y,v)$ can be chosen arbitrarily
	close to
	$d(y,F(x))$ and we get $d(x,F^{-1}(y))\le r^{-1}d(y,F(x))$, that is
	$r\cdot\reg_{\ga}F\le 1$.
	On the other hand $r$ can be chosen arbitrarily close to
	$\sur_{\ga}F$ and we can conclude that $\sur_{\ga}F\cdot\reg_{\ga}F\le 1$ so
	that
	$$
	1\ge \sur_{\ga}F(U|V)\cdot\reg_{\ga}F(U|V)\ge
	\sur_{\ga}F(U|V)\cdot\lip_{\ga}F(V|U)\ge 1
	$$
	which completes the proof of the theorem.     \endproof

	The most important example of the  horizon function is \
	$m(x)=d(x,X\backslash U)$. The meaning is that we need not look at points beyond $U$.
	 We shall call $F$ {\it Milyutin regular} on $U\times
	V$ if it is $m$-regular.
	(This is actually the type of regularity implicit in the definition given in
	\cite{DMO}.)
	In what follows we shall deal only with Milyutin regularity when speaking
	about non-local matters.
	
	\begin{exercise} {\rm Prove that} $F$ is regular near $\xyb\in\gr F$ if and
		only if it is Milyutin regular on $\bo(\xb,\ep)\times \bo (\yb,\ep)$ for all
		sufficiently small $\ep$.
	\end{exercise}
	
	We conclude the section with a useful result (a slight modification of the
	corresponding result in \cite{AI00}) showing that, as far as metric regularity
	is concerned, any set-valued mapping can be equivalently and in a canonical way
	replaces  by a single-valued mapping continuous on its domain.

	\begin{proposition}[single-valued reduction]\label{sinval}  Let $X\times Y$ be
		endowed with the $\xi$-metric. Let
		$F$ be Milyutin regular on $U\times V$ with $\sur_mF(U|V)\ge r>0$.
		Consider the mapping $\Pi_F: \gr F\to Y $ which is
		the restriction to $\gr F$ of the Cartesian projection $(x,y)\to y$.
		
		Then $\Pi_F$ is Milyutin regular on $(U\times Y)\times V$ and
		$\sur _m \Pi_F(U\times Y|V) = \sur_m F(U|V)$
		if $X\times Y$ is considered e.g. with the $\xi$-metric.
	\end{proposition}

	A few bibliographic comments.
	To begin with, it is worth mentioning   that in the classical theory no interest
	to metric estimates can be traced. The covering property close to the covering
	part of Milyutin regularity was introduced in \cite{DMO} and attributed to
	Milyutin.
	An  estimate of metric regularity type
	first time appeared in Lyusternik's paper   \cite{LAL} but for $x$ restricted
	to the kernel of the derivative.
	In Ioffe-Tikhomirov \cite{IT} metric regularity  was proved under the
	assumptions of the Graves theorem. Robinson  was probably the first to consider set-valued
	mappings. In \cite{SMR76}
	he proved metric regularity of the mapping $F(x)= f(x)+K$ (even of the
	restriction of this
	mapping to a convex closed subset of $X$), under the assumptions that $f:X\to Y$ is continuously
	differentiable and $K\subset Y$ is a closed convex cone, under certain
	qualification condition extending Lyusternik's $\im F'(x) = Y$. The definition
	of $\ga$-regularity was given in \cite{AI11a}.
	
	Equivalence of covering and metric
	regularity was explicitly mentioned (without proof) in the paper of
	Dmitruk-Milyutin-Osmolovski  \cite{DMO} that marked the beginning of systematic
	study
	of the regularity phenomena, in particular in metric spaces, and  Ioffe in
	\cite{AI81}
	stated a certain equivalence result (Proposition 11.12 -- see \cite{AI98a}
	for its proof) which, as was much later understood, contains even more precise
	information about the connection of the
	covering and metric regularity properties. And the pseudo-Lipschitz property
	was introduced by Aubin in  \cite{JPA84}.
	
	This was the sequence of events prior to the
	proof of the equivalence of the three properties by Borwein-Zhuang
	\cite{BZ88}
	and Penot \cite{JPP89}. It has to be mentioned that in both papers more
	general
	"nonlinear" properties were considered. In this connection we also mention
	the
	paper by Frankowska \cite{HF90} with a short proof of nonlinear openness and
	some pseudo-H\"older property.

	\section{Metric theory. Regularity criteria.}

	This section is central.
	Here we prove  necessary and sufficient conditions
	for
	regularity. The key results is Theorems \ref{gencrit},  \ref{secmil} and \ref{secmilmod}
containing   general regularity criteria.
	The criteria (especially the first of them)   will serve us as a basis for
	obtaining various qualitative
	and quantitative   characterizations of regularity in this and subsequent
	sections. The criteria are very simple to prove and, at the same time,
	provide us with  an instrument of analysis which is both powerful and easy to
	use.
	We shall see this already in this section and many times in what follows.
	In the second subsection we consider infinitesimal criteria for local regularity
	based on the concept of {\it slope}, the central  in the local theory.
	
	Given a set-valued mapping $F: X\rra Y$, we associate with it the following
	functions that will be systematically used in connection with the criteria and
	their applications:
	$$
	\vf_y(x,v)=\left\{\begin{array}{lc} d(y,v),&{\rm if}\; v\in F(x);\\ +\infty,&
	{\rm otherwise}\end{array}\right.;\quad \psi_y(x) = d(y,F(x));
	\quad \bpsi_y(x) = \liminf_{u\to x}\psi (u).
	$$
	Note that $\vf_y$  is  Lipschitz continuous  on $\gr F$,  hence it is lower
	semicontinuous whenever $\gr F$ is a closed set.
	
	\subsection{General criteria.}
	Given a $\xi>0$, we define the $\xi$-metric on $X\times Y$ by
	$$
	d_{\xi}((x,y),(x',y'))=\max\{d(x,x'),\xi d(y,y')\}.
	$$
	
	\begin{theorem}[criterion for Milyutin regularity]\label{gencrit}
		Let $U\subset X$ and $V\subset Y$ be open sets, and let $F: X\rra Y$ be
		a set-valued mapping whose graph is  complete in the product metric. Let
		further
		$r>0$ and there be a $\xi>0$ such that  for any
		$x\in U$,  $y\in V$, $v\in F(x)$ with  $0<d(y,v)<r m(x)$,
		there is a pair $(u,w)\in\gr F$  different from  $(x,v)$ and such that
		\begin{equation}\label{3.1}
		d(y,w)\le d(y,v)-rd_{\xi}((x,v),(u,w)).
		\end{equation}
		Then  $F$ is  Milyutin regular on $U\times V$   with $\sur_mF(U|V)\ge r$.
		
		Conversely, if $F$ is Milyutin regular on $U\times V$, then for any positive
		$r<\sur_{\ga}F(U|V)$,
		any $\xi\in(0,r^{-1})$,  any  $x\in U$, $v\in F(x)$ and
		$y\in V$ satisfying $0<d(y,v)<r\ga(x)$,
		there is a pair $(u,w)\in\gr F$ different from $(x,v)$ such that (\ref{3.1})
		holds.
	\end{theorem}

	The theorem offers a very simple geometric interpretation of the regularity
	phenomenon: it  means that $F$ is regular  if
	for any $(x,v)\in\gr F$ and any $y\neq v$ there is a point in the graph
	whose $Y$-component is closer to $y$ (than $v$) and the distance from the new
	point
	to the original point $(x,v)$ is proportional to the gain in the distance to
	$y$.
	
	\proof We have to verify that, given $(\xb,\vb)\in \gr F$ with $\xb\in U$,
	$y\in V$ and $0<d(y,\vb)\le rt,\; t<m(\xb)$, there is a $u\in B(x,t)$ such that
	$y\in F(u)$.
	We have $\vf_y(\xb,\vb)\le rt$. By Ekeland's variational principle (see e.g.
	\cite{BZ})
	there is a pair $(\hat x,\hat v)\in\gr F$ such that $d_{\xi}((\hat x,\hat
	v),(\xb,\vb))\le t$
	and
	\begin{equation}\label{3.2}
	\vf_y(x,v)+ rd_{\xi}((x,v),(\hat x,\hat v))> \vf_y(\hat x,\hat v)
	\end{equation}
	if $(x,v)\neq (\hat x,\hat v)$. We claim that $\vf_y(\hat x,\hat v)=0$, that is
	$y=\hat v\in F(\hat x)$.  Indeed, $\hat x\in U$, so by the assumption if
	$y\neq \hat v$, there is a pair $(u,w)\neq (\hat x,\hat v)$
	and such that  (\ref{3.1}) holds with $(\hat x,\hat v)$ as $(x,v)$ which
	however contradicts (\ref{3.2}). This proves the first statement.
	
	Assume now that $F$ is Milyutin regular on $U\times V$ with the surjection
	modulus not smaller than $r$.
	Take a  positive $\xi<r^{-1}$ and  $x\in U$, $y\in V$,  $v\in F(x)$ with
	$d(y,v)<r\ga(x)$.
	Take a small $\ep\in (0,r)$  and choose a $t\in (0,m(x))$  such that
	$(r-\ep)t\le d(y,v)< rt$.
	By regularity there is a $u$ such that $d(u,x)< t$ and $y\in F(u)$.
	Note that $t>\xi d(y,v)$ by the choice of $\xi$. So setting  $w=y$ we have
	$t>\xi d(v,w)$ and
	$$
	d(y,w)=0\le d(y,v)-(r-\ep)t\le d(y,v)- (r-\ep)d_{\xi}((x,v),(u,w)).
	$$
	Since $\ep$ can be chosen arbitrarily small, the result follows.  \endproof

	\begin{theorem}[second criterion for Milyutin regularity]\label{secmil}
		Let $X$ be a complete metric space, $U\subset X$ and $V\subset Y$ open sets
		and
		$F: X\rra Y$  a set-valued mapping with closed graph.
		Then $F$ is Milyutin regular on $U\times V$ with $\sur_mF(U|V)\ge r$
		if and only if for any $x\in U$ and any $y\in V$ with $0<\bpsi_y(x)<rm(x)$
		there is a $u\neq x$
		such that
		\begin{equation}\label{3.3}
		\bpsi_y(u)\le \bpsi_y(x)-r d(x,u).
		\end{equation}
		
	\end{theorem}

	\proof   The proof of sufficiency is similar to the proof of the first part of
	the previous theorem.

	To prove that (\ref{3.3}) is necessary for Milyutin regularity take
	$x\in U,\; y\in V$ such that  $0<d(y,F(x))< rm(x)$.  Take $\rho<r$ such that
	still $d(y,F(x))< \rho m(x)$,
	and let $\rho<\rho'<r$.
	Let $x_n\to x$ be such that $d(y,F(x_n))\to \bpsi_y(x)$.
	We may assume that $d(y,F(x_n))< r m(x)$ for all $n$.
	Choose positive $\del_n\to 0$ such that $d(y,F(x_n))\le (1+\del_n) \bpsi_y(x)$,
	and let $t_n$
	be defined by $\rho' t_n=(1+\del_n)\bpsi_y(x)$.  Then $y\in \bo
	(F(x_n),\rho't_n)$, $t_n< m(x_n)$
	(at least for large $n$)   and due to
	the regularity assumption on  $F$ for any $n$ we can find a
	$u_n$ such that $d(u_n,x_n)<t_n$ and $y\in F(u_n)$. Note that $u_n$ are bounded
	away from $x$
	for otherwise (as $\gr F$ is closed) we would inevitably conclude that $y\in
	F(x)$ which cannot
	happen as $\bpsi_y(x)>0$. This means that $\la_n=d(u_n,x_n)/d(u_n,x)$ converge
	to one. Thus
	$$
	\begin{array}{lcl}
	\bpsi_y(u_n)=0&=& \bpsi_y(x)-\bpsi_y(x)=\bpsi_y(x)-\dfrac{\rho't_n}{1+\del_n}\\
	&\le& \bpsi_y(x)-\dfrac{\rho'}{1+\del_n}d(u_n,x_n)\\  &=&
	\bpsi_y(x)-\dfrac{\la_n\rho'}{1+\del_n}d(u_n,x)
	\le  \bpsi_y(x) -\rho  d(u_n,x),
	\end{array}
	$$
	the last inequality being eventually true as $\la_n\rho'>\rho(1+\del_n)$  for
	large $n$.
	\endproof
	
The theorem is especially convenient when $\psi_y$ is lower semicontinuous for
every $y\in V$. Otherwise, the need for  preliminary calculation of
$\bpsi_y$, the lower closure of $\psi_y$, may cause difficulties. It is possible however to
modify the condition of the theorem and get a statement that requires verification of
(\ref{3.3})-like inequality for $\psi$ rather than $\bpsi$, although at the expense  of
some additional uniformity assumption.

\begin{theorem}[modified second criterion for Milyutin regularity]\label{secmilmod}
Let $X$, $Y$, $F$, $U$ and $V$ be as in Theorem \ref{secmil}. A necessary and sufficient condition
for $F$ to be Milyutin regular on $U\times V$ with $\sur F(\xb|\yb)\ge r$ is that there is a $\la\in (0,1)$ and for any
$x\in  U$ and $y\in V$ with $0<\psi_y(x)< rm(x)$ there is a $u\neq x$ such that
\begin{equation}\label{3.3m}
		\psi_y(u)\le \psi_y(x)-r d(x,u),\quad \ \la \psi_y(u)\le\la \psi_y(x).
		\end{equation}
\end{theorem}

\proof The key for understanding the theorem is the following implication
\begin{equation}\label{impcom}
\bpsi_y(x)=0\ \Rightarrow\ y\in F(x)
\end{equation}
of course valid, under the condition of the theorem for $x\in U,\ y\in V$. Indeed,
$\bpsi_y(x)=0$ means that there is a sequence $(x_n)$ converging to $x$ such that
$\psi_y(x_n)\to 0$. This in turn implies the existence of $v_n\in F(x_n)$ converging to
$y$. As the graph of $F$ is closed, it follows that $(x,y)\in \gr F$ as claimed.

Now we can verify that under  the assumptions of the theorem, the condition
of Theorem \ref{secmil} holds. So let
$x\in U,\ y\in V$ and $0<\al=\bpsi_y(x)$. Take $x_n\to x$  such that
$\psi_y(x_n)=\al_n\to \al$ and for each $n$ a $u_n$ such that $\psi_y(u_n)\le\la \al_n$
and $\psi_y(u_n)\le\psi_y(x_n)-rd(x_n,u_n)$. An easy calculation shows that
$$
\psi_y(u_n)\le\bpsi_y(x)- rd(x,u_n) + \ep_n,
$$
where $\ep_n\to 0$. As $d(x,u_n)$ are bounded away from zero by a positive constant,
we have $\ep_n=\del_nd(x,u_n)$, where $\del_n\to 0$. Combining this with the above inequality,
we conclude that for any $r'<r$ that $u_n\neq x$ and inequality
$$
\bpsi_y(u_n)\le \bpsi_y(x)- r'd(x,u_n)
$$
holds for sufficiently large $n$. This  allows to apply Theorem \ref{secmil} and conclude
(by virtue of (\ref{impcom})) that there is a $w\in B(x,(r')^{-1})$ such that $y\in F(x)$, that is
 $\sur_mF(U|V)\ge r'$. \endproof

Note that the proof of necessity in the two last theorems does  not differ from the proof of Theorem \ref{gencrit}. 	Corresponding criteria for local regularity are  immediate.
	\begin{theorem}[criterion for local regularity]\label{locrit}
		Let $F: X\rra Y$ be a set-valued mapping with closed graph, and let $\xyb\in\gr
		F$.
		Then $F$ is regular
		near $\xyb$ if and only if there are $\ep> 0$, $\xi>0$ and $r>0$ such that
		for any $x, \ v$ and $y$
		satisfying $d(x,\xb)<\ep,\; d(y,\yb)<\ep,\; v\in F(x)$ and
		$0<d(y,v)<\ep$ either of the following two properties is valid:
		
		(a) $\gr F$ is locally complete and there is a pair $(u,w)\in \gr F$,
		$(u,w)\neq (x,v)$ such
		that (\ref{3.1}) holds.
		
		(b) $X$ is a complete metric space, the graph of $F$ is closed  and either (\ref{3.3}) or (\ref{3.3m}) holds true.
		
		\noindent Moreover, in either case  $\sur F(\xb|\yb)\ge r$.
	\end{theorem}

	Theorem \ref{gencrit} is a particular case of the criterion for $\ga$-regularity
	proved in \cite{AI11a}.
	Theorem \ref{locrit} is a modification of the result
	established in \cite{AI00}. Theorem \ref{secmil} is  new
	but it was largely stimulated by a recent result of Ngai-Tron-Thera \cite{NTT12}
	(see Theorem \ref{slsc} later  in this section) and by a much earlier
	observation by  Cominetti \cite{RC90} that $\bpsi_y(x)=0$ implies that $y\in
	F(x)$.
	Surprisingly, it has been recently discovered that
	sufficiency in statement of the part (a) of the local criterion (Theorem \ref{locrit}) is present as a remark in a
	much earlier paper by Fabian and Preiss \cite{FP87}.

	The completeness assumption in the first theorem differs from the corresponding assumption of the other two theorems. So the natural
	question is if and how they are connected. It is an easy matter to see, in view
	of Proposition \ref{sinval}, that Theorem \ref{gencrit}
	follows from Theorem \ref{secmil}.  On the other hand,  Theorem \ref{gencrit} is
	easier to use as it does not need a priori calculation of any limit or verification of the existence of $\la$ as in the third theorem.
	However, if the functions $d(y,F(\cdot))$ are lower semicontinuous, the second
	criterion
	may be more convenient. It should also be observed that
	the theorems can be equivalent in some cases (as follows from Proposition 1.5 in
	\cite{AI00}).

	\subsection{An application: density theorem.}
	Here is the first example demonstrating how handy and powerful the criteria are.
	\begin{theorem}[density theorem \cite{DMO,AI11a}]\label{dens}
		Let $U\subset X$ and $V\subset Y$ be open sets,  let $F: X\rra Y$ be a
		set-valued mapping
		with complete graph.  We assume that
		whenever $x\in U$,  $v\in F(x)$ and $t< m(x)$, the set $F(B(x,t))$ is a $\ell
		t$-net in
		$B(v,rt)\bigcap V$, where $0\le \ell<r$.  
		Then $F$ is Milyutin
		regular on $U\times V$ and $\sur_{m}F\ge r-\ell$ .     In particular, if
		$F(B(x,t))$
		is dense in $B(F(x),rt)\bigcap V$ for $x\in U$ and $t<\ m(x)$, then
		$\sur_{m}F(U|V)\ge r$ .
	\end{theorem}
	
	\proof
	Take $x\in U$ and suppose $y\in  V$  is such that $d(y,F(x))< rm(x)$.
	Take a $v\in F(x)$ such that
	$d(y,v)<rm(x)$ and set $t=d(y,v)/r$.  Then $t<m(x)$ and
	by the assumption we can choose $(u,w)\in\gr F$ such that $d(x,u)\le t$ and
	$d(y,w)\le \ell t= (\ell/r)d(y,v)$. Then
	$$
	d(v,w)\le d(y,v)+d(y,w)\le \big(1+\frac{\ell}{r}\big)d(y,v)\le 2d(y,v).
	$$
	Take a $\xi>0$ such that  $\xi r \le 1/2$. Then $\xi d(v,w)<2\xi r t\le  t$ and
	therefore
	$$
	d(y,w)\le \ell t=rt-(r-\ell)t = d(y,v)-(r-l)t \le
	d(y,v)-(r-\ell)d{\xi}((x,v),d(u,w)).
	$$
	Apply Theorem \ref{gencrit}\endproof
	
	\begin{exercise}\label{denscl}
		{\rm Prove the theorem under the assumptions of Theorem \ref{secmil} rather
			than Theorem \ref{gencrit}. }
	\end{exercise}
	
	\begin{exercise} {\rm Prove Banach-Shauder open mapping theorem using the
			density theorem
			(and the Baire category theorem)}
	\end{exercise}
	
	The specification of Theorem \ref{dens}  for local regularity at $\xyb$ is
	
	\begin{corollary}[density theorem - local version]\label{lodense}
		Suppose there are $r>0$,  and $\ep>0$ such that $F(B(x,t))$ is an
		$\ell t$-net   in $B(v,rt)$ whenever
		$d(x,\xb)<\ep,\; d(v,\yb)<\ep,\; v\in F(x)$ and $t<\ep$. Then $\sur
		F(\xb|\yb)\ge r-\ell$. Thus if
		$B(v,rt)\subset \cl F(B(x,t))$ for all $x,\ v$ and $t$ satisfying the specified
		above conditions,
		then $B(v,rt)\subset F(B(x,t)$ for the same set of the variables.
	\end{corollary}
	
	The density phenomenon
	was extensively discussed, especially  at the early stage of developments.
	Results   in the spirit of Corollary \ref{lodense} were
	first considered in
	Ptak \cite{VP74}, Tziskaridze \cite{KST75} and Dolecki \cite{SD78,SD78a} in
	mid-1970s.
	The very idea
	(and to a large extent the techniques used) could be traced back
	to Banach's proof of the closed graph/open mapping theorem.
	Some of the subsequent studies
	(e.g. \cite{BZ88,CU05})  were primarily concentrated on results of such type.
	We refer to \cite{ACP98}
	for detailed discussions and many references.
	Dmitruk-Milyutin-Osmolovskii in \cite{DMO} made a substantial step forward when
	they replaced
	(in the global context) the density
	requirement by the assumption that   $F(B(x),t)$ is an $\ell t$-net in
	$B(F(x),rt)$. This opened
	way to proving the Milyutin perturbation theorem (see the next section). A
	similar advance in the framework of the infinitesimal
	approach (for mappings between  Banach spaces)
	was made by Aubin \cite{JPA81}.

	\subsection{Infinitesimal criteria.}
	The main tool of the infinitesimal regularity theory in metric spaces is
	provided by the concept of (strong) slope -- which is just  the maximal speed of
	descent of
	the function from a given point -- introduced in
	1980 by DeGiorgi-Marino-Tosques \cite{DMT} and since then widely used in various
	chapters of metric analysis.
	
	\begin{definition}[slope]\label{dslope} {\rm Let $f$ be an extended-real-valued
			function
			on $X$ which is finite at $x$. The quantity
			$$
			|\nabla f|(x) =\limsup_{u\to x\atop{u\neq
					x}}\frac{(f(x)-f(u))^+}{d(x,u)}
			$$
			is called the {\it (strong) slope} of $f$ at $x$. We also agree to
			set $|\nab f|(x)=\infty$ if $f(x)=\infty$.  The function is called
			{\it calm} at $x$ if $|\nabla f|(x)<\infty$.}
	\end{definition}

	We shall consider only local regularity in this subsection (although it is
	possible
	to give slope-based characterizations of Milyutin regularity as well).
	It is easy to observe that $|\nabla f|(x)>r$  means that  arbitrarily close to
	$x$
	there are $u\neq x$ such that  $f(x)> f(u) + rd(x,u)$. This allows to
	reformulate the sufficient part
	of the regularity criteria of  Theorem \ref{locrit}.
	in infinitesimal terms.  To this end set as before
	$$
	\vf_y(x,v)=d(y,v)+i_{\gr F}(x,v),\quad\psi _y(x) = d(y,F(x)),\quad
	\overline{\psi}_y(x)= \liminf_{u\to x}\psi_y(u),
	$$
	and let $\nabla_{\xi}$ stand for the slope
	of functions on $X\times Y$ with respect to the $d_{\xi}$-metric:
	$d_{\xi}((x,v),(x',v'))=\max\{d(x,x'),\xi d(v,v'))$.
	
	Things are more complicated with the necessity part: to prove it, an additional
	assumption  on the target space  is needed. Namely, let us say that
	a metric space $X$ is {\it locally coherent} if for any $x$
	$$
	\lim_{u,w\underset{u\neq w}\to x}|\nabla d(u,\cdot)|(w)=1.
	$$
	It can be shown that a convex set and a smooth manifold in a Banach space are
	locally coherent in the induced metric \cite{AI07a} and that any length metric
	space (space whose metric is defined by minimal lengths of curves connecting
	points) is locally
	coherent (as follows from \cite{AC08}).
	
	\begin{theorem}[local regularity criterion 1 \cite{AI07a}]\label{critgen} Let
		$X$
		and $Y$ be  metric spaces, let $F:\ X\rra Y$ be a set-valued mapping,
		and let $\xyb\in \gr F$.  We assume that $\gr F$ is locally complete at
		$\xyb$.
		Suppose further that there are $\ep> 0$, and $r>0$ such that for some $\xi>0$
		\begin{equation}\label{3.21}
		|\nabla_{\xi}\vf_y|(x,v)>r
		\end{equation}
		if
		\begin{equation}\label{3.22}
		v\in  F(x),\quad d(x,\xb)<\ep,\quad d(y,\yb)<\ep,\quad
		d(v,\yb)<\ep,\quad v\neq y.
		\end{equation}
		Then $F$ is regular near $\xyb$ with $\sur F\xyb\ge r.$

		Conversely, let \ $Y$ \ be  locally coherent  at $\yb$. Assume that
		$\sur F(\xb|\yb)> r>0$. Take a $\xi<r^{-1}$. Then for any
		$\del>0$ there is an $\ep>0$  such that $|\nab_{\xi}
		\vf_y|(x,v)\ge(1-\del)r$ whenever $(x,y,v)$ satisfy (\ref{3.22}). Thus, in
		this case
		\begin{equation}\label{3.23}
		\sur F\xyb=\liminf_{{(x,v)\  \underset{{\rm Graph}F}\to\
				(\xb,\yb)}\atop{y\to\yb,\ y\neq v}}|\nabla_{\xi} \vf_y|(x,v).
		\end{equation}
	\end{theorem}

	For mappings into metrically
	convex  spaces (for any two points there is a shortest path connecting the
	points) the final  statement of Theorem \ref{critgen} can be slightly improved.

	\begin{corollary}\label{critgp} Suppose under the conditions of Theorem
		\ref{critgen} that
		$Y$ is metrically convex. Then for any neighborhood $V$ of $\yb$
		\begin{equation}\label{3.26}
		\sur F\xyb= \liminf_{(x,v)\underset{\rm Graph F}\to
			(\xb,\yb)}~\inf_{y\in V\backslash\{v\}}|\nabla_{\xi} \vf_y|(x,v)
		\end{equation}
	\end{corollary}

	\begin{theorem}[local regularity criterion 2]\label{slsc}
		Suppose that $X$ is complete and the graph of $F$ is closed. Assume further
		that there are neighborhood $U\subset X$ of $\xb$ and $V\subset Y$ of $\yb$,
		$r>0$ and $\ep>0$ such that
		that $|\nabla\bpsi_y|(x)>r$  for all $(x, y)\in U\times V$  such that
		$\ep>\bpsi_y(x)>0$.
		Then $\sur F(\xb|\yb)\ge r$.
		
		Conversely,
		if in addition  $Y$ is a length space and $\sur F(\xb|\yb)>r>0$, then there is a
		neighborhood of $\xyb$ and an  $\ep>0$ such that $|\nabla\bpsi_y|(x)\ge  r$
		for all $(x,y)$ of the neighborhood such that $y\not\in F(x)$ and
		$0<\bpsi_y(x)<\ep r$. Thus in this case
		$$
		\sur F(\xb|\yb)=\liminf_{{(x,y)\to \xyb}\atop{0\neq d(y,F(x))\to
				0,}}|\nabla\bpsi_y|(x).
		$$
		
		In particular,  if $\psi_y=d(y,F(\cdot))$ is lower semicontinuous at every $x$
		of a neighborhood of $\xb$ and for every
		$y\not\in F(x)$ close to $\yb$,  then
		$$
		\sur F(\xb|\yb)=\liminf_{{(x,y)\to \xyb}\atop{0\neq d(y,F(x))\to 0
			}}|\nabla\psi_y|(x).
			$$
			
		\end{theorem}

		The starting point for developing slope-based regularity theory was the paper
		by Az\'e-Corvellec-Lucchetti \cite{ACL} (its first version was circulated in 1998)
		who obtained a global error bound in terms of "variational pairs" that
		include slope on a metric space as a particular case. Theorem
		\ref{critgen}, and specifically the fact that the slope estimate is precise,
		was proved in \cite{AI00}
		 under somewhat stronger condition (equivalent to $Y$ being a
		length space). We refer to \cite{DA06} for a systematic exposition of the slope-based approach to local regularity.
		Theorem \ref{slsc} is a slightly modified version
		of the mentioned result of Ngai-Tron-Thera \cite{NTT12} (proved originally for
		$Y$ being a Banach space).
		
		To explain how the additional assumption on $Y$ is used to get necessity e.g. in
		Theorem \ref{critgen}, let us consider, following the original argument in
		\cite{AI00},
		$(x,y,v)$ sufficiently close to $\xb$ and $\yb$
		respectively and such that $y\neq v\in F(x)$. For any $n$ take $\del_n=
		o(n^{-1})$
		and a $v_n$ such that $d(v_n,v)\le (n^{-1}+\del_n)d(y,v)$
		and $d(v_ny)\le (1-n^{-1}+\del_n)d(y,v)$. If $Y$ is a length space such $v_n$
		can be found.
		As $F$ covers near
		$\xyb$ with modulus greater than $r$, there is a $u_n$ such that
		$v_n\in F(u_n)$ and $d(u_n,x)\le  r^{-1}d(v_n,v)\to 0$ when $n\to\infty$. We
		have
		$|d(y,v) - (d(y,v_n) + d(v,v_n))|= o(d(v_n,v))$. Therefore (as $r\xi<1$)
		$$
		|\nabla\vf_y|(x,v) \ge \lim_{n\to\infty}
		\frac{\vf_y(x,v) - \vf_y(u_n,v_n)}{\max\{d(u_n,x),\xi d(v_n,v)\}}
		\ge  \lim_{n\to\infty}
		\frac{d(v_n,v)}{r^{-1}d(v_n,v)}= r.
		$$
		Similar argument, modified as the definition of $\bpsi_y$ includes a limit
		operation,
		can be used also for the proof of necessity in Theorem \ref{slsc}.

        It should be observed that the class of locally coherent spaces is strictly bigger than the
		class of length spaces. For instance a smooth manifold in a Banach space with the induced
        metric is a locally coherent space but not a length space (unless it is a linear manifold).

		\subsection{Related concepts:  metric subregularity, calmness,
			controllability, linear recession}

		In the definitions of the local versions of the three main regularity properties
		we
		scan entire neighborhoods of the reference point of the graph
		of the mapping. Fixing one or both components of the point leads to
		new weaker concepts that differ from regularity in many respects. Subregularity
		and calmness attract much attention last years. We refer to \cite{DR} for a
		detailed study
		of the concepts mainly for mappings between finite dimensional spaces, and begin
		with
		parallel concepts relating to linear openness which are rather new in the
		context of variational analysis. We skip (really elementary) proofs of almost
		all results in this subsection.

		\begin{definition}[controllability]\label{contr}
			{\rm A set valued mapping $F:\ X\rra Y$ is said to be {\it (locally)
					controllable at} $\xyb$ if there are $\ep>0, \ga >0$ such that
				\begin{equation}\label{3.1d}
				B(\yb,r t)\subset F(B(\xb,t)),\quad {\rm if}\; 0\le t< \ep.
				\end{equation}
				The upper bound of such $r$ is the {\it rate} or
				{\it modulus of controllability} of $F$ {\it at} $\xyb$. We shall
				denote it $\cont F(\xb|\yb)$ and $\cont F(\xb)$ if $F$ is
				single-valued.}
		\end{definition}

		\begin{proposition}[Regularity vs.  controllability]\label{regvscont}
			Let $X$ and $Y$ be metric spaces, let $F: X\rra Y$ have locally complete graph,
			and let
			$\xyb\in\gr F$. Then
			\begin{equation}\label{3.d}
			\sur F(\xb|\yb)=\lim_{\ep\to 0}\inf\{{\rm contr}F(x|y):\; (x,y)\in\gr F,\;
			\max\{d(x,\xb),d(y,\yb)\}<\ep\}.
			\end{equation}
		\end{proposition}
		

		\begin{definition}[linear recession]\label{recess}
			{\rm Lets us say that $F$ {\it recedes from $\yb$ at $\xyb$ at a linear rate}
				if there are $\ep>0$ and $K\ge 0$ such that
				\begin{equation}\label{3.2d}
				d(\yb, F(x))\le K d(x,\xb),\quad {\rm if}\; d(x,\xb)<\ep.
				\end{equation}
				We shall call the lower bound of such $K$ the {\it speed of
					recession} of $F$ from $\yb$ at $\xyb$ and denote it $\ress F(\xb|\yb)$}
		\end{definition}

		The other possible way to  ``pointify'' the
		Aubin property is to fix $\xb$ and allow $(x,y)$ to change
		within $\gr F$.  Then, instead of (\ref{3.2d}) we get the
		inequality
		\begin{equation}\label{3.4d}
		d(y, F(\xb))\le K d(x,\xb) .
		\end{equation}
		
		\begin{definition}[calmness] {\rm It is said that $F: X\rra Y$ is {\it calm} at
				$\xyb$ if there are $\ep>0$, $K\ge 0$ such that (\ref{3.4d}) holds if
				$d(x,\xb)<\ep, \ d(y,\yb)<\ep$ and $y\in F(x)$. The lower bound of
				all such $K$ will be called the {\it modulus of calmness} of $F$
				at $\xyb$. We shall denote it by $\calm F(\xb|\yb)$ ($\calm
				F(\xb)$ if $F$ is single-valued).}
		\end{definition}

		Again we can easily see that {\it uniform calmness}, that is
		calmness at every $(x,y)$ of the intersection of $\gr F$ with a
		neighborhood of $\xyb$ with the same $\ep$ and $K$ for all such
		$(x,y)$, is equivalent to the Aubin property of $F$ near $\xyb$.

		\begin{definition}[subregularity]\label{subdef}{\rm
				Let $F: X\rra Y$ and $\yb\in F(\xb)$. It is said that $F$ is
				{\it (metrically) subregular} at $\xyb$ if there is a $K>0$
				such that
				\begin{equation}\label{3.6d}
				d(x,F^{-1}(\yb))\le Kd(\yb,F(x)) \quad {\rm if}\; d(x,\xb)<\ep.
				\end{equation}
				for all $x$ of a neighborhood of $\xb$. The lower bound of such $K$
				is called the {\it rate} or {\it modulus of subregularity} of $F$ at $\xyb$.
				It will be denoted $\subreg F(\xb|\yb)$.

                We say that $F$ is {\it strongly subregular} at $\xyb$ if it is subregular at the point and $\yb\not\in F(x)$ for $x\neq \xb$ of a neighborhood of $\xb$.	}
		\end{definition}

		\begin{proposition}\label{contsub} The equalities
			$$
			\subreg F(\xb|\yb)= \calm F^{-1}(\yb|\xb),\quad
			\cont F(\xb|\yb)\cdot \ress F^{-1}(\yb|\xb)=1
			$$
			always hold. If moreover, $F$ is strongly subreglar at $\xyb$, then
			$$
			\cont F(\xb|\yb)\cdot \subreg F(\xb|\yb)\ge 1.
			$$
		\end{proposition}

		\begin{theorem}[slope criterion for calmness]\label{slcalm}
			Let $X$ and $Y$ be arbitrary metric spaces, let $F: X\rra Y$ be a set-valued
			mapping
			with closed graph and let $\xyb\in\gr F$. Then
			$$
			\calm F(\xb|\yb)\ge\limsup_{y\to \yb}|\nabla \psi_y|(\xb),
			$$
			where, as earlier, $\psi_y(x)= d(y,F(x))$.
		\end{theorem}
		
		\proof Let $K>\calm F(\xb|\yb)$ then there is an $\ep >0$ such that (\ref{3.4d})
		holds, provided $d(x,\xb)<\ep$ and $y\in F(x)$. To prove the theorem, it
		is sufficient to show that $|\nabla\vf_y|(\xb)\le K$ for all $y$
		sufficiently close to $\yb$. To this end, it is sufficient to verify that
		there is a $\del>0$ such that the inequality
		$$
		d(y,F(\xb))-d(y,F(x))\le Kd(x,\xb)
		$$
		holds for all $x,y$ satisfying $d(x,\xb)<\del,\; d(y,\yb)<\del$.
		
		If $y\in F(x)$, then (\ref{3.6d}) reduces to (\ref{3.4d}). Take a positive
		$\del<\ep/2$,
		and let $x$ and $y$ be such that $d(x,\xb)<\del,\; d(y,\yb)<\del$.
		If $d(y,F(x))\ge \del$, then (\ref{3.6d}) obviously holds. If $d(y,F(x))<\del$,
		we can choose a $v\in F(x)$ such that $d(y,v)<\del$. Then
		$d(v,\yb)<\ep$ and therefore $d(v,F(\xb))\le d(x,\xb)$. Thus
		$$
		\begin{array}{lcl}
		d(y,F(\xb))-d(y,F(x))&\le& d(y,v)+d(v,F(\xb))-d(y,F(x))\\   \\
		&\le& Kd(x,\xb) + d(y,v)-d(y,F(x))
		\end{array}
		$$
		and (\ref{3.6d}) follows as $d(y,v)$ can be arbitrarily close to
		$d(y,F(x))$.\endproof

		\begin{theorem}[slope criterion for subregularity]\label{calmcrit}
			Assume that $X$ is a complete metric space. Let $F:\ X\rra Y$ be a
			closed set-valued mapping and $\xyb\in \gr F$. Assume that
			the  function $\psi(x)=d(\yb,F(x))$ is lower semicontinuous and
			there are $\ep>0$ and $ r>0$ such that
			$$
			|\nabla\psi_{\yb}|(x)=|\nabla d(\yb,F(\cdot))|(x) \ge   r,
			$$
			if $d(x,\xb)<\ep$ and $0<d(\yb, F(x))<\ep$. Then $F$ is subregular
			at $\xyb$ with modulus of subregularity (and hence the modulus of calmness
			of $F^{-1}$ at $(\yb,\xb)$)  not greater than $r^{-1}$.
		\end{theorem}

		\section{Metric theory. Perturbations and stability.}
		In this section we concentrate on two  fundamental questions:
		
		(a) what happens with regularity (and subregularity) properties of $F$ if the
		mapping is slightly perturbed?
		
		(b) how the set of solutions of the inclusion $y\in F(x,p)$ (where $F$ depends
		on a parameter $p$) depends on $(y,p)$?
		
		\noindent The answer to the second question leads us to a fairly general
		implicit
		function theorems. The key point in both cases  is that we have to require a
		certain amount of Lipschitzness of perturbations to get desirable results.
		
		\subsection{Stability under Lipschitz perturbation}

		\begin{theorem}[stability under Lipschitz perturbation]\label{stablip}
			Let $X$, $Y$ be metric spaces, let $U\subset X$ and $V\subset Y$ be open sets.
			Consider  a set-valued mapping $\Psi: X\times X\rra Y$ with closed graph
			assuming that either $X$ or the graph of $\Psi$ is complete. Let
			 $F(x)=\Psi(x,x)$.  Suppose  that
			
			(a) for any $u\in U$ the mapping $\Psi(\cdot,u)$ is Milyutin regular on $(U,V)$
			with modulus of surjection greater than $r$, that is for any $x\in U$, any $v\in
			\Psi(x,u)$ and any $y\in \bo(v,rt)\cap V$ with $t< d(x,X\backslash U)$ there is
			an $x'$
			such that   $d(x,x')\le r^{-1}d(y,v)$ and $ y\in F(x');$
			
			(b) for any $x\in U$ the mapping $\Psi(x,\cdot)$ is pseudo-Lipschitz on $(U,V)$
			with modulus $\ell<r$,
			that is for any $u,w \in U$
			$$
			{\rm ex}(\Psi(x,u)\cap V,\Psi(x,w))< \ell d(u,w).
			$$
			\noindent Then $F(x)=\Psi(x,x)$ is Milyutin regular on $(U,V)$ with
			$\sur_{m}F(U|V)\ge r-\ell$.
		\end{theorem}
		
		
		\proof  We shall consider only the case of complete $\gr\Psi$.
		According to the general regularity criterion of Theorem \ref{gencrit}
		all we have to show is that there is a $\xi>0$ such that, given $(x,v)\in gr F$
		and $y$ such that
		$x\in U$, $y\in V$ and  $0<d(y,v)< rm(x)$, there is another point
		$(x',v')\neq (x,v)$ in the graph of $F$ such that

		\centerline{$d(y,v')\le d(y,v)-(r-\ell)\max\{d(x,x'),\xi(v,v'))\}$.}
		
		\vskip 1mm
		
		\noindent We have by (a): $B(v,rt)\cap V\subset \Psi(B(x,t),x)$ if $t<m (x)$. As
		$d(y,v)< rm (x)$, it follows that there is a $x'\in \B(x,t)$ such that
		$y\in \Psi(x',x)$ and $d(x,x')\le r^{-1}d(y,v)$.
		
		Clearly,  $x'\in U$. Therefore by (b) $d(y,\Psi(x',x'))< \ell d(x,x')$. This
		means that
		there is a $v'\in F(x')$ such that
		$$
		d(y,v')\le \ell d(x,x')\le \frac{\ell}{r} d(y,v).
		$$
		Take $\xi<(r+\ell)^{-1}$. Then
		$$
		\xi d(v,v')\le (r+\ell)^{-1}(d(v,y)+d(y,v'))\le
		(r+\ell)^{-1}\Big(1+\frac{\ell}{r}\Big)d(y,v)
		= \frac{1}{r}d(y,v).
		$$
		Thus $\max\{d(x,x'),\xi d(v,v')\}\le r^{-1}d(y,v)$ and we have
		$$
		d(y,v')< (\ell/r)d(y,v)= d(y,v)- \frac{r-l}{r}d(y,v)\le
		d(y,v)-(r-l)\max\{d(x,x'),\xi d(v,v')\}.
		$$
		as needed.\endproof
		

		
		\begin{corollary}[Milyutin's perturbation theorem \cite{DMO}]\label{milt1}
			Let $X$ be a metric space, let $Y$ be a normed space and  $F: X\rra Y$ and
			$G: X\rra Y$ We assume that either the graphs of $F$ and $G$ are complete or $X$ is a complete space.
			Let further $U\subset X$ be an open set such that $F$ is Milyutin regular on
			$U$ with $\sur F(U)\ge r$ and $G$ is (Hausdorff) Lipschitz with $\lip
			G(U)\le\ell<r$. If either $F$ or $G$ is single-valued continuous on $U$, then   $F+G$ is Milyutin regular on $U$ and
			$\sur (F+G)(U)\ge r-\ell$.
		\end{corollary}
		
		\proof Apply the theorem to $\Psi(x,u)= F(x)+G(u)$.\endproof

		To state a local version of the theorem, we need the following
		
		\begin{definition}[uniform regularity]\label{unreg}{\rm
				Let $P$ be a topological space, let   $F: P\times X\rra Y$, let $\pb\in P$,
				and
				let $\xyb\in\gr F(\pb,\cdot)$. We shall say that $F$ is regular near $\xyb$
				{\it uniformly}
				in $p\in P$ near $\pb$  if for any $r<\sur F(\pb,\cdot)(\xb|\yb)$
				there are $\ep>0$ and a neighborhood $W\subset P$ of $\pb$ such that for any
				$p\in W$ and
				any $x$ with $d(x,\xb)<\ep$
				$$
				B(F(p,x),rt)\cap B(\yb,\ep)\subset F(p,B(x,t)), \quad {\rm if}\; 0\le t<\ep.
				$$
			}
		\end{definition}
		
		\begin{theorem}[stability under Lipschitz perturbations: local
			version]\label{stabliploc}
			Let $X$, $Y$, $\Psi: X\times X\rra Y$ and $F(x)= \Psi(x,x)$ be as in Theorem \ref{stablip}, and let
			$\xyb\in\gr F$. We assume that
			
			(a) $\Psi(\cdot, u)$ is regular  near $\xyb$ uniformly in $u$ near $\xb$;
			
			(b) $\Psi(x,\cdot)$ is pseudo-Lipschitz near $\xyb$  uniformly in $x$ near
			$\xb$.
			
			\noindent If $\lip \Psi(\xb,\cdot)(\xb|\yb)<\ell<r<\sur
			\Psi(\cdot,\xb)(\xb|\yb)$,
			then $F$ is regular near $\xyb$ with modulus of surjection greater than
			$r-\ell$.
			
		\end{theorem}
		The last theorem in turn immediately implies
		Milyutin's theorem and its versions correspond to $\Psi(x,y)= F(x)+g(y)$ with
		$g$ being single-valued
		Lipschitz. The following corollary from the theorems is straightforward

		\begin{theorem}[Milyutin's perturbation theorem - local version]\label{milt}
			Let $X$ be a metric space, let $Y$ be a normed space, and let $F: X\rra Y$ and
			$G: X\rra Y$.
			Given $\xb\in\dom F\cap\dom G$,  $\yb\in F(\xb),\; \zb\in G(\xb)$, we assume
			that $F$
			is regular near $\xyb$ with $\sur F(\xb|\yb)\ge r$ and $G$ has the Aubin
			property  near
			$(\xb,\zb)$ with $\lip G(\xb|\zb)\le\ell$. If either $F$ or $G$ is
			single-valued continuous on its
			domain and the graph of the other is complete in the product metric, then
			
			$$
			\sur (F+G)(\xb,\yb+\zb)\ge r-\ell.
			$$
		\end{theorem}
		
		\proof Set $\Psi(x,y)=F(x)+ G(y)$. It is an easy matter to check that the
		conditions of
		Theorem \ref{stabliploc} are valid.
		\endproof
		
		As an immediate consequence of the last theorem we mention a stronger version of the Lyusternik-Graves theorem stating
		that its condition is not only sufficient but also necessary for regularity is an
		immediate corollary of the
		last theorem.

		\begin{corollary}[Lyusternik-Graves from Mulyutin]\label{neclg}
			Let $X$ and $Y$ be Banach spaces,  and let $F:X\to Y$ be strictly
			differentiable at $\xb$.
			Then   $\sur F(\xb)=C(F'(\xb))$.
		\end{corollary}
		
		\proof Indeed, let $X,Y$
		be Banach spaces, and let
		$F:X\to Y$ be strictly differentiable at $\xb$. Set $g(x) = F(x) -
		F'(\xb)(x-\xb)$. As $F$
		is strictly differentiable at $\xb$, the Lipschitz constant of $g$ at $\xb$ is
		zero which by
		Milyutin's theorem means that the moduli of
		surjection of $F$ at $\xb$ and $F'(\xb)$ coincide.
		\endproof

		We observe next that in Theorem \ref{milt}  one of the mappings is assumed
		single-valued.
		This assumption is essential. With both mappings set-valued the result may be
		wrong as the following example shows.

		\begin{example}[cf. \cite{DR}]\label{contr1}{\rm
				Let $X=Y=\R$, $G(x,y)=\{x^2,-1\},\; F(x)=\{-2x,1\}$.
				It is easy to see that $F$ is regular near $(0,0)$ and $G$ is
                Lipschitz in the Hausdorff metric.
				On the other hand,
				$$
				\Phi(x)=\{x^2-2x,x^2+1,-2x-1,0\}
				$$
				is not even regular at $(0,0)$. Indeed $(\xi,0)\in\gr \Phi$ for any $\xi$.
				However, if $\xi\neq 0$,
				then  the $\Phi$-image of a sufficiently small neighborhood of $\xi$ does not
				contain points
				of a small neighborhood of zero other than zero itself. }
		\end{example}
		
		Perturbation analysis of regularity properties was initiated by Dmitruk-Milyutin-Osmolovski in \cite{DMO} with a proof of  a global version
		of Theorem \ref{milt1} (attributed in \cite{DMO} to Milyutin)
		with both the mapping and the perturbation single valued. The first perturbation
		result for set-valued mappings was proved probably by Ursescu \cite{CU96}
		(see also \cite{AI00}). Observe that  global theorems are valid for Lipschitz set-valued perturbations as well.
		
		Till very recently the main attention was devoted to additive perturbations
		into a linear range space, especially in connection with implicit function theorems for generalized equations - see e.g. \cite{AB08,DR}. Interest to non-additive Lipschitz set-valued perturbations of set-valued mappings appeared just a few years ago, partly in connection with fixed point and coincidence theorems
		\cite{AAGDO,DF12,AI11a,AI14}
		
		The Graves theorem can  be viewed as a perturbation theorem for a {\it linear} regular operator.
		For that reason in some publications (e.g. \cite{DF11,DR})  this theorem is called
		"extended Lyusternik-Graves theorem". I believe the name "Milyutin theorem"
		is  adequate.  It is quite obvious that Graves did not have in mind the perturbation issue and was interested only in a quality of approximation needed to get the result. (Tikhomirov and I a had similar idea when proving the metric regularity counterpart of the Graves theorem for \cite{IT} without any knowledge of the Graves' paper.) And the fact that the Lipschitz
		property of the perturbation as the key for the estimate was explicitly
		emphasized in \cite{DMO}. Note also that even Corollary \ref{neclg} cannot be obtained from the Graves theorem.

		Milyutin's theorem can also be viewed as a  regularity result for  a composition
		$\Phi(x,F(x))$, where $\Phi(x,y)= G(x)+y$. Theorems \ref{stablip} and \ref{stabliploc} can be applied to prove regularity of more general compositions, with arbitrary $\Phi$, just by taking $\Psi(x,u)=\Phi(x,F(u))$.
		However,  a certain caution is needed to guarantee that such a $\Psi$
		satisfies the required assumptions (as say in \cite{AI11a} where $\Phi(x,\cdot)$
		is assumed to be an isometry or in \cite{DS12} where a certain "composition stability" is a priori assumed). Corollary \ref{neclg} was
		first stated in \cite{ALD96}  with a direct proof, not using Milyutin's theorem.

		\subsection{Strong regularity and metric implicit function theorem.}
		
		Generally speaking, the essence of the inverse function theorem
		is already captured by the main Equivalence Theorem \ref{equiv3}.
		But in  view of the very special role of the inverse and implicit
		function theorems in the classical theory, it seems appropriate to
		make the connection with the classical results more transparent.
		
		So let $F(x,p): X\times P\rra Y$. We shall view $P$ as a parameter space.
		Let $S(y,p)=\{ x\in X:\; y\in F(x,p)\}$  stand for the solution mapping
		of the inclusion $y\in F(x,p)$. In all
		theorems to follow we consider $Y\times P$ with an $\ell^1$-type distance
		$$
		d_{\al}^1((y,p),(y',p'))= \al d(y,y')+ d(p,p'),
		$$
		where $\al$ will be further determined  by Lipschitz moduli of mappings
		involved.
		
		\begin{theorem}[general proposition on implicit functions]\label{genimp}
			We assume that $\yb\in F(\xb,\pb)$ and
			$F$ satisfies the following conditions: there are constants $K>0$, $\al >0$ and
			a
			sufficiently
			small $\ep> 0$ such that  the following
			relations hold:
			
			\vskip 1mm
			
			(a) $F(\cdot,p)$ is regular near ($\xyb,\pb$) uniformly in $p$  with the rate of
			metric regularity not grater than $K$;
			\vskip 1mm
			
			(b) $F(x,\cdot)$  is pseudo-Lipschitz near $(\xb,(\pb,\yb))$ uniformly in $x$
			with the Lipschitz modulus not greater than $\al$.

			\vskip 1mm
			
			\noindent   Then $S$ has the Aubin property near $((\yb,\pb),\xb)$ with the
			Lipschitz modulus
			with respect to the metric $d_{\al}^1$ in $Y\times P$ not greater than
			$\reg F(\cdot,\pb)(\xb,\yb)$.
			
			In particular, if we are interested in solutions of the inclusion $\yb\in
			F(x,p)$  (with  fixed $\yb$),
			then  under the assumption of the theorem the solution mapping
			$p\mapsto S_{\yb}(p)$  has the Aubin property near $(\pb,\xb)$ with Lipschitz
			modulus not exceeding $K\al$.
		\end{theorem}
		
		\proof  As $F(\xb,\pb)\neq\emptyset$, the uniform pseudo-Lipschitz property
		implies that
		$S(y,p)\neq\emptyset$ for $(y,p)$ close to $(\yb,\pb)$. If now $y\in F(x,p)$, then
		$$
		\begin{array}{lcl}
		d(x,S(y',p'))&\le& K d(y',F(x,p'))\le K\big(d(y,y') + d(y,F(x,p'))\big)\\
		&\le&  K\big(d(y,y') +\al d(p,p')\big)= K d_{\al}^1((y,p),(y',p'))\\
		&=&K\al(d(p,p')+\al^{-1}d(y,y')),
		\end{array}
		$$
		and the proof is completed.\endproof

		\begin{definition}\label{sregdef} {\rm Let $F: X\rra Y$, and let $\yb\in
				F(\xb)$. We say that $F$ is {\it strongly (metrically)
					regular} near $\xyb\in\gr F$  if for some $\ep>0,\ \del >0$ and $K\in[0,\infty)$
				\begin{equation}\label{2.4.1}
				B(\yb,\del)\subset F(B(\xb,\ep)) \quad \&\quad d(x,u)\le Kd(y,F(x))
				\end{equation}
				\noindent whenever $x\in B(\xb,\ep)$, $u\in B(\xb,\ep)$  and $y\in F(u)\bigcap
				B(\yb,\del)$.

				We shall also say following \cite{DR} that $F$ {\it has a single-valued
					localization near} $\xyb$ if there are $\ep> 0,\ \del>0$ such that
				the restriction of $F(x)\cap B(\yb,\del)$ to $B(\xb,\ep)$ is
				single-valued.  If in addition, the restriction is Lipschitz continuous, we say
				that  $F$
				has {\it Lipschitz localization} near $\xyb$ .}
		\end{definition}
It is obvious from the definition that strong regularity implies regularity: the second relation in
(\ref{2.4.1}) is clearly stronger than metric regularity.

		\begin{proposition}[characterization of strong regularity]\label{sregelem}
			Let $F: X\rra Y$ and $\xyb\in\gr F$. Then the following properties
			are equivalent
			
			\vskip 1mm
			
			(a) $F$ is strongly regular near $\xyb$;
			
			\vskip 1mm
			
			(b)  there are $\ep>0$ and $\del>0$ such that $B(\yb,\del)subset F(B(\xb,\ep))$ and
			\begin{equation}\label{2.4.2}
			F(x)\bigcap F(u)\bigcap B(\yb,\del)=\emptyset,
			\end{equation}
			\noindent whenever $u\neq x$ and both $x$ and $u$ belong to $B(\xb,\ep)$;

			\vskip 1mm
			
			(c) $F$ is regular near $\xyb$ and there  are $\ep>0,\ \del>0$ such that
			$F^{-1}$ has a
			
			single-valued localization near $(\yb,\xb)$;
			
			\vskip 1mm
			
			(d) $F^{-1}$ has a Lipschitz localization  $G(y)$ near $(\yb,\xb)$.
			In particular $y\in F(G(y))$
			
			for all $y$ of a neighborhood of $\yb$.

			\vskip 1mm
			
			\noindent Moreover, if $F$ is  strongly regular near $\xyb$, then the
			lower bound of $K$ for which the second part  of (\ref{2.4.1}) holds and the
			Lipschitz modulus of its Lipschitz localization $G$ at $\yb$ coincide with
			$\reg F(\xb|\yb)$.
		\end{proposition}

		\begin{theorem}[persistence of strong regularity under Lipschitz
			perturbation]\label{strlip}
			We  consider a set-valued mapping $\Phi:X\rra Y$  with   complete graph, and a
			(single-valued) mapping
			$G: X\times Y\to Z$. Let $\yb\in\Phi(\xb)$ and $\zb=G(\xb,\yb)$. We assume that
			
			(a) $\Phi$ is strongly regular near $\xyb$ with $\sur \Phi(\xb|\yb)>r$;
			
			(b) $G(x,\cdot)$ is an isometry from $Y$ onto $Z$ for any $x$ of a neighborhood
			of $\xb$;
			
			(c) $G(\cdot,y)$ is Lipschitz with constant $\ell<r$ in a neighborhood of
			$\xb$, the same for
			
			all  $y$ of a neighborhood of $\yb$.
			
			\noindent Set $F(x)=G(x,\Phi(x))$. Then $F$ is strongly regular near
			$(\xb,\zb)$.
			
			In particular, if $Y$ is a normed space, $\Phi$ is strongly regular near $\xyb\in\gr \Phi$ and $G(x,y)= g(x)+ y$ with
			$\lip g(\xb)<\sur \Phi(\xb|\yb)$, then $F(x) = \Phi(x)+g(x)$
			is strongly regular near $(\xb,\yb+g(\xb))$.
		\end{theorem}

		\begin{remark}\label{remstrong}{\rm
				It is to be observed in connection with the last theorem that strong
				regularity is not preserved
				under  set-valued perturbations like those in  Theorem \ref{stablip}. Here is
				a simple example:
				$$
				\Psi (x,u)= x+ u^2[-1,1] \; (x,u\in \R),\quad \xb = 0.
				$$
				Clearly $\Psi(\cdot,0)$ is strongly regular but $F(x)= x+x^2[-1,1]$ is of
				course regular but not
				strongly regular.
				
				It follows that strong regularity is somewhat less robust compare to the
				standard regularity.}
		\end{remark}

		\begin{theorem}[implicit function theorem - metric version]\label{mimpl}
			Assume in addition to the assumptions of Theorem \ref{genimp} that
			\begin{equation}\label{2.4.3}
			F(x,p)\cap F(x',p)\cap\bo (\yb,\ep)=\emptyset\quad\forall \;
			x,x'\in\bo(\xb,\ep),\;,x\neq x',\; p\in\bo(\pb,\ep).
			\end{equation}
			Then the solution map $S$ has a Lipschitz localization $G$ near
			$((\pb,\yb),\xb)$ with $\lip G(\pb,\yb)\le K$
			(with respect to the $d_{\al}^1$-metric in $Y\times P$. In particular $z\in
			F(S(p,y),y)$ for all $(p,y)$
			of a neighborhood of $(\pb,\yb)$.
		\end{theorem}

		The conclusion is already very
		similar to the conclusion of the classical implicit function theorem.
		Indeed, it contains precisely the same information about the solution,
		namely its uniqueness in a neighborhood and its Lipschitz continuity (replacing
		differentiability) with the Equivalence Theorem \ref{equiv3}
		providing, along with the concluding part of Proposition \ref{sregelem} an
		estimate for the Lipschitz constant of the solution map (replacing the formulas
		for partial derivative in the classical theorem).
		Moreover, the proof below is based on the same main idea as the proof
		of the classical theorem, say the second proof in \cite{DR}.

		\proof
		Consider the set-valued mapping $\Phi$ from $X\times P$ into $P\times Y$.
		defined by
		$$
		\Phi(x,p)=\{p\}\times F(x,p).
		$$
		Then $(\pb,\yb)\in\Phi(\xb,\pb)$. We claim that $\Phi$ is strongly regular near
		$((\xb,\pb),(\pb,\yb))$. Indeed, we have for $x,\ p,\  y$ sufficiently close to
		$\xb,\pb,\yb$
		\begin{equation}\label{2.4.4}
		\Phi^{-1}(x,y)=\{ p\}\times  S(p,y)
		\end{equation}
		By Theorem \ref{genimp} $S$ has the Aubin property at $((\pb,\yb),\xb)$.  This
		obviously implies
		that  $\Phi^{-1}$ has the Aubin property at  $((\pb,\yb),(\xb,\pb)$.  The latter
		means that
		$\Phi$ is  regular at $((\xb,\pb),(\pb,\yb))$.
		
		On the other hand,  $(p,y)\in\Phi(x,p)\cap\Phi(x',p')$  means that
		$p=p'$ and $y\in F(x,p)\cap F(x',p)$, so that  (\ref{2.4.3})  may happen only if
		$x=x'$.
		This proves the claim.
		
		By Proposition \ref{sregelem} there is a Lipschitz localization
		of $\Phi^{-1}$ defined in a neighborhood of $(\pb,\yb)$.
		By (\ref{2.4.3}) this localization has the form $(p,G(p,y))$, where
		$G(p,y)\in S(p,y)$. Thus $G$ is a Lipschitz localization
		of $S$ and  by Theorem \ref{genimp} its Lipschitz constant is not greater
		than $K$.
		\endproof

		\begin{theorem}[metric infinitesimal implicit function theorem]\label{impinf}
			Let $\yb\in F(\xb,\pb)$,  and assume that there are $\xi>0,\ r>0,\ \ell>0, \
			\ep >0$ are such that for all  $x,y,p,v$ satisfying
			$$
			d(x,\xb)<\ep,\; d(y,\yb)<\ep,\; d(p,\pb)<\ep,
			$$
			either $\gr F$ is complete and

			(a$_1$) \ $|\nabla_{\xi}\vf_y(\cdot,p)|(x,v)>r$ \quad {\rm if} $v\in F(x,p)$
			and $ d(y,v)>0$
			
			\noindent or $X$ is a complete space and
			
			(a$_2$) \  $|\nabla\bpsi_y(\cdot,p)|(x)> r$ \quad {\rm if} $\; \bpsi_y(x,p)>0$
			
			\vskip 1mm
			
			\noindent holds along with
			
			\vskip 1mm
			
			(b) \ $|\nabla\psi_y(x,\cdot)|(p) <\ell d(p,p') $, if $y\in F(x,p')$ for some
			$p'\in\bo (\pb,\ep)$.
			
			\vskip 1mm

			\noindent Then $S$ has the Aubin property near $(\yb,\pb)$ with $\lip
			S((\yb,\pb)|\xb) \le r^{-1}$
			if   $Y\times P$ is considered with the distance $d_{\ell}^1((y,p),(y',p'))=
			\ell d(p,p')+ d(y,y')$.
			
		\end{theorem}
		
		The proof of the theorem consists in verifying the assumptions of Theorem
		\ref{genimp}
		for all $(x,y,p)$ of a neighborhood of $(\xb,\pb,\yb)$ and $p'$ close to $\pb$.

		The next theorem is an infinitesimal counterpart of Theorem \ref{mimpl}.
		
		\begin{theorem}\label{simpinf}
			In addition to the conditions of Theorem \ref{impinf} we assume that
			
			\vskip 1mm
			
			(c) \ $|\nabla\psi_y(\cdot,p)|(x) >0$ \quad if $y\in F(x',p)$ for some $x'\neq  x$.
			
			\vskip 1mm
			
			\noindent Then $S$ has a Lipschitz localization $G$ in a neighborhood of
			$(\pb,\zb)$ with
			$G(\pb,\yb)=\xb$ and the Lipschitz constant (with respect to the
			$d_{\ell}^1$-metric in
			$P\times Y$) not exceeding $r^{-1}$.
		\end{theorem}
		
		\proof Indeed, it follows from (c) that $y\not\in F(x,p)$ that is
		$(F(x,p)\cap F(x',p))\cap\bo(\yb,\ep)=\emptyset$
		for   $x, x'$ close to $\xb$ and $p$ close to $\pb$ and the reference to
		Theorems
		\ref{impinf} and \ref{mimpl}
		completes the proof.\endproof
		
		There have been numerous publications extending, one way or another,  the implicit function theorem to
		settings of variational analysis, see e.g	 \cite{AB08,DR,BDG01,AI00,LZ99,NT04,NTT12}.
		Most of them deal with Banach spaces and/or specific classes of mappings, e.g. associated with generalized equations.
		It should be also said that some results named ``implicit function theorem"
		are rather parametric regularity or subregularity theorems giving uniform
		(w.r.t parameter) estimates for regularity rates of a mapping depending on a parameter.
		
		The concept of strong regularity was introduced by Robinson in \cite{SMR80}.
		A number of characterizations of strong regularity can be found in
		\cite{DR}. It is appropriate to mention (especially because we do not discuss
		these questions in the paper) that there are certain important classes of mappings for which regularity and strong regularity are equivalent. Such are
		monotone operators, in particular subdifferentials of convex functions, or
		Kojima mappings associated with constrained optimization \cite{DR,KK02}.

		\section{Banach space theory.}
		Needless to say that the vast majority of applications of the theory of metric regularity
relate to problems naturally stated in Banach spaces. Variational analysis and metric regularity theory  in Banach spaces are distinguished by

(a) the existence of an approximation mechanisms, both primal and dual, using homogeneous mappings
(graphical derivatives and coderivatives) in case of set-valued mappings or directional subderivatives and subdifferentials for functions;

(b) the possibility of separable reduction for metric regularity that allows to reduce much of analysis to mappings between separable spaces;

(c) the existence of a class of linear perturbations, most natural and interesting in many cases.

		\subsection{Techniques of variational analysis in Banach spaces.}
		
		\subsubsection{Homogeneous set-valued mappings.}
		\begin{definition}\label{normsvm}
			{\rm A set valued mapping $\ch: X\rra Y$ is
				{\it homogeneous} if its graph is a pointed cone. The latter means that $0\in
				H(0)$.
				The mapping 
				$$
				\ch^*(y^*)=\{ x^*:\; \lan x^*,x\ran-\lan y^*,y\ran \le 0,\;\forall\;
				(x,y)\in\gr \ch\}
				$$
				is called {\it adjoint} or {\it dual} to $\ch$ (or the {\it dual convex
					process} as it is often called
				for the reasons to be explained in the next chapter).
				It is an easy matter to see, that
				
				\vskip 1mm
				
				\centerline{$\gr \ch^*=\{(y^*,x^*):\; (x^*,-y^*)\in \big(\gr \ch)^{\circ}\}.$}
				
				\vskip 1mm
				
				With every homogeneous mapping $\ch$ we associate  the {\it upper norm}
				
				\vskip 1mm
				
				\centerline{$\| \ch\|_+= \sup\{ \| y\|:\; y\in \ch(x),\; x\in\dom\ch,\; \| x\|\le 1\},$}
				and the {\it lower norm}
				
				\centerline{$\| \ch\|_-=\sup_{x\in B\cap\dom\ch}\inf \{\| y\|:\; y\in \ch(x)\}=\sup_{x\in B\cap\dom\ch} d(0,\ch(x))$.}
				
				\vskip 2mm
				
				\noindent
				For single-valued mappings with $\dom\ch= X$ both quantities coincide and we
				may speak about the {\it norm} of
				$\ch$. The mapping $\ch$ is {\it bounded} if $\| \ch\|_+<\infty$. This
				obviously means that
				there is an $r>0$
				such that $\ch(x)\subset r\| x\|B_Y$ for all $x$.
				
				Very often however, in the context of regularity estimates, it is more
				convenient to
				deal with different quantities defined by way of the norms as follows:
				$$
				C(\ch) = \| \ch^{-1}\|_-^{-1}\quad {\rm and}\quad C^*(\ch)= \|
				\ch^{-1}\|_+^{-1}.
				$$
				The quantities are  respectively called the {\it Banach constant} and the {\it  dual Banach constant}	of $\ch$.
				To justify the terminology, note that for linear operators they coincide with
				the Banach constants
				introduced for the latter in the first section.}
		\end{definition}
		The proposition below containing important geometric interpretation of the concepts shows   that the Banach constants are actually very natural objects..

		\begin{proposition}[cf. Proposition \ref{calca}]\label{dualb}For any homogeneous
			$\ch: X\rra Y$
			$$
			\begin{array}{l}
			C(\ch)={\rm contr}\ch(0|0)=\sup\{ r\ge 0:\; rB_Y\subset \ch(B_X)\};\\       \\
			C^*(\ch)= (\subreg\ch(0|0))^{-1}= \inf \{ \| y\|:\; y\in \ch(x),\; \| x\|= 1\}
			=\dis\inf_{\| x\|=1} d(0,\ch(x)).
			\end{array}
			$$
		\end{proposition}

		\proof The equality ${\rm contr}\ch(0|0)=\sup\{ r\ge 0:\; rB_Y\subset
		\ch(B_X)\}$ follows
		from homogeneity of $\ch$. On the other hand, saying that
		$rB_Y\subset \ch(B_X)$ is the same as saying  that for any $y$ with $\|
		y\|=r$
		there is an $\| x\|$ with $\| x\|\le 1$ such that $x\in\ch^{-1}(y)$ which means
		that $\| \ch^{-1}\|_-\le r^{-1}$ and therefore $C(\ch)\ge {\rm contr}\ch(0|0) $.
		Likewise,  $\|\ch^{-1}\|_- < r^{-1}$ means that for any $y$ with $\| y\|=1$
		there is an
		$x$ with $\| x\|\le r^{-1}$ such that $y\in\ch(x)$ from which we get
		that $rB_Y\subset \ch(B_X)$ and the first equality follows.

		To prove the second equality, consider first the case $C^*(H)<\infty$. Then
		$$
		\begin{array}{lcl}
		C^*(\ch)&=&\dis\inf_{\| y\|=1}\inf\{\| x\|^{-1}: \; x\in \ch^{-1}(y)\}\\
		&=&\inf\{\| y\|:\; y\in \ch(x), \| x\| = 1\}.
		\end{array}
		$$
		
		If $C^*(\ch)=\infty$, and therefore $\| \ch^{-1}\|_+=0$, then for any $y$ the
		set
		$\ch^{-1}(y)$ is either empty
		(recall our convention: $\inf \emptyset =\infty, \; \sup\emptyset =0$)
		or contains only the zero vector.   Hence  the domain of $\ch$ is a singleton
		containing
		the origin. It follows that
		$\inf\{ \| y\|:\; y\in \ch(x),\; \| x\|=1\}=\inf\emptyset=\infty$.

This proves the left equality. Consider again the case
$C^*(\ch)>0$. Then $\| \ch^{-1}\|_+<\infty$ and consequently, $\ch^{-1}(0)=\{ 0\}$.
It follows that $d(x,\ch^{-1}(0))= \| x\|$. Setting $K=(C^*(\ch))^{1}$, we get
for any $x$ with $\| x\|=1$:
$$
Kd(0,\ch(x))\ge 1= \| x\|= d(x,\ch^{-1}(0)
$$
and on the other hand for any $K'<K$ we can find an $x$ with $\| x\|=1$ such that $K'd(0,\ch(x))<1$.
It follows that $K= \subreg \ch(0|0)$. The case $C^*(\ch)=0$ is treated as above.
\endproof

		\begin{corollary}\label{constin}
			For any homogeneous mappings $\ch: X\rra Y$ and $\ce: Y\rra Z$
			$$
			C(\ce\circ\ch)\ge C(\ce)\cdot C(\ch).
			$$	
		\end{corollary}
		
		\proof Take $\rho< C(\ch)$. Then $\rho(B_Y)\subset \ch(B_X)$ and therefore
		$$
		\begin{array}{lcl}
		C(\ce\circ \ch)&=& \sup \{ r\ge 0: \; r B_Z\subset (\ce\circ\ch)(B_X)\}\\
		&\ge& \sup \{ r\ge 0: \; r B_Z\subset \ce(\rho B_Y)\} =\rho C(\ce)
		\end{array}
		$$
		and the result follows. \endproof

		We shall see that the tangential (primal) regularity estimates
		are stated in terms of Banach constants  of contingent derivatives of
		the mapping while the subdifferential  estimate need dual Banach constants
		of coderivatives. The following theorem is the first indicator that
		(surprisingly!)
		the dual estimates can be better.

		\begin{theorem}[basic inequality for Banach constants]\label{bancon} For any
			homogeneous
			set-valued mapping $H: X\rra Y$
			$$
			C^*(\ch^*)\ge C(\ch)\ge C^*(\ch).
			$$
		\end{theorem}
		
		\noindent Note that for linear operators we have equality -- see Proposition
		\ref{calca}. In
		the next section
		we shall see that the equality also holds for convex processes and some other
		set-valued mappings.
		
		\proof
		The right inequality is immediate from the definition.
		If $C(\ch)=\infty$, that is $\| \ch^{-1}\|_-=0$, then for any $y\in Y$ there is
		a sequence
		$(x_n)\subset X$ norm converging to zero and such that $y\in \ch(x_n)$. It is
		easy to see
		that in this case
		\begin{equation}\label{5.1}
		\ch^*(y^*)=\left\{\begin{array}{lcl}\emptyset,&{\rm if}& y^*\neq 0;\\
		X^*,&{\rm if}& y^*=0,
		\end{array}\right.
		\end{equation}
		that is $(\ch^*)^{-1}\equiv \{0\}$, $\| (\ch^*)^{-1}\|*^+=0$ and hence
		$C^*(H^*) =\infty$.
		
		Let now  $\infty> C(\ch)= r >0$. Set $\la = r^{-1}$. Then
		$\| \ch^{-1}\|_-=\la$ so that
		for any $y$ with $\|y\|=1$ and any
		$\ep>0$ there
		is an $x$ such that $\| x\|\le \la+\ep$ and $y\in \ch(x)$. Let now
		 $x^*\in\ch^*(y^*)$, that
		is $\lan x^*,x\ran -\lan y^*,y\ran\le 0$  if $y\in \ch(x)$. Take $y\in S_Y$ such
		that
		$\lan y^*,y\ran\le (-1+\ep)\| y^*\|$ and choose an $x\in\ch^{-1}(y)$ with $\| x\|\le \la +\ep$.
		Then
		$$
		-(\la+\ep)\| x^*\|\le \lan x^*,x\ran\le  \lan y^*,y\ran \le (-1+\ep)\|y^*\|
		$$
		that is $(\la +\ep)\| x^*\|\ge (1-\ep)\| y^*\|$. As $\ep$ can be chosen arbitrarily close to zero this implies that
		$\| (\ch^*)^{-1}\|_+\le r^{-1}$ and
		therefore  $C^*(\ch^*)\ge r= C(\ch)$.  \endproof
		
		The following property plays an essential role in future discussions.
		\begin{definition}[non-singularity]\label{sing}{\rm
				We say that $\ch$ is {\it non-singular} if $C^*(\ch)>0$. Otherwise we shall
				call
				$\ch$ {\it singular}.}
		\end{definition}
		
		We conclude the subsection with showing that regularity of a homogeneous mapping near the origins implies its global regularity.
		\begin{proposition}\label{globa}
			Let $X$ and $Y$ be two Banach spaces, and let $F: X\rra Y$ be a homogeneous set-valued
			mapping. If $F$ is regular near $(0,0)$, then it is globally regular with the same rates.
		\end{proposition}
		
		\proof
		By the assumption, there are $K>0$ and $\ep >0$ such that
		$d(x,F^{-1}(y))\le Kd(y,F(x))$ if $\max\{\| x\|,\| y\|  \}<\ep$. Let now $(x,y)$ be an arbitrary point of the graph. Set $\| m\|= \max\{\| x\|,\| y\|  \}$, and let
		$\mu<\ep/m$. Then
		$$
		\mu d(x,F^{-1}(y))=d(\mu x,F^{-1}\mu y)\le d(\mu y, F(\mu x))=\mu d(\mu y,F(\mu x))
		$$
		whence $d(x,F^{-1}(y))\le Kd(y,F(x))$.\endproof

		The norms for homogeneous multifunctions were originally introduced
		first by Rockafellar \cite{RTR67} and Robinson \cite{SMR72} in the context of convex processes (lower norm)
		and then by Ioffe \cite{AI81} (upper norm for arbitrary homogenous maps)   and Borwein  \cite{JMB83} (upper norm and duality for convex processes -see also
		\cite{JMB86a,BoL,DR}).
		The dual Banach constant $C^*$  was also introduced in \cite{AI81}.
		The meaning of the primal constant has undergone some evolution since it first
		appeared in \cite{AI81}. The $C(\ch)$ introduced here is reciprocal to that
		in \cite{AI87} mainly because  the connection of Banach constants with the norms
		of homogeneous mappings makes the present definition more natural.

		\subsubsection{Tangent cones and contingent derivatives}

		Given a set $Q\subset X$ and an $\xb\in Q$. The {\it tangent (or contingent)  cone}
		$T(Q,\xb)$ is the collection of $h\in X$ with the following property: there are
		sequences of $t_k\searrow 0$ and $h_k\to h$
		such that $\xb+t_kh_k\in Q$ for all $k$.
		If $F: X\rra Y$ then the {\it contingent} or {\it graphical  derivative} of $F$ at $\xyb$
		is the set-valued mapping
		$$
		X\ni h\mapsto DF\xyb(h)=\{v\in Y:\; (h,v)\in T(\gr F,\xyb)\}.
		$$
		
		Let now $f$ be a function
		on $X$ finite at $\xb$. The function
		$$
		h\mapsto f^-(\xb;h)=\liminf_{(t,h')\to (0^+,h)}t^{-1}(f(\xb+th')-f(\xb))
		$$
		is called the
		{\it Dini-Hadamard lower directional derivative} of $f$ at $\xb$. This function
		is either lsc and equal to zero at the origin or identically equal to $-\infty$.
		The latter of course cannot happen if $f$ is Lipschitz near $\xb$.
		
		The connection between the two concepts is very simple: $h\in T(Q,\xb)$
		if and only if
		$d^-(\cdot,Q)(\xb;h)=0$ and
		$\alpha = f^-(\xb;h)$ if and only if $(h,\alpha)\in T(\epi f,(\xb,f(\xb)))$.

		If $F: X\rra Y$ then the {\it contingent derivative} of $F$ at $\xb$
		is the set-valued mapping
		$$
		X\ni h\mapsto DF(\xb;h)=\{v\in Y:\; (h,v)\in T(\gr F,(\xb,F(\xb)))\}.
		$$
		
		The contingent tangent cone and contingent derivative were introduced by Aubin
		in
		\cite{JPA81} (see \cite{AF} for detailed comments concerning genesis of the
		concept.)

		\subsubsection{Subdifferentials, normal cones  and coderivatives.}
		{\it From now on, unless the opposite is explicitly said, all spaces are assumed
			separable.} Thanks to the separable
		reduction theorem to be proved in the next subsection such a restriction is
		justifiable in the context of regularity theory. On the other hand, it provides
		for a substantial economy of efforts, especially in the non-reflexive (or to be
		precise, non-Asplund) case.
		
		Subdifferential is among the most fundamental concepts in local variational
		analysis. Essential for the
infinite dimensional variational analysis are five types of subdifferentials: Fr\'echet subdifferential, Dini-Hadamard subdifferential (the two are sometimes called ``elementary subdifferentials"), limiting Fr\'echet subdifferential, $G$-subdifferential and
the generalized gradient. In Hilbert space
		there is one more convenient construction, ``proximal subdifferential". We shall introduce it in \S~7,

		So let $f$ be a function on $X$ which if finite at $x$. The sets
		
		\vskip 2mm
		
		\centerline{
			$\sd_H f(x) =\{ x^*\in X^*:\; \lan x^*h\ran\le f^-(x;h),\; \forall h\in X\}$}
		
		\noindent and
		
		\centerline{$\sd_Ff(x)=\{ x^*\in X^*:\; \lan x^*,h\ran\le f(x+h)-f(x) + o(\|
			h\|) \}$}
		
		\vskip 2mm
		
		\noindent are called respectively the {\it Dini-Hadamard} and {\it Fr\'echet
			subdifferential} of $f$ at $x$. The corresponding {\it limiting} subdifferential
		at $x$
		(we denote them for a time being $\sd_{LH}$ and $\sd_{LF}$)
		is defined as the collection of $x^*$ such that there is a sequence
		$(x_n,x_n^*)$
		with $x_n$ norm converging to $x$ and $x_n^*$ weak$^*$-converging to $x^*$.
The essential point in the definition of the limiting subdifferentials is that only
{\it sequential} weak$^*$-limits of elements of elementary subdifferentials are considered.
The limiting Dini-Hadamard subdifferential is basically an intermediate product in the definition of the $G$-subdifferential.
		Given a set $Q\subset X$, the {\it $G$-normal cone} to $Q$ at $x\in Q$ is
		$$
		N_G(S,x) = \bigcup_{\la\ge 0}\la\sd_{LH}d(\cdot,Q)(x).
		$$
		The {\it G-subdifferential} of $f$ at $x$ is defined as follows
		$$
		\sd_Gf(x)=\{x^*:\; (x^*,-1)\in N_G(\epi f,(x,f(x)) \}.
		$$
		The cone $N_C(Q,x)=\cl(\conv N_G(Q,x))$ is {\it Clarke's normal cone} to $Q$ at
		$x$
		and the set
		$$
		\sd_Cf(x) = \{x^*:\; (x^*,-1)\in N_C(Q,x)\}
		$$
		is the {\it subdifferential} or {\it generalized gradient of Clarke}.

		\begin{proposition}[some basic properties of subdifferentials]\label{bprop}
			The following statements hold true:

			(a) for any lsc function $\sd_Hf(x)\neq\emptyset$ on a dense subset of $\dom f$;
			
			(b) the same is true for $\sd_F$ if there is a Fr\'echet differentiable (off the
			origin) norm  in    $X$ (that is if $X$ is an Asplund space);
			
			(c) if $f$ is Lipschitz near $x$, then $\sd_Gf(x)\neq\emptyset$ and the
			set-valued mapping $x\mapsto \sd_G f(x)$ is compact-valued (see (f) below) and
			upper semicontinuous;
			
			(d) if $f$ is continuously (or strictly) differentiable at $x$, then
			$\sd f(x)=\{f'(x)\}$ for any of the mentioned subdifferentials;
			
			(e) if $f$ is convex, then all mentioned subdifferentials coincide with the
			subdifferential in the sense of convex analysis: $\sd
			f(x)=\{x^*:\;f(x+h)-f(x)\ge\lan x^*,h\ran,\; \forall\; h  \}$;
			
			(f) if $f$ is Lipschitz near $x$ with Lipschitz constant $K$, then $\| x^*\|\le
			K$
			for any $x^*\in \sd f(x)$ and any of the mentioned subdifferentials;
			
			(g) if $f$ is Lipschitz near $x$, then $\sd_{LH}f(x)=\sd_Gf(x)$ and
$\sd_Cf(x) = \cl(\conv\sd_G(x))$;
			
			(h) if $f$ is lsc and  $X$ is an  Asplund space, then $\sd_{LF}f(x)=\sd_Gf(x)$
			for any $x$;
			
			(i) if $f(x,y)=\vf(x)+\psi(y)$, then $\sd f(x,y)= \sd \vf(x)+\sd\psi(y)$, where
			$\sd$
			any of $\sd_F,\ \sd_H, \sd_G$ (but not $\sd_C$).
			
		\end{proposition}
		
		\begin{remark}{\rm
				It should be observed in connection with the proposition that
				
				$\bullet$ \ $\sd_{LH}$ has little interest for non-Lipschitz functions: it may
				be too big to contain any useful information about the function.
				
				$\bullet$ \ If $X$ is not Asplund, $\sd_{LF}f(x)$ may be identically empty even
				for a very
				simple Lipschitz function (e.g. $-\| x\|$ in $C[0,1]$). In terminology of the
				subdifferential calculus this means that $\sd_F$ {\it cannot be trusted} on
				non-Asplund
				spaces.
			}
		\end{remark}
		
			We do not need here a formal definition for the concept of a subdifferential
			trusted on a space or a class of spaces (see e.g. \cite{AI12a}).
			Loosely speaking this means that a version of the fuzzy variational principle is valid for
			the subdifferentials of lsc functions on the space.
			Just note that
			the Fr\'echet subdifferential is trusted on Asplund spaces and only on them,
			Dini-Hadamard subdifferential is trusted on G\^ateaux smooth spaces
			and the G-subdifferential and the generalized gradient are trusted on all Banach spaces.
		
		There is one more important property of subdifferentials that has not been
		mentioned in the proposition. This property is called {\it tightness} and it
		characterizes a reasonable quality of lower approximation provided by the
		subdifferential (see \cite{AI12a}).
		It turns out that the Dini-Hadamard, Fr\'echet and $G$-subdifferentials are
		tight
		but Clarke's generalized gradient is not. This determines a relatively small
		role played by generalized gradient in the regularity theory. On the other hand,
		generalized gradient typically is much easier to compute and work with.
		Moreover, convexity
		of the generalized gradient makes it the only  subdifferential that can be used
		in the critical point theory associated with the concept of ``weak slope", not
		considered here.
		
		We do not need here the general theory of subdifferentials.
		Just mention
		in connection with the property (h) in Proposition \ref{bprop} that in
		separable spaces the $G$-subdifferential is a unique subdifferential having a
		certain collection of properties
		(including tightness, (c), (e), (f) and "exact calculus" as defined in the
		proposition below). It is to be again emphasized that we assume all spaces separable.

		\begin{proposition}[basic  calculus rules]\label{basrul}
			Let $f(x) = f_1(x)+f_2(x)$, where both functions are lsc and one of them  is
			Lipschitz near $\xb$. Then the following  statements are true
			
			1. {\rm Fuzzy variational principle:}  If $f$ attains a local minimum at $\xb$,
			then there are sequences $(x_{in})$ and $(x_{in}^*)$, $i=1,2$ such that
			$x_{in}\to \xb$, $x_{in}^*\in\sd_H f_{in}(x_{in})$ and $\| x_{1n}^*+
			x_{2n}^*\|\to 0$,

			2. {\rm Fuzzy sum rule}: if $X$ is Asplund and $x^*\in\sd_Ff(\xb)$, then	
			there are sequences $(x_{in})$ and $(x_{in}^*)$, $i=1,2$ such that
			$x_{in}\to \xb$, $x_{in}^*\in\sd_H f_{in}(x_{in})$ and $\| x_{1n}^*+
			x_{2n}^*-x^*\|\to 0$.	
			
			3. {\rm Exact sum rule}: $\sd_Gf(\xb)\subset \sd_G f_1(\xb)+\sd_G f_2(\xb)$.	
		\end{proposition}

		Let $Q\subset X$ and $x\in Q$. Given a subdifferential $\sd$, the set
		$$
		N(Q,x)=\sd i_Q(x),
		$$
		always a cone, is called the {\it normal cone} to $Q$ at $x$ {\it associated
			with $\sd$}.
		It is an easy matter to see that in case of $\sd_G$ this definition coincides
		with
		the given earlier. For normal cones associated with $\sd_H$ and $\sd_F$ we use
		notation
		$N_H$ and $N_F$.
		
		Let $F: X\rra Y$ and $\yb\in F(\xb)$. Given a subdifferential $\sd$ and normal
		cone associated with $\sd$, the set-valued mapping
		$$
		y^*\mapsto D^*F\xyb(y^*)=\{x^*:\; (x^*,-y^*)\in N(\gr F,\xyb) \}
		$$
		is called the {\it coderivative} of $F$ at $\xyb$ {\it associated with $\sd$}.
		We use notation $D_H^*,\ D_F^*$ and $D_G^*$ for the coderivatives, associated
		with the mentioned subdifferentials.

There is a number of monographs and survey articles in which subdifferentials
		are studied at various levels of generality: \cite{RW} (finite dimensional
		theory),
		\cite{BZ,BM,JPP,WS} (Asplund spaces), \cite{AI12a,JPP} (general Banach spaces),
		\cite{FHC83,CLSW} (generalized gradients).
Concerning the sources of the main concepts: Clarke's subdifferential was first to appear -
it was introduced in Clarke's 1973 thesis \cite{FHC73} and in printed form first appeared
in \cite{FHC75}, it is not clear where the Fr\'echet subdifferential first appeared,
probably in \cite{BGN74},
the Dini-Hadamard subdifferential was introduced by Penot in \cite{JPP74}, the
sequential
limiting Fr\'echet subdifferential for functions on Fr\'echet smooth spaces
was introduced by Kruger in mimeographed paper
\cite{AK81} in 1981 (not in \cite{KM80} as stated in e.g. \cite{MS95,BM} and many other publications-
the definition given in \cite{KM80}
is purely topological and does not involve sequential weak$^*$-limits)
and in printed form appeared in \cite{AK85} (see \cite{AI12a} for  details). The $G$-subdifferential was
first defined in \cite{AI81c}
but its definition was later modified in \cite{AI89a}.

		\subsection{Separable reduction.}
		In this subsection $X$ and $Y$ are general Banach spaces, not necessarily
		separable.
		Recall that by $\cs(X)$ we denote the collection of separable subspaces of $X$.
		\begin{proposition}\label{inher}
			Assume that $\sur F(\xb|\yb)> r$. Then for any $L_0\subset \cs(X)$
			and $M\subset \cs(Y)$ there is an  $L\in\cs(X)$ containing $L_0$
			such that for sufficiently small $t\ge 0$
			$$
			y+rt(B_Y\cap M)\subset \cl \big(F(x+t(1+\del)(B_X\cap L))\big),
			$$
			if $\del >0$ and the pair
			$(x,y)\in(\gr F)\cap(L\times M)$ is sufficiently close to $\xyb$.
		\end{proposition}
		
		\proof
		
		Take an $\ep>0$ to guarantee that the inclusion below holds for $x\in
		B(\xb,\ep)$.
		\begin{equation}\label{6.2}
		F(x)\cap B(\yb,\ep)+  trB_Y\subset F(B(x,t)).
		\end{equation}
		We shall prove that there is a nondecreasing sequence $(L_n)$ of separable
		subspaces of $X$ such that:
		\begin{equation}\label{6.1}
		y+rt(B_Y\cap M)\subset \cl\big( F(x+t(1+\del)(B_X\cap L_{n+1}))\big),
		\end{equation}
		for all $\del >0$ and all $(x,y)\in(\gr F)\cap(L_n\times M)$ sufficiently close
		to $\xyb$.
		Then to complete the proof, it is sufficient to set $L=\cl(\cup L_n)$.
		
		Assume that we have already $L_n$ for some $n$.  Let $(x_i,y_i)$ be a dense
		countable subset of
		the intersection of $( \gr F)\cap(L_n\times M)$ with the neighborhood of $\xyb$
		in
		which (\ref{6.2})
		is guaranteed, let $(v_j)$ be a dense countable subset of $B_Y\cap M$, and let
		$(t_k)$
		be a dense countable subset of $(0,\ep)$. For any $i,j,k=1,2,\ldots$ we find
		from
		(\ref{6.2}) an $h_{ijk}\in B_X$
		such that  $y_i+rt_kv_j\in F(x_i+t_kh_{ijk})$, and let $\hat L_n$ be the
		subspace of $X$ spanned
		by the union of $L_n$ and the collection of all $h_{ijk}$.
		
		If now $(x,y)\in(\gr F)\cap (L_n\times M)$, $t\in (0,1)$, $v\in B_Y$ and
		$(x_{i_m},y_{i_m})$,
		$t_{k_m}$, $v_{j_m}$ converge  respectively to $(x,y)$, $t$  and $v$, then as
		$x_{i_m}+t_{k_m}(B_X\cap M_n)\subset x+ t(1+\del)(B_X\cap M_n)$ for sufficiently
		large $m$,
		we conclude that (\ref{6.1}) holds with $\hat L_n$ instead of $L_{n+1}$.
		\endproof

		\begin{theorem}[separable reduction of regularity \cite{AI13b}]\label{sepredreg}
			Let $X$ and $Y$ be Banach spaces.
			A set-valued mapping $F: X\rra Y$ with closed graph  is
			regular at $\xyb\in\gr F$
			if and only if  for any separable subspace $M\subset Y$ and any separable
			subspace
			$L_0\subset X$
			with $\xyb\in L_0\times M$ there exists a bigger separable subspace $ L\in
			\cs(X)$   such that
			the mapping $F_{L\times M}: L\rra M$ whose graph is the intersection of $\gr F$
			with $L\times M$
			is regular at $\xyb$.
			Moreover, if  $\sur F(\xb|\yb)> r$, we can choose $L\in\cs(X)$ and $M\in\cs(Y)$ containing respectively $\xb$ and $\yb$ to make sure that also
			$\sur F_{L\times M}(\xb|\yb)\ge r$.  Conversely, if  there is an $r>0$ such
			that for any
			separable $M_0\subset Y$ and $L_0\subset X$
			there are bigger separable subspaces $M\supset M_0$ and  $L\supset L_0$ such
			that
			$\sur F_{L\times M}(\xb|\yb)\ge r$,
			then $F$ is regular at $\xyb$ with $\sur F(\xb|\yb)\ge r$.
		\end{theorem}
		
		\proof
		So assume that $F$ is regular at $\xyb$ with $\sur F(\xb|\yb)> r$.
		Then, given $L_0$ and $M$, we can find a
		closed separable subspace $L\subset X$ containing $L_0$ such that
		(\ref{6.1}) holds for any $\del>0$,
		any $(x,y)\in(\gr F)\cap (L\times M)$ sufficiently close to $\xyb$ and any
		sufficiently small $t>0$.
		
		By the Density theorem we can drop the closure operation,
		so that $F_{L\times M}$ is indeed regular near $\xyb$ with
		$\sur F_{L\times M}(\xb|\yb)\ge (1+\del)^{-1}r$. As $\del$ can be arbitrarily
		small we get
		the desired estimate for the modulus of surjection of $F_{L\times M}$.
		
		On the other hand,
		if $F$ were not regular at $\xyb$, then we could find
		a sequence $(x_n,y_n)\in \gr F$ converging to $\xyb$ such that
		$y_n+ (t_n/n)v_n\not\in F(B(x_n,t_n))$ for some $t_n<1/n$  and $v_n\in B_Y$
		(respectively $y_n+ t_n(r-\del)v_n\not\in F(B(x_n,t_n))$
		for some $\del>0$). Clearly this carries over to any closed separable subspace
		$L\subset X$ and $M\subset Y$ containing respectively all $x_n$, all $y_n$
		and all $v_n$, so that
		no such $F_{L\times M}$ cannot be regular at $\xyb$
		(with the modulus of surjection $\ge r$) contrary to the assumption.\endproof

		\subsection{Contingent derivatives and primal regularity estimates}

		The following simple proposition establishes connection between slope of $f$ and
		its lower directional derivative.
		\begin{proposition}
			\label{sldinad}
			For any function $f$ and any x at which $f$ is finite
			$$
			|\nabla f|(x)\ge-\inf_{\| h\|=1}f^-(x;h).
			$$
		\end{proposition}
		
		\proof Take an $h$ with $\| h\|=1$. We have
		$$
		|\nabla f|(x)=\lim_{t\searrow 0}\sup_{\| u\|= 1}\frac{(f(x)-f(x+tu))^+}{t}\ge
		\limsup_{(t,u)\to (0+,h)}\frac{f(x)-f(x+tu)}{t}=-f^-(x;h)
		$$
		as claimed.\endproof
		
		The following result is now immediate from the proposition and
		Theorem \ref{secmil}.

		\begin{theorem}[tangential regularity estimate 1]\label{tancrit1}
			Let $\xyb\in\gr F$. Assume that there are neighborhoods $U$ of $\xb$ and $V$ of
			$\yb$ such that for any $y\in V$  the function
			$\psi_y$ is lower semicontinuous $U$ and $\inf_{\| h\|=1}\psi_y'(x;h)\le -r$
			for $x\in U$ and $y\in V$. Then
			\begin{equation}\label{tan}
			\sur F(\xb|\yb)\ge r.
			\end{equation}
		\end{theorem}
		\noindent (Of course a similar estimate can be obtained from Theorem
		\ref{critgen}.)
		
		\begin{theorem}[tangential regularity estimate 2]\label{tancrit3}	
			Suppose there are
			a neighborhood $U$ of $\xyb$ and two numbers $c>0$ and $\la\in[0,1)$ such that
			for
			any $(x,y)\in U\cap\gr F$
			\begin{equation}\label{5.2}
			{\rm ex}(S_Y,DF(x,y)(cB_X))\le \la,
			\end{equation}
			then
			\begin{equation}\label{5.3}
			\sur  F(\xb|\yb)\ge\frac{1-\la}{c}.
			\end{equation}
			
		\end{theorem}
		
		\proof
		Take an $(x,v)\in U\cap\gr F$ with $v\neq y$ and set $z=\| y-v\|^{-1}(y-v)$.
		By the assumption for any $\la'>\la$ there is a pair $(\tilde h,\tilde w)$ with
		$\tilde w\in DF(x,v)(\tilde h)$ such that
		$\|\tilde h\|= c$ and $\| z-\tilde w\|\le \la'$. As $(\tilde h,\tilde w)$
		belongs to the contingent cone
		to the $\gr F$ at $(x,v)$, we can find (for sufficiently small $t>0$) vectors
		$h(t)$ and $w(t)$ norm converging to $\tilde h$ and $\tilde w$ respectively and
		such that
		$v+tw(t)\in F(x+th(t))$. We have
		\begin{equation}\label{5.6}
		\begin{array}{lcl}
		\|y-(v+tw(t))\|&=& \| y-v-t\tilde w\| +o(t)\\
		&\le&\| y-v-tz\| +t\| z-\tilde w\|+o(t)\\
		&\le& \|y-v\|\big(1- \dis\frac{t}{\| y-v\|}\big) +t\la' +o(t),
		\end{array}
		\end{equation}
		so that
		$$
		\vf_y^-((x,v);(\tilde h,\tilde w))
		\le\lim_{t\to+0}\frac{\|y-t(v+w(t))\|-\|y-v\|}{t}\le -(1-\la').
		$$
		
		Take a $\xi >0$ such that $\xi(1+\la)<c$ and consider the $\xi$-norm in $X\times
		Y$,
		Then $\|(\tilde h,\tilde w)\|_{\xi}\le\max\{c,\xi(1+\la')\}=c$ (if $\la'$ is
		sufficiently close to $\la$) and we get from (\ref{5.4})
		$$
		\inf\{\vf_y^-((x,v);( h,w)): \; \|(h,w)\|_{\xi}\le 1\}
		\le\frac{1}{c}\vf_y^-((x,v);(\tilde h,\tilde w)) \le -\frac{1-\la'}{c}.
		$$
		It remains to refer to Proposition \ref{sldinad} and
		Theorem \ref{critgen}.  \endproof

		\begin{theorem}[tangential regularity estimate 3]\label{tancrit2}
			Let $X$ and $Y$ be Banach spaces, and let  $F: X\rra Y$ be
			a set-valued mapping with locally closed graph. Let finally $\yb\in F(\xb)$.
			Then
			\begin{equation}\label{5.4}
			\sur F(\xb|\yb)\ge\lim_{\ep\to 0}\inf \{ C(DF(x,y)):\; (x,y)\in(\gr F)\bigcap
			B(\xyb,\ep)\},
			\end{equation}
			or equivalently,
			$$
			\begin{array}{lcl}
			\reg F(\xb|\yb)&\le&\dis\lim_{\ep\to 0}\sup\{ \|(DF(x,y))^{-1}\|_-: y\in F(x),
			\ \| x-\xb\|<\ep,\ \| y-\yb\|<\ep\}\\
			&=&\dis\lim_{\ep\to 0} \big\{\dis\sup_{\| v\|=1}\inf\{\| h\|:\; v\in
			DF(x,y)(h)\}:\\
			& &\qquad\qquad\qquad \qquad\qquad\qquad(x,y)\in(\gr F)\bigcap
			B(\xyb,\ep)\big\}.
			\end{array}
			$$
		\end{theorem}

		\proof  We first note that $DF(x,v)(B_X)$ is a star-shaped set as it contains zero and
		$z\in DF(x,v)(h)$ implies that $\la z\in DF(x,v)(\la h)$ for $\la >0$.
		On the other hand,
		by Proposition \ref{dualb} $C(DF(x,v))> r>0$ means
		that $r B_Y\subset DF(x,v)(B_X)$. It follows that  $B_Y\subset DF(x,v)(r^{-1}B_X)$.
		If this is true for all $(x,v)\in\gr F$ close to $\xyb$,
		this in turn means that the condition of Theorem \ref{tancrit3} is satisfied with $c= r^{-1}$ and $\la =1$, whence the theorem. \endproof

		\begin{remark}{\rm
		In fact the last two theorems are equivalent. Indeed, let the conditions of Theorem \ref{tancrit3} be satisfied. Then $(1-\la)B_Y\subset DF(x,v)(cB_X)$
		for all $(x,v)\in\gr F$ close to $\xyb$ and setting $r=c^{-1}(1-\la)$
		we get $rB_Y\subset DF(x,v)(B_X)$ for the same $(x,v)$.}
		\end{remark}
		
		It follows from the proofs that the estimate provided by Theorem \ref{tancrit1} is never worse than the estimates given by the other two theorems. But it can actually be
		strictly better (unless both spaces are finite dimensional). Informally, this is
		easy to understand: the quality of approximation provided by the contingent
		derivative for a map into an infinite dimensional spaces maybe much lower than
		for a real-valued function. The following  example illustrates the phenomenon.
		
		\begin{example}\label{ex1}{\rm
				Let $X=Y$ be a separable Hilbert space, and let $(e_1,e_2,\ldots)$  an
				orthonormal basis in
				$X$. Consider the following mapping from $[0,1]$ into $X$:
				$$
				\eta(t)=\left\{
				\begin{array}{cl}  0,&{\rm if}\; t\in\{0,1\}\\
				2^{-(n+2)}e_n,&{\rm if}\; t= 2^{-n},
				\end{array}\right.
				$$
				and $\eta (\cdot)$ is linear on every segment $[2^{-(n+1)},2^{-n}],\;
				n=0,1,\ldots$. Define a mapping from the unit ball of $\ell_2$ into $\ell_2$ by
				$$
				F(x) = x-\eta(\| x\|).
				$$
				
				It is an easy matter to see that $x\mapsto \eta(\| x\|)$ is $(\sqrt
				5/4)$-Lipschitz, hence by Milyutin's perturbation theorem $F$ is
				open near the origin with the rate of surjection at least $1-(\sqrt 5/4)$.
				
				Let us look what we get applying both statements of the theorem
				for the mapping. If $\| h\|= 1$ and $t\in (2^{-(n+1)},2^{-n}]$, then
				$F(th)= th - (t/2)(e_n-e_{n+1}) - 2^{-(n+2)}(2e_{n+1}-e_n))$,
				and it is easy to see that for no sequence $(t_k)$ converging to zero
				$t_k^{-1}F(t_k)$ converge.  Hence
				the tangent cone to the graph of $F$ at zero consists
				of a single point $(0,0)$ and the first statement gives
				$\sur F(0)\ge 0$ - a trivial conclusion.
				
				Now take an $x$ with $\| x\|<1$ and a $y\neq F(x)$. We have
				$$
				\begin{array}{lcl}
				\|F(x+th)-y\|&=&\|x+th-\eta(\|x+th\|)-y\|\\
				&\le& \|x+th-\eta(\| x\|)-y\|+\|\eta(\|x+th\|)-\eta(\| x\|)\|\\
				&\le &  \|F(x)+th-y\| +(3/4)t\| h\|.
				\end{array}
				$$
				
				Taking $h= (y-F(x))/\| y-F(x)\|$, we get
				$$
				\vf_y^-(x;h)\le \lim_{t\searrow
					 0}t^{-1}\Big(\Big(1-\frac{t}{\|F(x)-y\|}\Big)\|F(x)-y\|-\|F(x)-y\|\Big)
				+\frac{\sqrt 5}{4}
				=-\frac{4-\sqrt 5}{4}
				$$
				which gives $\sur F(x)\ge 1- (\sqrt 5/4)$ for all $x$ with $\| x\|<1$.
			}
		\end{example}

		A tangential regularity estimate, similar to but somewhat weaker than that in Theorem \ref{tancrit3}
		was first obtained  by Aubin in \cite{JPA84}
		(see also \cite{AF}) under the same assumptions.  The very estimate (\ref{5.3})
		was obtained in   \cite{AI87}.
		Theorem \ref{tancrit2} was proved by Dontchev-Quincampoix-Zlateva in
		\cite{DQZ06}.
		Theorem \ref{tancrit1} seems to have been state for the first time   in \cite{CFI}. Example
		\ref{ex1} has also been borrowed from that paper.
		
		
		\subsection{Dual regularity estimates.}
		This is the part of the local regularity theory  that attracted
		main attention in the 80s and 90s. The role of  coderivatives was in the center
		of the studies. Further developments, however, that followed the discovery of the
		role of slope open gates for potentially stronger (and often easier to apply) results
		involving subdifferentials of the functions $\vf_y$ and $\psi_y$.

		\subsubsection{Neighborhood estimates}
		There is a simple connection between slopes and norms of elements of subdifferentials.
			
			\begin{proposition}[slopes and subdifferentials]\label{slsub}
				Let $f$ be  lsc,  and let an open set $U$ have nonempty intersection with $\dom	f$. Then
				for any subdifferential $\sd$
				$$
				\inf_{x\in U}d(0,\sd f(x))\le \inf_{x\in U}|\nabla f|(x).
				$$
				On the other hand, $\| x^*\|\ge |\nabla f|(x)$ if $x^*\in\sd_Ff(x)$.
			\end{proposition}

		Combining this with Theorems \ref{critgen} and \ref{slsc}, we get
		
		\begin{theorem}[subdifferential regularity estimate 1]\label{subest}
			Let $X$ and $Y$ be Banach spaces,  let $F: X\rra Y$ have a locally closed
			graph, and let $\sd$ be a subdifferential
			trusted on a class of Banach spaces containing both $X$ and $Y$.
			Then for any $\xyb\in\gr F$ and any $\xi>0$
			\begin{equation}\label{est1}
			\sur F(\xb|\yb)\ge\liminf_{{(x,v)\  \underset{{\rm Graph}F}\to\
					(\xb,\yb)}\atop{y\to\yb,\ y\neq v}}\inf\{ \| x^*\|+\xi^{-1}\| v^*\|:\;
			(x^*,y^*)\in \sd\vf_y(x,v)\}.
			\end{equation}
			and
			\begin{equation}\label{est2}
			\sur F(\xb|\yb)\ge\liminf_{{(x,y)\to \xyb}\atop{y\not\in F(x)}}
			d(0,\sd \bpsi_y(x)).
			\end{equation}
		\end{theorem}

		\begin{theorem}[subdifferential regularity estimate 2]\label{subcrit2}
			Let $\xyb\in\gr F$. Assume that there are neighborhoods $U$ of $\xb$ and $V$ of
			$\yb$ such that for any $y\in V$  the function
			$\psi_y$ is lower semicontinuous and $\| x^*\|\ge r$ if $x^*\in\sd_H\psi_y(x)$
			for all $x\in U$ and $y\in V$. Then
			\begin{equation}\label{tan}
			\sur F(\xb|\yb)\ge r.
			\end{equation}
		\end{theorem}
		
		The obvious inequality $\| x^*\|\ge -f^-(x;h)$ if $x^*\in\sd_Hf(x)$ and $\|
		h\|=1$
		shows that {\it the estimate provided by the last theorem cannot be worse that
			the
			estimate of Theorem \ref{tancrit1}.}
		
		Our next purpose is to derive coderivative estimates for regularity rates.

		\begin{theorem}[coderivative regularity estimate 1]\label{subcrit}
            Let $F: X\rra Y$ be a set-valued mapping
			with locally closed graph containing $\xyb$. 		Then
			$$
			\begin{array}{lcl}
			\sur F(\xb|\yb)&\ge&\dis\dis\lim_{\ep\to 0} \inf\{ C^*(D_H^*F(x,y)): y\in F(x),
			\ \| x-\xb\|<\ep,\ \| y-\yb\|<\ep\}\\
			&=&\dis\dis\lim_{\ep\to 0} \inf\{\| x^*\|:\; x^*\in  D_H^*F(x,y)(y^*),\; \|
			y^*\|=1,\\
			& &\qquad\qquad\qquad \qquad\qquad\qquad(x,y)\in(\gr F)\bigcap B(\xyb,\ep)\},
			\end{array}
			$$
			or equivalently,
			$$
			\begin{array}{lcl}
			\reg F(\xb|\yb)=\lip F^{-1}(\yb|\xb)&\le&\dis\lim_{\ep\to
				0}\sup\{\|D_H^*F^{-1}(x,y)\|_+: \\
			& & \qquad\qquad\qquad \quad(x,y)\in(\gr F)\bigcap B(\xyb,\ep)\}\\
			&=&\dis\lim_{\ep\to 0}\sup \{\| y^*\|:\; x^*\in D_H^*F(x,y)(y^*),\;\| x^*\|=1,\\
			& &\qquad\qquad\qquad \quad(x,y)\in(\gr F)\bigcap B(\xyb,\ep)\}.
			\end{array}
			$$
		\end{theorem}
		
		To furnish the proof we can use either any of the estimates of the preceding
		theorem
		or apply directly the slope-based results of Theorems \ref{critgen} and
		\ref{slsc} via (\ref{slsub}). We choose the second option as it actually leads
		to a shorter proof. The first approach requires to work with weak$^*$ neighborhoods to estimate subdifferential
		of a sum of functions (that inevitably appears in the course of calculation)
		which
		makes estimating norms of subgradients difficult (if possible at all).
		
		\proof We only need to show that, given $(x,w)\in\gr  F$, for any neighborhoods
		$U\subset X$ and $V\subset Y$ of $x$ and $y$
		$$
		\inf\{\| x^*\|:\; x^*\in  D^*F(u,v)(y^*),\;(u,v)\in\gr F\cap(U\times V),\; \|
		y^*\|=1\}\le m.
		$$
		if $|\nabla_{\xi}\vf_y|(x,w)<m$ for small $\xi$. Then the theorem is immediate
		from Theorem \ref{critgen} in view of  Proposition \ref{slsub}.
		
		So let $|\nabla_{\xi}\vf_y|(x,w)<m$. Take an $m'<m$ but still greater than
		$|\nabla_{\xi}\vf_y|(x,v)$ and set
		$$
		\begin{array}{lcl}
		f(u,v)&=&\vf_y(u,v) + m'\max\{\|u-x\|,\xi\|v-w\| \}\\ &=& \|v-y\|+ i_{\gr
			F}(u,v)+
		m'\max\{\|u-x\|,\xi\|v-w\| \}.
		\end{array}
		$$
		Then $f$ attains a local minimum at $(x,w)$.
		
		We thus can apply Proposition \ref{basrul}: given a $\del>0$, there are
		$v_i,\ i=0,1,2$, $u_i,\ i=1,2$ with $(u_1,v_1)\in\gr F$ and
		$v_0^*\in \sd \|\cdot\|(y-v_0)$, $(u_1^*,v_1^*)\in N(\gr F,(u_1,v_1))$ and
		$(u_2^*,v_2^*)$ with $\| u_2^*\|+\xi^{-1}\| v_2^*\|\le m'$ such that
		$$
		\| v_i-w\|<\del,\quad \| u_i-x\|<\del,\quad \| u_1^*+ u_2^*\|<\del, \quad \|
		v_0^*+v_1^*+v_2^*\|<\del.
		$$
		Take $\del<\| y- w\|$, $(1+2\del)m'<m$ and $\xi$ so small that $\xi m'<\del$.
		Then $y\neq v_0$, so that
		$\| v_0^*\|=1$, $\| x_2^*\|\le m'$ and  $\| v_2^*\|<\del$. We thus have
		$\| x_1^*\|\le m'+\del< m$ and  $|\| v_1^*\|- 1|<1+2\del$. It remains to set
		$y^*= v_1^*/\| v_1^*\|$, $x^*= x_1^*/\| v_1^*\|$ to complete the proof.\endproof

		\begin{theorem}[coderivative regularity estimate 2]\label{subcrit1}
			If in addition to the assumptions of Theorem \ref{subcrit}
			both $X$ and $Y$ are Asplund spaces, then
			$$
			\begin{array}{lcl}
			\sur F(\xb|\yb)&=&\dis\dis\lim_{\ep\to 0} \inf\{ C^*(D_F^*F(x,y)): y\in F(x),
			\ \| x-\xb\|<\ep,\ \| y-\yb\|<\ep\}\\
			&=&\dis\dis\lim_{\ep\to 0} \inf\{\| x^*\|:\; x^*\in  D_F^*F(x,y)(y^*),\; \|
			y^*\|=1,\\
			& &\qquad\qquad\qquad \qquad\qquad\qquad(x,y)\in(\gr F)\bigcap B(\xyb,\ep)\},
			\end{array}
			$$
			or equivalently,
			$$
			\begin{array}{lcl}
			\reg F(\xb|\yb)=\lip F^{-1}(\yb|\xb)&=&\dis\lim_{\ep\to
				0}\sup\{\|D_F^*F^{-1}(x,y)\|_+: \\
			& & \qquad\qquad\qquad \quad(x,y)\in(\gr F)\bigcap B(\xyb,\ep)\}\\
			&=&\dis\lim_{\ep\to 0}\sup \{\| y^*\|:\; x^*\in D_F^*F(x,y)(y^*),\;\|
			x^*\|=1,\\
			& &\qquad\qquad\qquad \quad(x,y)\in(\gr F)\bigcap B(\xyb,\ep)\}.
			\end{array}
			$$
		\end{theorem}
		
		\proof If the spaces are Asplund, then the same arguments as in the proof of the preceding theorem lead to the same conclusion with $D_H^*$ replaced by $D_F^*$.
		So we have to show that the opposite inequality holds. This however is an elementary consequence of the
		definition.
		Indeed, fix certain $(x,y)\in \gr F$ close to $\xyb$ and let
		$$
		m= \inf\{\| x^*\|:\; x^*\in  D_F^*F(x,y)(y^*),\; \| y^*\|=1\}.
		$$
If $\sur D_F^*F(\xb|\yb)$= 0 or $D_F^*F(x,y)(y^*)=\emptyset$ (in which case $m=\infty$ by the general convention), the inequality is trivial.
So we can take  a positive $r<\sur F(\xb|\yb)$ in which case we may assume that $B(y,rt)\subset F(B(x,t))$ for small $t$ and $y$ close to $\yb$, and suppose that
		$m<\infty$. Take a $x^*\in D_F^*(x,y)(y^*)$ with $\| y^*\|=1$ and  $\| x^*\|<m+\del$ for some $\del>0$. Then
		$\lan x^*,h\ran-\lan y^*,v\ran\le o(\|h\|+\| v\|)$ whenever $(x+h,y+v)\in\gr F$.
		Now take $v(t)\in B(y,rt)$ such that $\lan y^*,v(t)\ran\le - (1-t^2)\| v(t)\|$
and an $h(t)$ with $\| h(t)\|\le t$ such that $(x+th(t),y+v(t))\in \gr F$. Then
$$
-t\| x^*\| + (1-t^2)rt\le \lan x^*,h(t)\ran-\lan y^*,v(t)\ran \le o(\| h(t)\|+\| v(t)\|)= o(t)
$$
which implies that $r\le m$ and the result follows.
\endproof
		
		\begin{remark}
			{\rm Note that the just given proof  (that the inequality $\le$ holds) works in
				any space, not necessarily Asplund. In other words, the  part of the theorem
				that
				incorporates essential properties of the space (that is that the Fr\'echet subdifferential is trusted) is contained in Theorem
				\ref{subcrit}.
			}
		\end{remark}

		Comparing the last theorem with Example \ref{ex1}, we conclude that in Asplund
		spaces the coderivative estimate using Fr\'echet coderivative  can be strictly
		better than the tangential estimate provided  by Theorem \ref{tancrit2}. What
		about connection of the estimates from Theorems \ref{tancrit2} and
		\ref{subcrit}?

		\begin{proposition}[DH-coderivative vs. tangential criterion]\label{dhvstan}
			The regularity estimate involving Dini-Hadamard coderivative  (Theorem
			\ref{subcrit})
			is never worse than tangential estimate provided by Theorem \ref{tancrit2}.
		\end{proposition}
		
		\proof Indeed, by definition $D_H^*F(x,y)=(DF(x,y))^*$ and we only need to
		recall that
		$C^*(D_H^*F(x,y))\ge C(DF(x,y))$ for any $(x,y)\in\gr F$ by Theorem
		\ref{bancon}.
		\endproof

		Theorem \ref{subcrit} was proved in \cite{AI87} for subdifferentials satisfying
		a bit stronger requirements than the subdifferential of Dini-Hadamard. However a
		minor change in the proof allows to extend it to all subdifferentials trusted on
		the given
		Banach space (see e.g. \cite{AI00,AI12a} also for a proof) , in particular to
		the DH-subdifferential on any G\^ateaux smooth space.
		Likewise, Theorem \ref{subcrit1} was proved in \cite{AK88}, in a somewhat
		different form and  in terms of $\ep$-Fr\'echet subdifferential on Fr\'echet
		smooth spaces. And again,
		a minor change is needed to extend the proof to standard Fr\'echet
		subdifferentials.  Theorem \ref{subcrit1} as stated was proved  in \cite{MS96a}
		(see also \cite{BM} for a proof, for all Asplund spaces, not necessarily
		separable). This extension can be viewed  as a consequence of the Fr\'echet smooth spaces
		version of the theorem and the separable reduction theorem of Fabian-Zhivkov
		\cite{FZ85} (and actually was proved that way).  Proposition \ref{dhvstan}
		seems to have never been mentioned earlier. It sounds rather surprising with all its simplicity. It would be interesting to find an example
with a Dini-Hadamard coderivative estimate strictly better, than the tangential estimate
(or to prove that the estimates are equal).
It is still unclear whether  strict inequality is possible. The general
		consideration (the dual object cannot contain more information that its original
		predecessor) suggests that this is rather unlikely. But no proof is available for
		the moment.
It should be mentioned however that the tangential estimate is valid in all Banach spaces while the Dini-Hadamard coderivative makes sense basically in G\^ateaux smooth spaces.

		\subsubsection{Perfect regularity and linear perturbations}
		The main  inconvenience of the  regularity criteria that have been just established,
		no matter primal or dual,
		comes from the necessity to scan an entire neighborhood of the point of
		interest. Below we define what can be viewed as an ideal situation.

		\begin{definition}{\rm
				We shall say that $F$ is {\it perfectly regular} at $\xyb\in\gr F$ if }
			\begin{equation}\label{5.16}
			\sur F(\xb|\yb)= C^*(D_G^*F\xyb)=\min\{\| x^*\|:\; x^*\in D_G^*F\xyb(y^*),\ \|
			y^*\|=1\}.
			\end{equation}
		\end{definition}
		
		Later we shall come across some classes of perfectly regular mappings and
		meanwhile
		consider an important class of  additive linear perturbations of maps.
		\begin{definition}{\rm
				Given a set-valued mapping $F: X\rra Y$ and an $\xyb\in\gr F$. The {\it radius
					of regularity}
				of $F$ at $\xyb$ is the lower bound of norms of linear continuous operators
				$A: X\to Y$ such that
				$\sur (F+A)(\xb,\yb+A\xb)=0$. We shall denote it $\rad F(\xb|\yb)$.}
		\end{definition}
		
		By Milyutin's theorem $\sur F(\xb|\yb)\le \rad F(\xb|\yb)$. It turns out that
		for perfectly regular mappings the equality holds. To show this we need the
		following proposition, not very difficult to prove.

		\begin{proposition}\label{sumlin} Let $X$ and $Y$ be normed  spaces, let  $F:
			X\rra Y$
			be set-valued mapping with closed graph,
			and let $A\in\cll (X,Y)$. Assume that $F$ is regular at $\xyb\in\gr F$and set
			$G=F+A$
			(that is $G(x)=F(x)+Ax$).
			Then
			$$
			D_G^*G(\xb|\yb+A\xb)=  D_G^*F\xyb + A^*
			$$
		\end{proposition}
		
		\noindent Note that the equality is elementary in case of Dini-Hadamard or
		Fr\'echet subdifferentials.
		
		\begin{theorem}[perfect regularity and radius formula]\label{radform}
			Assume that $X$ and $Y$ are Banach spaces, $F: X\rra Y$, $\xyb\in\gr F$   and
			$F+A$ is  perfectly
			regular at
			$(\xb,\yb+A\xb)$ for any $A\in\cll(X,Y)$ of rank 1. Then
			\begin{equation}\label{5.7a}
			\sur F(\xb|\yb) = \rad F(\xb|\yb).
			\end{equation}
			Moreover, for any $\ep >0$ there is a linear operator $A_{\ep}$ of rank one
			such that
			$\| A_{\ep}\|\le \sur F(\xb|\yb)+\ep$ and  $\sur (F+A)(\xb,\yb+A\xb))=0$.
		\end{theorem}
		
		In the sequel we call (\ref{5.7a}) the {\it radius formula}.
		
		\proof  Set $r=\sur F(\xb|\yb)$. The theorem is obviously valid if $r=0$. So we
		assume that $r>0$.
		Take an $\ep >0$ and  find a $y_{\ep}^*$ and an
		$x_{\ep}^*\in D_G^*F\xyb(y_{\ep}^*)$ such that $\| y_{\ep}^*\|=1,\;
		\|x_{\ep}^*\|\le (1+\ep)r$.
		Let further
		$x_{\ep}\in X$ and $y_{\ep}\in Y$ satisfy
		\begin{equation}\label{5.8a}
		\|x_{\ep}\|= \|y_{\ep}\|=1,\quad \lan x_{\ep}^*,x_{\ep}\ran\ge (1-\ep)\|
		x_{\ep}^*\|.
		\quad \lan y_{\ep}^*,y_{\ep}\ran\ge (1-\ep).
		\end{equation}
		We use these four vectors to define an operator $A_{\ep}: X\to Y$  as follows:
		$$
		A_{\ep}x=-\frac{\lan x_{\ep}^*,x\ran}{\lan y_{\ep}^*,y_{\ep}\ran}y_{\ep}.
		$$
		Then $\| A_{\ep}\|\le\dfrac{1+\ep}{1-\ep}r$ and
		$$
		A_{\ep}^*y^*=-\frac{\lan y^*,y_{\ep}\ran}{\lan y_{\ep}^*,y_{\ep}\ran}x_{\ep}^*.
		$$
		In particular we see that  $-x_{\ep}^*= A_{\ep}^*y_{\ep}^*$.  Combining this
		with
		Proposition \ref{sumlin} we get
		$0= x_{\ep}^*-A_{\ep}^*y_{\ep}^*\in D_G^*(F+A)(\xb,\yb+A\xb)(y_{\ep}^*)$
		and therefore by the prefect regularity assumption, $\sur
		(F+A)(\xb|\yb+A\xb)=0$, that is
		$\rad F\xyb\le \|A_{\ep}\|\to r$ as $\ep\to 0$. \endproof
		
		Let $S(y,A)$ be the set of solutions of the inclusion
		\begin{equation}\label{5.10a}
		y\in F(x)+Ax,
		\end{equation}
		where $A\in\cll(X,Y)$. Let $\xb$ be a nominal solution of (\ref{5.10a})
		with $y=\yb,\ A=\overline A$.
		The question we are going to consider concerns Lipschitz stability of $S$ with
		respect to small variations of both $y$ and $A$
		around
		the nominal value $(\yb,\overline A)$ and their effect on regularity rates.

		In other words, we  are interested in finding $\lip S((\yb,\overline A)|\xb)$.
		By the equivalence theorem, this is the same as finding the modulus of
		surjection of the mapping
		$\Phi=S^{-1}$ at  $(\xb,(\yb,\overline A))$. Clearly
		$$
		\Phi(x)=\{(y,A)\in Y\times\cll(X,Y):\; y\in F(x)+A(x)\}.
		$$

		We shall consider $Y\times \cll(X,Y)$ with the norm
		$\|(y,A)\|=\nu(\|y\|,\|A\|)$,
		where $\nu$ is a norm in $\R^2$.  The dual norm is $\nu^*(\|
		y^*\|,\|\ell\|)$,
		where $\ell\in(\cll(X\times Y))^*$ and $\nu^*$ is the norm in $\R^2$ dual to
		$\nu$:
		$\nu^*(u)=\sup\{\al\xi+\beta\eta:\; \nu (\al,\beta)\le 1\}$.
		As to the space dual to $\cll(X,Y)$, we only need the simplest elements
		of the space, rank one tensors  $y^*\otimes x$ whose action on $A\in\cll(X,Y)$
		is defined by $\lan y^*\otimes x,A\ran=\lan A^*y^*,x\ran$ and
		whose norm is $\| y^*\otimes x\|=\|y^*\|\|x\|$.
		
		The following theorem gives an answer to the question.

		\begin{theorem}[\cite{AI13}]\label{stablin}
			Let $X$ and $Y$ be Banach spaces, and let $F: X\rra Y$ be a set-valued mapping
			with closed graph.
			Let $\xyb\in\gr F$ and let $\overline A\in\cll(X,Y)$ be given. Then
			$$
			\lip S((\yb,\overline A)|\xb)\le  \nu^*(1,\|\xb\|)\reg  (F+\overline
			A)(\xb|\yb).
			$$
		\end{theorem}
		
		To prove the theorem we only  need to show that
		\begin{equation}\label{5.11ab}
		\sur\Phi(\xb|(\yb, \overline A))\ge \frac{1}{\nu^*(1,\|\xb\|)}\sur (F+\overline
		A)(\xb|\yb).
		\end{equation}
		So the proof (involving some calculation) can be obtained either from Theorem
		\ref{milt} or directly from the  general regularity criterion of Theorem
		\ref{gencrit},
		
		The concepts of perfect regularity and radius of regularity were introduced
		respectively in \cite{IS08} and \cite{DLR}. Theorem \ref{radform} is a new
		result.
		A finite dimensional version of Theorem \ref{stablin}
		for a class of $F$ with convex graph was proved in \cite{CGP10a}. We shall
		discuss
		the problems considered in this subsection in more details for finite
		dimensional mappings later in  Section 8.

			\section{Finite dimensional theory.}

			In this section we concentrate on characterizations of regularity, subregularity and transversality for set-valued mappings between finite dimensional spaces.
			There are several basic differences that make the finite dimensional case especially
			rich. The first is that the subdifferential calculus is much more efficient.
			In addition certain properties different
			in the general case appear to be identical in $\R^n$.   In
			particular, for a lower semicontinuous function the Dini-Hadamard subdifferential and the
			Fr\'echet
			subdifferential are identical. Therefore the usual notation used in the
			literature for this common subdifferential
			is $\hat{\sd}$ rather than $\sd_H$ or $\sd_F$. Likewise, as the limiting
			Fr\'echet
			and the $G$-subdifferentials are also equal, it is convenient to speak simply
			about
			{\it limiting subdifferential} and denote it simply by $\sd$.
			
			The second circumstance to be mentioned is the abundance of some special classes of objects of practical importance and definite theoretical interest. Enough to mention polyhedral and semi-algebraic  sets and mappings
			(to be considered in the second part of the paper),
			semi-smooth functions, prox-regular functions and sets etc..
			We do not discuss some interesting
				and important subjects, e.g.   Kummer's
				inverse function theorem and its applications
				(well presented in the literature: much on the subjects can be found in \cite{DR,KK02}) or semismooth mappings
				(see e.g. \cite{PF13}).

			
			\subsection{Regularity.}
			
			\begin{theorem}\label{F1}
				A set-valued mapping $F: \R^n\rra \R^m$ with locally closed graph is perfectly
				regular
				near any point of its graph.	
			\end{theorem}
			
			\proof This is immediate from Theorem \ref{subcrit1}.\endproof

			\begin{theorem}\label{F2}
				The radius formula holds at any point of the graph of a
				set-valued mapping $F: \R^n\rra \R^m$ with locally closed graph. Moreover, the
				lower
				bound in the definition of the radius of regularity is attained at a linear
				operator
				$A: \R^n\to \R^m$ of rank one.	
			\end{theorem}
			
			\proof This is immediate from Theorem \ref{radform}.\endproof
			
			\begin{theorem}\label{F3}
				Let $F: \R^n\rra \R^m$ be a  set-valued mapping  with locally closed graph, and
				let
				$\xyb\in\gr F$. Then
				\begin{equation}\label{eqf1}
				\sur F(\xb|\yb)=\lim_{\ep\to 0}\inf \{ C(DF(x,y)):\; (x,y)\in(\gr F)\bigcap
				B(\xyb,\ep)\}.
				\end{equation}
			\end{theorem}
			
			\proof In view of Theorem \ref{tancrit2}, it is enough to verify that
			$C(DF(x,y))\ge r$  if $B(y,tr)\subset F(B(x,t))$ for all sufficiently small $t$
			(of course for $(x,y)\in\gr F$).
			So take a $v\in \R^m$ with $\| v\|\le r$ and let $h(t)$ be such that $\|
			h(t)\|\le 1$ and
			$y+tv\in F(x+th(t))$. If now $h$ is any limiting point of $h(t)$ as $t\to 0$,
			then
			$v\in DF(x,y)(h)$. This shows that $rB_{\R^m}\subset DF(x,y)(B_{\R^n})$.
			\endproof

			Similarly, inequality can be replaced by equality in the estimate of Lipschitz
			stability
			of solutions of the inclusion
			\begin{equation}\label{7.1}
			y\in F(x)+Ax
			\end{equation}
			with both $y$ and $A$ viewed as perturbations (cf. Theorem \ref{stablin}).
			But first we have to do some preliminary job.
			As in 5.4.2 we denote by $S(y,A)$ the set of solutions of (\ref{7.1})
			and by $\Phi$ the inverse mapping
			$$
			\Phi(x)=\{(y,A):\; y\in F(x)+Ax\}.
			$$
			
			\begin{lemma} For any $x\in X$, let
				$E(x): Y\times \cll(X,Y)\to Y$ be the linear operator defined by $E(y,\Lambda)=	y-\Lambda x$. Then,
				under the assumptions of Theorem \ref{stablin}
				$$
				\nu(1,\|x\|)C(E(x)\circ D\Phi(x,(y,A))\le C(D(F+A)(x,y)),
				$$
				whenever $y\in F(x) + Ax$.
			\end{lemma}
			
			\proof By definition $(h,v,\Lambda)\in X\times Y\times \cll(X,Y)$ belongs to
			$T(\gr \Phi,(x,y,A))$ if
			there are sequences $(h_n)\to h,\; (v_n)\to v,\ (\Lambda_n)\to \Lambda$ and
			$(t_n)\to +0$
			such that
			$$
			y+t_nv_n-(A+t_n\Lambda_n)(x+t_nh_n)\in F(x+t_nh_n)
			$$
			or
			$$
			y+t_n(v_n-\Lambda_nx + t_n\Lambda_nh_n)\in (F+A)(x+t_nh_n).
			$$
			As $t_n\|\Lambda_nh_n\|\to 0$, it follows that
			$$
			T(\gr \Phi,(x,y,A))=\{(h,v,\Lambda):\; (h,v-\Lambda x)\in T(\gr(F+A),(x,y))\}
			$$
			which amounts to
			\begin{equation}\label{5.17a}
			E(x)\circ D\Phi(x,(y,A))=  D(F+A)(x,y).
			\end{equation}
			We have (Corollary \ref{constin}) $C(E(x))\cdot C(D\Phi(x,(y,A)))\le C(D(F+A)(x,y))$.
			On the other hand $E(x)^*(y^*)= (y^*, -y^*\otimes x)$ and therefore (Proposition
			\ref{calca})
			$$
			C(E(x))=\inf_{\|y^*\|=1}\| E(x)^* y^*\| = \nu (1,\| x\|).
			$$
			This completes the proof of the lemma.\endproof

			\begin{theorem}[linear perturbations - finite dimensional case]\label{F4}
				Let $F: \R^n\rra\R^m$ be a set-valued mapping with locally closed graph,
				and let $\yb\in F(\xb)$. We consider $\R^m\times \cll(\R^n,\R^m)$ with the norm
				$\nu(\| y\|,\| A\|)$, where $\nu$ is a certain  norm in $\R^2$.
				Then, given an $A\in\cll(\R^n,\R^m)$, we have
				$$
				\lip S((\yb, \overline A)|\xb) = \nu^*(1,\|\xb\|)\reg (F+\overline A)(\xb|\yb).
				$$
			\end{theorem}

			\proof Immediate from the lemma and
			Theorem \ref{stablin}.\endproof
			
			Finally, we have to mention that {\it a continuous single-valued mapping $f:\R^n\to\R^m$ can be strongly regular only if $m=n$.} This is a simple consequence of Brouwer's invariance of domain theorem (see e.g. \cite{KK02}).

			Theorem \ref{F1} was announced by Mordukhovich  in a somewhat different form \cite{BM88} (see also \cite{BM93}). But the lower  estimate
			for the modulus of surjection (which is actually the major step in the proof) is immediate from Ioffe \cite{AI84}. Theorem \ref{F2} was proved by Dontchev-Lewis-Rockafellar in \cite{DLR} and Theorem \ref{F3} by Dontchev-Quincampoix-Zlateva \cite{DQZ06}. Theorem \ref{F4} is a slightly generalized version of already mentioned result of
			C\'anovas, G\'omez and Senent-Parra \cite{CGP10a}.

			\subsection{Subregularity and error bounds.}
			
			Let $f$ be an extended-real-valued lsc function on $\R^n$. We can associate with this function the epigraphic map
			$$
			Epi f(x) = \{\al\in\R^:\; \al\ge f(x)   \}
			$$
			Subregularity of such a mapping at a point $(\xb,\al)$ (if $\al=f(\xb)$ is
			finite) means that there is a $K>0$ such that
			$$
			d(x,[f\le\al])\le K(F(x)-\al)^+
			$$
			for all $x$ close to $\xb$. The constant $K$ in this case is usually called a {\it local error bound} for $f$ at $x$. We shall say more about error bounds in the second part of the paper.

			To characterize the subregularity property of epigraphic maps
			we define
			{\it outer limiting subdifferential} of $f$ at $x$ as follows:
			$$
			\sd^{>}f(x)=\{\dis\lim_{k\to\infty} x_k^*: \; \exists \;
			x_k\underset{f}\to x,\; f(x_k)> f(x),\; x_k^*\in \hat{\sd} f(x_k)\}.
			$$
			
			\begin{theorem}[error bounds in $\R^n$]\label{erbr}
				Let $f$ be a lower semicontinuous function on $\R^n$ that is finite at $\xb$.
				Then $K>0$ is a local error bound of $f$ at $\xb$ if  either of the following
				two equivalent conditions is satisfied:
				
				\vskip 1mm
				
				(a) \ $K\cdot \dis\lim_{\ep\to 0} \inf \{|\nabla f|(x):\;
				\| x-\xb\|<\ep, \; f(\xb)<f(x)< f(\xb) +K\ep \} \ge 1$;
				
				\vskip 1mm
				
				(b) \ $K\cdot d(0,\sd^{>} f(\xb)) \ge 1.$
				
				\vskip 1mm
				
				Thus, if $F:\R^n\rra \R^m$ has locally closed graph and $\xyb\in\gr F$, then
				$$
				\subreg F(\xb|\yb)\le [\inf\{ \| x^*\|:\;  x^*\in\sd^{>}d(\yb,F(\cdot))(\xb)
				\}]^{-1}.
				$$
				
			\end{theorem}
			
			\proof  If (a) holds, then $K$ is a local error bound by Lemma \ref{baslem}
			to be proved in the next section.
			To prove that (a)$\Rightarrow$(b),  let $x^*\in\sd^{>}f(\xb)$. This
			means that there are sequences $(x_k)$
			and $(x_k^*)$ such that $x_k\to_f \xb$, $f(x_k)>f(\xb)$,  $x_k^*\to x^*$
			and $x_k^*\in\sd f(x_k)$. Choose $\ep_k\downarrow 0$ such that
			$\| x_k-\xb\|<\ep_k$ and $f(x_k)- f(\xb)<K\ep_k$. If (a) holds, then
			$K\cdot \liminf |\nabla f|(x_k)\ge 1$. But $\| x_k^*\|\ge \nabla f|(x_k)$
			(Proposition \ref{slsub}) and (b) follows.
			
			The opposite implication (b)$\Rightarrow$(a) also follows from Proposition
			\ref{slsub}. Indeed, denote by $r$ the value of the limit in the left side of
			(a),
			take an $\ep>0$ and let $x$ satisfy the bracketed inequalities in (a)
			along with $|\nabla f|(x)<r+\ep$. This means that $f+(r+\ep)\|\cdot-x\|$
			Applying the fuzzy variational principle, we shall find
			$u$ and  $u^*\in\sd_F(u)$ such that $\|u-x\|<\ep$, $f(u)<f(\xb)+\ep/K$ and
			$\| u^*\|<r+2\ep$. This means that there is a sequence of pairs $(x_k,x_k^*)$
			such that $x_k\to_f\xb$, $x_k^*\in\sd_Ff(x_k)$ and $\limsup \| x_k^*\|\le r$. As
			(b) holds, it follows that $Kr\ge 1$.  \endproof
			
			Conditions (a) and (b) are not necessary  for $K$ to be an error bound of $f$ at $\xb$.
			\begin{example}\label{conter}{\rm
					Consider
					$$
					f(x) = \left\{\begin{array}{cl} 0,&{\rm if}\;x \le 0;\\
					x+x^2\sin x^{-1},&{\rm if}\; x>0.\end{array}\right.	
					$$			
					
					It is an easy matter to see that any $K>1$ is an error bound for $f$ at zero
					but at the same time $0\in\sd^{>} f(0)$.}
			\end{example}
			Such a pathological  situation, however, does not occur if the function is
			''not too nonconvex" near $\xb$.
			
			\begin{proposition}\label{erbr1}
				Let $f$ be a lower semicontinuous function on $\R^n$ finite at $\xb$. Suppose
				there are a $\theta>0$ and a function $r(t)=o(t)$such that
				$$
				f(u)-f(x)\ge \lan x^*,u-x\ran- r(\| u-x\|)
				$$
				for all $x,\ u$ of a neighborhood of $\xb$, provided $f(\xb)<
				f(x)<f(\xb)+\theta$ and $x^*\in\hat{\sd}(x)$. If under these conditions,
				$K>0$ is an error bound of $f$ at $\xb$, then the conditions (a) and  (b) of
				Theorem \ref{erbr} hold.	
			\end{proposition}
			
			\proof Assume the contrary. Then there are $\ep>0$ and a  sequence of pairs
			$(x_k,\ x_k^*)\in\hat{\sd}f(x_k))$ such that $x_k\to_f\xb$, $f(x_k)> f(\xb)$
			and $\| x_k^*\|\le K^{-1}-\ep$. For any $k$ take an $\xb_k\in [f\le f(\xb)]$
			closest to $x_k$. Then $\xb_k\to f(\xb)$ and by the assumption
			$$
			f(\xb_k)-f(x_k)\ge \lan x_k^*,\xb_k-x_k\ran- r(\|\xb_k-x_k\|).
			$$
			As $\| \xb_k-x_k\|\to 0$, for large $k$ we have $r(\|
			\xb_k-x_k\|)\le(\ep/2)\| \xb_k-x_k\|$. For such $k$
			$$
			f(x_k)\le f(\xb_k)+ (\| x_k^*\|+(\ep/2))\| \xb_k-x_k\|.
			$$
			It follows that
			$$
			d(x_k,[f\le f(\xb)])=\| \xb_k-x_k\|\ge\frac{1}{\| x_k^*\|+(\ep/2)}f(x_k),
			$$
			that is $(K^{-1}-(\ep/2))d(x_k,[f\le f(\xb)])\ge f(x_k)$ contrary to the
			assumption.\endproof

			The last result of this subsection contains infinitesimal characterization of
			strong subregularity.
			
			\begin{theorem}[characterization of subregularity and strong subregularity]\label{charssub}
				Let again $F:\R^n\rra \R^m$ have locally closed graph and $\xyb\in\gr F$. Then
				
				$\bullet$ \ $F$ is subregular at $\xyb\in\gr F$ if $d(0,\sd^{>}\psi_{\yb}(\xb))>0$;

				$\bullet$ \ a necessary and sufficient condition for $F$ to be strongly subregular at $\xyb$
				is that $DF\xyb$ is nonsingular, that is \ 	$C^*(DF\xyb)>0$.
			\end{theorem}
			
			\proof

			The first statement is a consequence of Theorem \ref{erbr}. To prove the second,
			assume first that $F$ is strongly subregular at $\xyb$, that is
			there is a $K>0$ such that $\| x-\xb\|\le Kd(\yb,F(x))$ for $x$
			sufficiently close to $\xb$. If $DF(\xb,\yb)$ were singular,
			Proposition \ref{dualb} would guarantee the existence of sequences $(h_k)\subset
			\R^n$ and $(v_k)\subset \R^m$
			such that $\| h_k\|=1$, $\| v_k\|\to 0$ and $\yb+t_kv_k\in F(\xb+t_kh_k)$, so
			that for large $k$
			$$
			\|\xb+t_kh_k - \xb\|= t_k> Kt_k\| v_k\|=K\|\yb+t_kv_k-\yb\|\ge
			Kd(\yb,F(\xb+t_kh_k)),
			$$
			contrary to our assumption.
			
			Let now $DF\xyb$ be nonsingular. This means that $\| v\|\ge \kappa>0$
			whenever $v\in DF\xyb(h)$ with $\| h\|=1$. It immediately follows that, say,
			$\|y-\yb\|\ge (\kappa/2)\| x-\xb\|$ whenever $y\in F(x)$ and $x$
			is sufficiently close to $\xb$ which is strong subregularity of $F$ at
			$\xyb$.\endproof
			
			Literature on local error bounds in $\R^n$ is very rich - see e.g. the monograph by Facchinei and Pang \cite{FP} that summarizes developments prior to 2003.
			Theorem \ref{erbr} and Proposition \ref{erbr1} seem to be new as stated but they are closely connected with the results of Ioffe-Outrata \cite{IO08}
			and Meng and Yang \cite{MY12} among others. The second part of
			Theorem \ref{charssub} as well as  other  results relating to strong subregularity and applications can  be found in \cite{DR} and \cite{KK02}. (In \cite{KK02} the authors use the term "locally upper Lipschitz" property. The term "strong subregularity" seem to have appeared later.)  Another sufficient condition for subregularity was suggested by
			Gfrerer \cite{HG11}. It would be interesting to understand how the two are connected. It should also be noted  that no characterization for strong subregularity in terms of coderivatives is so far known.

			\subsection{Transversality.}
			We have mentioned already that the classical concepts of transversality and
			regularity are closely connected.
			To see how the concept of transversality can be interpreted in the context of
			variational analysis, we first consider  the case of two intersecting manifolds
			in a Banach space.
			
			Let $X$ be a Banach space and $M_1$ and $M_2$  smooth
			manifolds in $X$, both containing some $\xb$. As was mentioned in Subsection 1.4,
			the manifolds are transversal at $\xb$ if
			either $\xb\not\in M_1\cap M_2$ or
			the sum of the tangent subspaces to
			the manifolds at $\xb$ is the whole of $X$:
			$T_{\xb}M_1+ T_{\xb}M_2 = X.$
			The following simple lemma is the key to interpret this in regularity
			terms in a way suitable for extensions to the settings of variational analysis.
			
			\begin{lemma}\label{translem}
				Let $L_1$ and $L_2$ be closed subspaces of a Banach space $X$ such that $L_1+L_2=X$.
				Then for any $u, v\in X$ there is $h\in X$ such that $u+h\in L_1$ and $v+h\in
				L_2$.
			\end{lemma}
			
			\proof If $u+h\in L_1$, then $h\in -u+L_1$, so if the statement were wrong, we
			would have $(v-u+L_1)\cap L_2=\emptyset$. In this case there is a nonzero $x^*$
			separating $v-u+L_1$ and $L_2$, that is such that
			$\lan x^*,x\ran=0$ for all $x\in L_2$ and $\lan x^*,v-u+x\ran\ge 0$
			for all $x\in L_1$. But this means that $x^*$ vanishes on $L_1$ as well. In
			other words, both $L_1$ and $L_2$ belong to
			the annihilator of $x^*$ and so their sum cannot be the whole of $X$.\endproof
			
			The lemma effectively says that the linear mapping $(u,v,h)\mapsto (u+h,v+h)$
			maps $L_1\times L_2\times X$ {\it onto} $X\times X$, that is this mapping is
			regular. If $\xb\in M_1\cap M_2$, then  applying the density theorem
			(Theorem \ref{dens}), we get as an immediate corollary that the
			set-valued mapping $\Phi(x)= (M_1-x)\times (M_2-x)$ from $X$ into $X\times X$ is
			regular at zero. This justifies the following definition
			
			\begin{definition}\label{vartrans}{\rm Let  $S_i\subset X,\ i=1,\ldots,k$ be
					closed subsets of $X$. We say that $S_i$
					are {\it transversal} at $\xb\in X$ if either $\xb\not\in \cap S_i$
					or $\xb\in\cap S_i$ and the set-valued mapping
					$$
					x\mapsto F(x)= (S_1-x)\times\cdots\times (S_k-x)
					$$
					from $X$ into $X^k$ is regular near $(\xb,0,\ldots,0)$. In the latter case, we
					also say that {\it $S_i$ have transversal intersection} at $\xb$.
				} 	
			\end{definition}
			
			This definition may look strange at the first glance but the following
			characterization theorem shows that it is fairly natural.
			
			\begin{theorem}\label{F5}
				Let $S_i\subset \R^n, \ i=1,\ldots,k$ and $\xb\in\cap S_i$. Then the following
				statements are equivalent
				\vskip 1mm
				
				(a) \ $S_i$ are transversal at $\xb$;
				
				\vskip 1mm
				
				(b) \ $x_i^*\in N(S_i,\xb),\; x_1^*+\cdots + x_k^*= 0\  \Rightarrow \
				x_1^*=\ldots=x_k^*=0$;
				
				\vskip 1mm
				
				(c)	\ $d(x,\dis\bigcap_{i=1}^k(S_i-x_i)\le K\max_id(x,S_i-x_i)$ if $x_i$ are
				close to zero
				and $x$ is close to $\xb$.
			\end{theorem}
			
			\proof
			It is not a difficult matter to compute the limiting coderivative of $F$: if
			$(x_1,\ldots,x_k)\in F(x)$, then
			$$
			D^*F(x|(x_1,\ldots,x_k))=\left\{\begin{array}{ll}\dis\sum_{i=1}^k x_i^*,&
			{\rm if}\; x_i^*\in N(S_i,x_i+x);\\ \emptyset,&{\rm
				otherwise}.\end{array}\right.
			$$
			Combining this with Theorem \ref{F1}, we prove equivalence (a) and (b).
			
			Furthermore, $F^{-1}(x_1,\ldots,x_k)= (S_1-x_1)\cap\cdots\cap(S_k-x_k)$, whence
			equivalence of (a) and (c).\endproof
			
			Note that implicit in (c) is the statement that the intersection of $S_i-x_i$ is
			nonempty if $x_i$ are sufficiently small.
			In case of two sets one more  convenient characterization of transversality is
			available.
			
			\begin{corollary}\label{F6}
				Two sets $S_1$ and $S_2$ both containing $\xb$ are transversal at $\xb$ if and
				only if
				the set-valued mapping $\Phi: \R^n\times\R^n\rra \R^n$:
				$$
				\Phi(x_1,x_2)=\left\{ \begin{array}{cl} x_1-x_2,&{\rm if}\; x_i\in S_i;\\
				\emptyset,&{\rm otherwise}\end{array}\right.
				$$
				is regular near $(\xb,\xb,0)$.
			\end{corollary}
			
			\proof We have
			$T(\gr\Phi,((x_1,x_2),x_1-x_2)=\{(h_1,h_2,v):\; h_i\in T(S_i,x_i), \; v=
			h_1-h_2\}$, so that
			$$
			D^*\Phi((\xb,\xb),0)(x^*)=\{ (x_1^*,x_2^*):\; x_i^*\in N(S_i,\xb)+ x^* \}.
			$$
			If we consider the max-norm $\| (x_1,x_2)\|=\max\{\| x_1\|,\| x_2\|\|\}$ in
			$\R^n\times\R^n$, then it follows from Theorem \ref{F1} that
			$\Phi$ is regular near $(\xb,\xb,0)$ if and only if
			$$
			\inf \{\| x_1^*-x^*\|+\| x_2^*+ x^*\|:\; x_i^*\in N(S_i,\xn),\; \| x^*\|=1\}>0.
			$$
			This amounts to $N(S_1,\xb)\cap(-N(S_2,\xb))=\{ 0\}$, which is exactly the
			property in the part (b) of the theorem.\endproof

			In view of the equivalence between (a) and (c) in Theorem \ref{F5}, the
			following definition looks now very natural.
			
			\begin{definition}[subtransversality]\label{subtrand}{\rm
					We shall say that closed sets $S_1,\ldots,S_k$ are {\it subtransversal
						at} $\xb\in\cap S_i$ if there is a $K>0$ such that for any $x$ close to $\xb$
					$$
					d(x,\bigcap_{i=1}^k S_i)\le K\sum_{i=1}^kd(x,S_i).
					$$		
				}
			\end{definition}
			In a similar way, it is easy to see that subtrasversality is equivalent to
			subregularity
			of the same mapping $F$ and to get a sufficient  subtransversalty condition from
			Theorem \ref{erbr}.
			In the next section we shall be able to see the key role  subtransversality
			plays
			in some problems of optimization and subdifferential calculus.

			We conclude with a brief discussion of transversality of a mapping and a set.
			\begin{theorem}\label{tran1}
				Let $F: \R^n\rra \R^m$ have locally closed graph, and let $S\subset \R^m$
				be closed. Assume that $\yb\in F(\xb)\cap S$. Then the following
				statements are equivalent:
				
				(a) the set-valued mapping $\Phi: (x,y)\mapsto (F(x)-y)\times (S-y)$
				is regular near $(\xyb,(0,0))$;
				
				(b) the sets $\gr F$ and $\R^n \times S$ have transversal intersection
				near $\xyb$;
				
				(c) $0\in D^*F\xyb(y^*)\; \&\; y^*\in N(S,\yb)\; \Rightarrow\; y^*= 0$.
			\end{theorem}
			
			\proof  Equivalence of (b) and (c) follows from Theorem \ref{F5}. To prove that
			(a) and (b) are equivalent, set
			$\Psi (x,y) = (\gr F-(x,y))\times (\R^n\times S-(x,y))$.
			If $((\xi,\mu),(\eta,\nu))\in \Psi(x,y)$, then $(\mu,\nu)\in \Phi(u,y)$
			with $u= \xi+x$. Conversely,
			if $(\mu,\nu)\in\Phi(u,y)$, then
			$(u,\mu+y)\in \gr F$ and $(w,\nu+y)\in \R^n\times S$ for any $w\in\R^n$. Then for any
			$x$, we have, setting $\xi =u-x$, $\eta= w-x$ , that
			$((\xi,\mu),(\eta,\nu))\in \Psi(x,y)$.
			
			(b) $\Rightarrow$ (a).  If (b) holds, then $\Psi$
			is regular near $(\xyb,((0,0),(0,0)))$.
			So let $((\xi,\mu),(\eta,\nu))\in \Psi(x,y)$ with
			$(x,y)$  sufficiently close to $\xyb$ and $\xi,\mu,\eta,\nu$
			sufficiently close to zeros of the corresponding spaces. Take a small
			$t>0$ and let $\| \xi'-\xi\|<t$ etc. Then by (b) there is a $K>0$ and
			$(x',y')$ such that $\| x'-x\|\le Kt$, $\| y'-y\|\le Kt$ and
			$((\xi',\mu'),(\eta',\nu'))\in \Psi(x',y')$. We have
			$$
			\xi'=u'-x',\quad \mu'\in F(u')-y',\quad \eta'= w'-x',\quad,\nu'\in S-y'
			$$
			for some $(u',v')\in\gr F$ and  $w'\in \R^n$. We have therefore
			$\| u'-u\|\le \| x'-x\|+\| \xi'-\xi\|\le (K+1)t$.
			
			Thus, whenever $(\mu,\nu)\in\Phi(u,y)$ with $(u,y)$ close to $\xyb$
			and $(\mu,\nu)$ close to $(0,0)$ and $t>0$ is sufficiently small,
			for any $\mu',\nu'\in \R^m$ that differ from $\mu,\nu$ at most by $t$, there
			is a pair $(u',y')$ within $(K+1)t$ of $(u,y)$ such that
			$\mu'\in F(u')-y'$ and $\nu'\in S-y'$, that is (a).
			
			(a) $\Rightarrow$ (b). Here the arguments are similar, actually even a bit
			shorter. Let $((\xi,\mu),(\eta,\nu))\in \Psi(x,y)$ with $(x,y)$ close
			to $\xyb$ and $(\xi,\mu),(\eta,\nu))$ close to $((0,0),(0,0))$. Then
			as we have seen,  $(\mu,\nu)\in \Phi(u,y)$ with $u= \xi+x$, also close to $\xb$.
			Let further $\|\mu'-\mu\|<t,\ \| \nu'-\nu\|<t$. If $t$ is sufficiently small,
			then by (a) we can find $u', y'$ such that
			$\| u'-u\|\le Kt,\ \| y'-y\|\le Kt$ with some positive $K$ such that
			$(\mu',\nu')\in \Phi(u,y)$. Take $x'= x$, $\xi'=\ub'-x$, $\eta'=\eta$. Then as
			is
			immediate from what was explained in the first paragraph of the proof
			$((\xi,,\mu'),(\eta',\nu'))\in \Psi(x',y')$. Thus $\Psi$ is regular
			near  $(\xyb,((0,0),(0,0)))$.\endproof
			
			The proposition justifies the following definition.
			
			\begin{definition}\label{vartrans1}{\rm
					Let  $F: \R^n\rra \R^m$ have locally closed graph, let $S\subset \R^m$ be a
					closed set, and let
					$\xyb\in\gr F$. We say that $F$ is {\it transversal to $S$ at} $\xyb$
					if either $\yb\not\in S$ or $\yb\in S$ and $\gr F$ and $\R^n\times S$ are
					transversal at $\xyb$. We say that $F$ is {\it transversal} to $S$ if it is
					transversal to $S$ at any point of the graph.
					
					Likewise, if $\yb\in F(\xb)\cap S$, we shall say that $F$ is {\it subtransversal} to $S$ and $\xyb$, provided
					$$
					d((x,y),\gr F\cap(X\times S))\le Kd((x,y),\gr F) + d(y,S))
					$$
					for $(x,y)$ of a neighborhood of $\xyb$.
				}
			\end{definition}
			
			It is almost obvious from (a) that in case $\yb\in F*(\xb)\cap S$, transversality of $F$ to $S$ at $\xyb$	implies
			regularity of the mapping $x\mapsto F(x)-S$ near $(\xb,0)$. Without going into technical details the explanation is as follows. Suppose we wish to find an
			$x$ such that $z\in F(x)-S$. By (a) there are some $(x,y)$ such that
			$(0,z)\in \gr F- (x,y)$ and $(0,0)\in \R^n\times S - (x,y)$. This means that
			$z\in F(x)-y$, on the one hand, and $y\in S$, on the other hand, as required.
			
			The converse however
			does
			not seem to be valid at least for a set-valued $F$. The situation here is
			similar to that considered in Example
			\ref{contr1}. However there the converse is also true in one important case.
			
			\begin{theorem}\label{F7}
				Assume that $F: \R^n
				\to\R^m$ is Lipschitz in a neighborhood of $\xb$ and
				$C\subset \R^n$,  $Q\subset \R^m$ are nonempty and closed. Assume further  that
				$\yb=F(\xb)\in Q$. Let finally
				$$
				\Phi(x) =\left\{\begin{array}{cl} F(x)-Q,&{\rm if}\; x\in C;\\  \emptyset,&{\rm
					otherwise.}\end{array}
				\right. ;\quad
				F_C(x) =\left\{\begin{array}{cl} F(x),&{\rm if}\; x\in C;\\  \emptyset,&{\rm
					otherwise.}\end{array}
				\right.
				$$
				Then
				$	D^*\Phi(\xb,0)(y^*)=\sd(y^*\circ F_C)(\xb)$, if $y^*\in N(Q,0)$
				and $D^*\Phi(\xb,0)(y^*)=\emptyset$ otherwise.
				Thus
				$$
				\sur\Phi(\xb|0)=\min\{\| x^*\|:\; x^*\in\sd(y^*\circ F|_C)(\xb), \; y^*\in N(Q,\yb),\; \| y^*\|=1\}.
				$$		
			\end{theorem}
			\noindent (Here of course $(y^*\circ F|_C)(x)=\infty$ if $x\not\in C$.)
			If we compare this with Theorem \ref{tran1}, we see that transversality of $F_C$ to $Q$ at $\xb$ is equivalent to regularity of $F_C-Q$ near
			$(\xb,0)$. We note also the following simple corollary of the theorem
			
			\begin{corollary}\label{CF7}
				Under the assumption of the theorem	
				$$
				D^*\Phi(\xb,0)(y^*)\subset \sd(y^*\circ F)(\xb) +N(C,F(\xb)),
				\quad {\rm if} \; y^*\in N(Q,0).	
				$$
			\end{corollary}
			
			The set-valued mapping in Definition \ref{vartrans} was introduced in \cite{AI00} where it was shown that subtransversality of a collection of sets is
			equivalent to subregularity of the mapping.  Theorem \ref{F5} was partly proved in \cite{AK06} (equivalence of (a) and (c)) and partly in \cite{LLM09}
			(equivalence of (a) and (b)). We refer to \cite{AK06} for more equivalent
			descriptions (some looking very technical) of transversality and related
			properties. The results relating to transversality of  set-valued mappings and
			sets in the image space seem to be new. The exception is Theorem \ref{F7} that
			can be extracted from Theorem 5.23 of \cite{BM}.

		\vskip 1cm
		
		
		\centerline{\bf\Large Part 2. Applications}
		
		\section{Special classes of mappings}
		If additional information on the structure of a mapping is available, it is
		often possible
		to get stronger results and/or better estimates for regularity rates and to
		develop more
		convenient mechanisms to compute or estimate the latter. In this section we
		briefly discuss how this  can be implemented for
		three important classes of mappings.
		
		\subsection{Error bounds.}

		By an {\it error bound} for $f$ (at level $\al$) on a set $U$ we mean any
		estimate for the distance
		to $[f\le\al]$ in terms of $(f(x)-\al)^+$ for $x\in U$. We shall be mainly
		interested in estimates of the form
		\begin{equation}\label{6.1.1}
		d(x,[f\le\al])\le K(f(x)-\al)^+
		\end{equation}
		(which sometimes are called {\it linear} or {\it Lipschitz} error bounds).
		
		As follows from the definition, error bounds can be viewed as rates of metric subregularity
of the set-valued mapping ${\rm Epi}f(x)= [f(x),\infty)=\{\al: \; (x,\al)\in\epi f\}$
from $X$ into $\R$.
		
		\begin{lemma}[Basic lemma on error bounds]\label{baslem}
			Let $X$ be a complete metric space, let  $U\subset X$ be an open set, and  let
			$f$ be a
			lower semi-continuous function.  Suppose  that $|\nabla f|(x)>r>0$ for any
			$u\in U\backslash [f\le 0]$.
			Then for any $\xb\in U$ such that $f (\xb)<r d(\xb,X\backslash U)$ there is a
			$\ub$ such that $f(\ub)\le 0$ and
			$d(\ub,\xb)\le r^{-1}(f(\xb))^+$.
		\end{lemma}
		
		\proof Without loss of generality, we may assume that $f$ is nonnegative: just
		take
		$f^+$ instead of $f$. So
		take  an $\xb$ as in the statement. By Ekeland's principle there is a $\ub$
		such that
		$d(\ub,\xb)\le r^{-1}f(\xb)$ and $f (x) + rd(x,\ub)>f (\ub)$ if $x\neq \ub$.
		We claim that $f(\ub)\le 0$. Indeed, otherwise, by the assumption there would be
		an $x\neq \ub$ such that $f (\ub)-f(x)\ge rd(x,\ub)$ -- a
		contradiction.\endproof
		
		For simplicity we shall speak here mainly about {\it global} error bounds,
		corresponding to $U=X$, at the zero level.  We shall denote by $K_f$ the lower
		bound of $K$ such that (\ref{6.1.1}) holds for all $x$. We also set for brevity
		$$
		S= [f\le 0],\qquad S_0 = [f=0].
		$$

		\subsubsection{Error bounds for convex functions.} We shall start with the
		simplest
		case of a convex function $f$ (extended-real-valued in general) on a Banach
		space
		$X$.
		
		\begin{theorem}\label{conver}
			Let $X$ be a Banach space and $f$ a  proper closed convex function on $X$.
			Assume that $S=[f\le 0]\neq\emptyset$. Then
			\begin{equation}\label{6.1.2}
			K_f^{-1}= \inf_{x\not\in S}\ \sup_{ \| h\|\le 1}(-f'(x;h))
			= \inf_{x\not\in S}d(0,\sd f(x))
			= \inf_{x\not\in S}\dis\sur ({\rm Epi} f)(x,f(x)).
			\end{equation}
		\end{theorem}
		\noindent Here $\sd f(x)=\{x^*: f(x+h)-f(x)\ge \lan x^*,h\ran\}$ is the convex
		subdiffential.
		
		\proof Equality of the three quantities on the right is not connected with
		regularity
		and we omit the proof. To prove the first equality, we observe that the
		inequality $K_f^{-1}\le r=\inf_{x\in[f>0]}\sup_{\| h\|\le 1}(-f'(x;h))$ is
		immediate from
		Basic Lemma because for a convex function $|\nabla f|(x) = -\inf_{\| h\|\le
			1}f'(x;h)$.
		So it remains to prove the opposite inequality for which we can assume that
		$r>0$.
		
		Take a positive  $r'$ and $\del$ such that $\del< r'<r$ and let $TU(x)$ be the
		set of pairs $(u,t)$ satisfying
		\begin{equation}\label{6.1.7}
		\|u-x\|\le t,    \qquad                    f(u)\le f(x)-r't
		\end{equation}
		By Ekeland's variational principle for any $\del>0$ there is a  $(\ub,\tb)\in
		TU(x)$
		such that $f(u)+\del \| u-\ub\|$ attains its minimum at $\ub$. Clearly $\tb>0$
		(as $f(x)>0$).
		We claim that $f(\ub)=0$. Indeed,
		if $f(\ub)> 0$, then there is an $h$ with $\| h\|=1$ such that $-f'(\ub;h)>r'$,
		that is $f(\ub+th)<f(\ub)-r't$ for some $t>0$. Set $u=\ub+th$. Then
		$f(u)<f(\ub)-\del\| u-\ub\|$
		and we get a contradiction with the definition of $\ub$.
		
		Thus $f(\ub)= 0$ which means that
		$$
		d(x,S_0)\le \| \ub-x\|\le t\le \frac{1}{r'}f(x)
		$$
		and we are done as $r'$ can be chosen arbitrarily close to $r$ and $x$ is an
		arbitrary point of $[f>0]$. \endproof

		There is another way to characterize $K_f$ in terms of  normal cones to $[f\le
		0]$.
		
		\begin{theorem}\label{errng}
			For any continuous convex function $f$ on a Banach space $X$
			\begin{equation}\label{6.1.11}
			K_f= \inf_{x\in[f=0]}\inf \{\tau> 0:\; N([f=0],x)\bigcap B_{X^*}\subset
			[0,\tau]\sd f(x)\}.
			\end{equation}
		\end{theorem}

		\subsubsection{Some general results on global error bounds.} Let us turn now to
		the general case of
		a lsc function on a complete metric space.
		
		Denote now by $K_f(\al,\beta)$ (where $\beta >\al\ge  0$) the lower bound of $K$
		such that
		
		\vskip 2mm
		
		\centerline{$d(x,[f\le \al])\le K f(x)^+$  if $\al<f(x)\le \beta$.}
		
		\vskip 2mm
		
		\noindent Clearly, $K_f=\lim_{\beta\to\infty}K_f(0,\beta)$.
		
		\begin{theorem}\label{azcorv} Let $X$ be a complete metric
			space and $f$ a lower semicontinuous function on $X$. If $[f\le
			0]\neq\emptyset$, then
			$$
			\inf_{x\in[0<f\le\beta]} |\nabla f|(x)=
			\inf_{\al\in[0,\beta)}K_f(\al,\beta)^{-1}  .
			$$
		\end{theorem}
		
		\proof  Set $r=\inf_{x\in[0<f\le\beta]} |\nabla f|(x)$.  The inequality
		$K_f(\al,\beta)^{-1}\ge r$  for $0\le\al<\beta$
		is immediate from Lemma \ref{baslem}. This proves that the left side of
		the equality cannot be
		greater than the quantity on the right.
		To prove the opposite inequality it is natural to assume that
		$K_f(\al,\beta)^{-1}\ge\xi> 0$ for all $\al\in[0,\beta)$.
		For any $x\in [f>\al]$ and any $\ep >0$ such that $f(x)-\ep>\al$ choose
		a $u=u(\ep)\in [f\le f(x)-\ep]$ such that $d(x,u)\le (1+\ep)d(x,[f\le
		f(x)-\ep])\le (1+\ep)\xi^{-1}\ep$
		and therefore $u\to x$ as $\ep\to  0$.  On the other hand, $\xi d(x,u)\le
		f(x)-f(u)$ which (as $u\neq x$)
		implies that $\xi\le |\nabla f|(x)$, whence $\xi\le |\nabla f|(x)$,  and  the
		result follows.
		\endproof
		
		As an immediate consequence we get
		\begin{corollary}\label{genin}
			Under the assumption of the theorem
			$$
			K_f^{-1}\ge \inf_{x\in[f> 0]} |\nabla f|(x).
			$$
		\end{corollary}
		
		A trivial example of a function $f$ having an isolated local minimum at a
		certain $\xb$ and such that
		$\inf f<f(\xb)$  shows that the inequality can be strict. This may happen of
		course even if
		the slope is different from zero everywhere on $[f>0]$. In this case  an
		estimate of another sort  can be obtained.  Set (for $\beta >0$)
		$$
		d_f(\beta)=\sup_{x\in [ f\le\beta]}d(x,[f\le 0])
		$$
		and define the functions
		$$
		\kappa_{f,\ep}(t)=\sup\{\frac{1}{|\nabla f|(x)}:\; |f(x)-t|<\ep\};\quad
		\kappa_f(t)=\lim_{\ep\to 0}\kappa_{f,\ep}(t).
		$$
		
		\begin{proposition}\label{intin} Let $\beta>0$. Assume that  $[f\le
			0]\neq\emptyset$ and $|\nabla f|(x)\ge r>0$
			if $x\in [0< f\le \beta]$.   Then
			$$
			d_f(\beta)\le\int_0^{\beta}\kappa_f(t)dt.
			$$
		\end{proposition}

		Following the pioneering 1952 work by Hoffmann \cite{AJH52} (to be proved later in this section), error bounds,
		both for nonconvex and, especially, convex functions have been intensively studied,
		especially during last 2-3 decades, both theoretically, in connection with
		metric regularity,
		and also in view of their role in numerical analysis, see e.g.
		\cite{CC09, FP,LP98,NY05,WS06,ZN04}. Basic lemma was proved in \cite{AI00}, its
		earlier version corresponding to $U=X$ was proved by Az\'e-Corvellec-Lucchetti
		and appeared in \cite{ACL}.
		A finite dimensional versions of
		Theorems \ref{conver} and \ref{errng}
		were proved in Lewis-Pang \cite{LP98}.
		Klatte and Li \cite{KL99}.
		 The equality $K_f^{-1}=\inf\{d(0,\sd
		f(x)):\; x\in[f> 0]\}$
		in Theorem \ref{conver} was proved by Zalinescu (see  \cite{CZ}).
		The first two equalities in the theorem can be found in \cite{AC02,AC04}
		and the third equality for polyhedral functions on $\R^n$ in \cite{MPR10}.
		Theorem \ref{errng} was proved by Zheng and Ng \cite{ZN04} and
		Theorem \ref{azcorv} by
		Az\'e and Corvellec in \cite{AC02}. The papers also contain sufficiently
		thorough bibliographic comments.
		Here we follow \cite{AI13} where  proofs of all stated  and some other results
		can be found.

		\subsection{Mappings with convex graphs.}

		\subsubsection{Convex processes.} We start with the simplest class of convex
		mappings known as
		convex processes. By definition a {\it convex process} is a set-valued mapping
		$\ca: X\rra Y$
		from one Banach space into another whose  graph is a  convex cone. A convex
		process is {\it closed}
		if its graph is a closed convex cone. The closure $\cl\ca$ of a convex process
		$\ca$ is defined by
		$\gr(\cl\ca) = \cl(\gr\ca)$.
		We shall usually work with closed convex processes.
		A convex process is {\it bounded} if there is an $r>0$ such that $\| y\|\le r\|
		x\|$
		whenever $y\in \ca(x)$.  A simplest nontrivial
		example of an unbounded  closed convex process is a densely defined closed unbounded
		linear operator, as say
		the mapping $x(\cdot)\mapsto \dot x(\cdot)$ from $C[0,1]$ into itself which
		associates
		with every continuously differentiable $x(\cdot)$ its derivative and the empty
		set with any other
		element of $C[0,1]$.
		
		According to Definition \ref{normsvm}, given a convex process $\ca: X\rra Y$,
		the {\it adjoint process}
		$\ca^*: Y^*\rra X^*$ (always closed) is defined by
		$$
		\ca^*(y^*)=\{ x^*\in X^*: \; \lan x^*,x\ran\le \lan y^*,y\ran,\; \forall\;
		(x,y)\in \gr \ca\}.
		$$
		By $\ca^{**}$ we denote a convex process
		from $X$ into $Y$ whose graph is the intersection of $-\gr(\ca^*)^*$ with
		$X\times Y$, that is
		$\ca^{**}(x)=\{y:\; -y\in(\ca^*)^*(-x)\}$. Simple separation arguments show that
		$\ca^{**}=\cl\ca$ for any convex process.

		\begin{proposition}\label{convproc}
			Let $A: X\rra Y$ be a convex process. Then $\ca(Q)$ is a convex set if so is
			$Q$ and
			for any $x_1,\ x_2\in X$
			$$
			\ca(x_1)+\ca(x_2)\subset \ca(x_1+x_2).
			$$
		\end{proposition}
		
		\begin{proposition}\label{tanconv} Let $K\subset X$ be a convex closed cone.
			Then
			for any $x\in K$ the tangent cone $T(K,x)$ is the closure of the cone generated
			by $K-x$.
			In particular $K\subset T(K,x) $.
		\end{proposition}
		\noindent The propositions are the key element in the proof of the following
		fundamental property
		of convex processes.

		\begin{theorem}[regularity moduli of a convex process]\label{regproc}
			For any closed convex process $\ca: X\rra Y$ from one Banach space into another
			$$
			C(\ca)=C^*(\ca^*)=\sur \ca(0|0)=\cont\ca(0|0).
			$$
		\end{theorem}
		\noindent Note that the left inequality is equivalent to $\|
		\ca^{-1}\|_-=\|(\ca^{-1})^*\|_+$ (cf \cite{BoL}).
		
		\proof We first observe that the right equality is a consequence of the other
		two in view of Proposition \ref{dualb}.
		The inequality $C^*(\ca^*)\ge C(\ca)$ follows from Theorem \ref{bancon}. The
		same theorem
		together with the definition of Banach constants implies that $$C^*(\ca^{**})\ge
		C^*((\ca^*)^*)\ge C(\ca^*)\ge C^*(\ca^*).$$
		But  $\ca^{**} = \ca$, as $\ca$ is closed,  so that
		$C^*(\ca^{**})=C^*(\ca))\le C(\ca)$ (see again Theorem \ref{bancon}). This
		proves the left equality.
		
		Passing to the proof of the middle equality, we first observe that by
		Proposition
		\ref{dualb} $C(\ca)={\rm contr}\ca(0|0)\ge \sur\ca(0|0)$ as the rate of
		surjection can never exceed the
		modulus of controllability. On the other hand,
		by Proposition \ref{tanconv}
		$D\ca(0,0)(h)\subset D\ca(x,y)(h)$ for all $(x,y)\in\gr \ca$ and all $h$. Hence
		by Theorem \ref{tancrit1}
		$\sur \ca(0|0)\ge C(D\ca(0,0))$. But $D\ca(0,0)(h)=\ca(h)$ as the tangent cone
		to a closed convex cone at zero coincides with the latter. Thus $\sur\ca(0|0)\ge
		C(\ca)$. \endproof

		\begin{corollary}[perfect regularity of convex processes]\label{prcp}
			Any closed convex process is perfectly regular at the origin.
		\end{corollary}
		
		Note that a convex process may be not perfectly regular outside of the origin.
		For instance, consider in the space $C[0,1]$ the mapping into itself defined by
		$A(x(\cdot))= x(\cdot) + K$
		where $K$ is the cone of nonnegative functions.

		We conclude this subsection by considering the effect of linear perturbations.
		If $\ca$
		is a convex process,
		then so is $\ca+A$ where $A$ is a linear bounded operator from $X$ into $Y$.
		Thus if $\ca$
		is closed, then
		$\ca+A$ is perfectly regular at the origin and we get as an immediate
		consequence of
		Theorem \ref{radform}
		
		\begin{theorem}[radius of regularity of a convex process]\label{convrad}
			If $\ca: X\rra Y$ is a closed convex process, then
			$$
			\rad\ca(0|0)=\sur\ca(0|0).
			$$
		\end{theorem}
	
			Convex processes were introduced by Rockafellar \cite{RTR67,RTR} as an extension of linear operators  and subsequently thoroughly studied by Robinson \cite{SMR72}, Borwein
			\cite{JMB83,JMB86a} and Lewis \cite{AL99,AL01}. In particular, \cite{SMR72}
			contains an extension to convex processes of Banach-Schauder open mapping theorem.
			Another remarkable result (which is actually a special case of Theorem 5
			in the paper) can be reformulated as follows: {\it let $X$ and $Y$ be Banach spaces, and let  $\ca: X\rra Y$ and $\ct: X\rra Y$ be closed convex processes. Then $C(\ca-\ct)\ge C(\ca)- \|\ct\|_-$}. The result equivalent to the equality $C(\ca)=C^*(\ca^*)$ (Theorem \ref{regproc}) was proved
			and further discussed in \cite{JMB83,JMB86a} and Theorem \ref{convrad} in \cite{AL99}
along with the equality of the radius and distance to infeasibility for convex processes..

		\subsubsection{ Theorem of Robinson-Ursescu.}

		\begin{theorem}[surjection modulus of a convex map]\label{qru}
			Let $X$ and $Y$ be Banach spaces, and let $F: X\rra Y$ be a set-valued mapping
			with
			convex and locally closed graph. Suppose there are $\xyb\in\gr F$,
			$\al>0$ and $\beta>0$ such that $F(B(\xb,\al))$ is dense in $B(\yb,\beta)$.
			Then
			\begin{equation}\label{con1}
			\sur F(\xb|\yb)\ge \frac{\beta}{\al}.
			\end{equation}
		\end{theorem}
		
		\proof  We can set $\xb=0$, $\yb=0$. It is clear that
		$F(t\al B_X)$ is dense in $t\beta B_Y$ for any $t\in (0,1)$.  Denote $r=\beta/\al$. We shall show
		that, given a $\ga >0$, there is an $\ep>0$ such that  $F(B(x,(1+\ga)t))$ is
		dense in $B(v,rt)$ if $\| x\|<\ep$,
		$\| v\|<\ep$ and $v\in F(x)$. The theorem then will follow from Corollary
		\ref{lodense}
		
		So take a small $\ep>0$, and  let $\| x_0\|<\ep$, $\| v_0\|<\ep$ and $v_0\in F(x_0)$.
		Let further $y\in B(v_0,rt)$ for some $t\in (0,\ep)$.
		Consider the ray emanating from $v_0$ through $y$ and let $y_1$ be the point of
		the ray with $\| y_1\|=\beta$, that is there is a $\la>0$ such that
		$$
		y= \frac{1}{1+\la}y_1 +\frac{\la}{1+\la}v_0,\quad \la\ge\frac{\beta -
		\ep}{rt}.
		$$	
		We have $\| y_1-y\|=\la\| v_0-y\|$, that is
		$$
		\la=\frac{\|y_1-y}{\| v_0-y\|}\ge \frac{\beta -\ep-rt}{rt};\qquad
		1+\la\ge\frac{\beta-\ep}{rt}
		$$
		In particular, if $\beta\ge (1+2r)\ep$, which we may assume, then
		$\la\ge 1$.

		Take a $\del>0$. By the assumption there is an $x_1\in \al B$ such that
		$\| y_1-v_1\|<\del$ for some $v_1\in F(x_1)$. Set
		$$
		v=\frac{1}{1+\la}v_1 +\frac{\la}{1+\la}v_0,\quad x=\frac{1}{1+\la}x_1
		+\frac{\la}{1+\la}x_0
		$$
		Then $v\in F(x)$ as $\gr F$ is convex. We have $\| y-v\|\le
		\del/(1+\la)\le\del/2$ and
		$$
		\| x-x_0\|\le \frac{1}{1+\la}\| x_1-x_0\|\le \frac{\al+\ep}{1+\la}\le
		\frac{\al+\ep}{\beta-\ep}rt.
		$$
		If $$ 1+\ga\ge \frac{\al+\ep}{\beta-\ep}\cdot\frac{\beta}{\al},$$ this completes
		the proof as $\del$ can be chosen arbitrary small.\endproof

		As a corollary we get
		\begin{theorem}[Robinson-Ursescu \cite{SMR76b,CU75}]\label{roburs}
			Let $X$ and $Y$ be  Banach spaces. If the graph of $F: X\rra Y$ is convex  and
			closed and
			$\yb\in\intr F(X)$, then $F$ is regular at any $\xyb\in\gr F$.
		\end{theorem}
		
		\proof Let $\yb\in F(\xb)$. We have to show that there are $\al>0$
		and $\beta>0$ such that $F(B(\xb,\al))$ is dense in $B(\yb,\beta)$
		which is easy to do with the help of the standard argument using Baire category.
		\endproof

		\subsubsection{Mappings with convex graphs.  Regularity rates.} Here we give two
		results containing  exact formulas for the rate of surjection of  set-valued
		mappings with
		convex graph.

		\begin{theorem}\label{modual}
			Let $F: X\rra Y$ be a set-valued mapping with convex and locally closed graph.
			If $\yb\in F(\xb)$,
			then
			$$
			\sur F(\xb|\yb)=\lim_{\ep\to +0}\ \inf_{\| y^*\|=1}\  \inf_{x^*}\Big(\| x^*\|
			+\frac{1}{\ep}S_{\gr (F-\xyb)}(x^*,y^*)\Big).
			$$
		\end{theorem}

		The theorem  was proved in Ioffe-Sekiguchi \cite{IS08}, see also for \cite{AI13}
		for a short proof. It allows
		to also get a "primal" representation for the rate of surjection
		of a convex set-valued mapping. The key to this development is the concept of
		{\it homogenization}  $\cq$ of a convex set $Q\subset X$ which is the
		closed
		convex cone in $X\times\R$ generated by the set $Q\times \{1\}$. It is an easy
		matter
		to verify (if $Q$ is also closed) that $(x,t)\in\cq$ if and only if $x\in tQ$ if
		$t >0$ and $x\in Q^{\infty}$,
		the recession cone of $Q$, if $t=0$. (Recall that
		$Q^{\infty}=\{h\in Q:\; x+h\in Q,\; \forall x\in Q\}$.)

		Given a set-valued mapping $F: X\rra Y$ with convex closed graph, we  associate
		with $F$  and any $\xyb\in X\times Y$ (not necessarily in the graph of $F$)
		a convex process $\cf_{\xyb}: X\times\R\rra Y$
		whose graph is
		the homogenization  of $\gr F-\xyb$.   It is easy to see that
		$$
		\cf_{\xyb}(h,t) =\left\{\begin{array}{cl} t\big(
		F(\xb+\dfrac{h}{t})-\yb\big),&	{\rm if}\; t>0,\\
		F^{\infty}(h),&{\rm if} \; t=0,\\  \emptyset,&{\rm if} \; t<0,\end{array}\right.
		$$
		where $F^{\infty}$ is the ``horizon'' mapping of $F$ whose graph is the
		recession cone
		of $\gr F$:
		$$
		\gr F^{\infty}=\{(h,v):\; (x+h,y+v)\in\gr F,\; \forall (x,y)\in\gr F\}.
		$$
		If $\xyb= (0,0)$, we shall simply write $\cf$ (without the subscript) and call
		this convex process the {\it homogenization} of $F$.

		In the theorem below we use the  $\ep$-norms in $X\times \R$: $\|
		(h,t)\|_{\ep}=\max\{\| x\|,\ep t\}$
		and denote by $C_{\ep}(\cf_{\xyb})$ the Banach constant of $\cf_{\xyb}$
		corresponding to this norm.

		\begin{theorem}[primal representation of the surjection modulus]\label{moprimal}
			If $F: X\rra Y$ is a set-valued mapping with convex and locally closed graph,
			then
			$$
			\sur F(\xb|\yb)= \lim_{\ep\to +0}C_{\ep}(\cf_{\xyb}).
			$$
		\end{theorem}
		
		\proof We have (setting below $h=t(x-\xb),\; v=t(y-\yb)$)
		$$
		\begin{array}{lcl}
		\gr \cf_{\xyb}^*&=&\{(x^*,y^*,\la): \lan x^*,h\ran -\lan y^*,v\ran +\la t\le
		0:\\
		 &&\qquad\qquad\qquad\qquad\qquad\qquad\qquad\qquad\forall\;(h,v,t)\in\gr\cf_{\xyb}\}\\
		&= &\{(x^*,y^*,\la): t[\lan x^*,x-\xb\ran -\lan y^*,y-\yb\ran +\la]\le 0:\\
		&&\qquad\qquad\qquad\qquad\qquad\qquad\qquad\qquad \forall \; (x,y)\in\gr F,\;
		t>0\}\\
		&=& \{(x^*,y^*,\la): s_{\gr F-\xyb}(x^*,-y^*)+\la\le 0\}.
		\end{array}
		$$
		As the support function of $\gr F-\xyb$ is nonnegative, it follows that $\la\le
		0$
		whenever $(x^*,y^*,\la)\in\gr\cf_{\xyb}$.
		The norm in $X^*\times \R$ dual to $\|\cdot\|_{\ep}$ is
		$\|(x^*,\la)\|_{\ep}=\| x^*\|+\ep^{-1}|\la|$.  Let $d_{\ep}$ stand for the
		distance
		in $X^*\times\R$ corresponding to
		this norm. Then
		$$
		\begin{array}{lcl}
		d_{\ep}(0,\cf_{\xyb}^*{\xyb}(y^*))&=& \inf\{\|x^*\|+\ep^{-1}|\la|:\;
		s_{\gr F-\xyb}(x^*,-y^*)+\la\le 0\}\\
		&=&\dis\inf_{x^*}(\| x^*\|+\ep^{-1}s_{\gr F-\xyb}(x^*,-y^*)).
		\end{array}
		$$
		It remains to compare this with Theorem \ref{modual} to see that
		$$
		\sur F(\xb|\yb)=\lim_{\ep\to +0}\inf_{\| y^*\|=1}d_{\ep}(0,\cf_{\xyb}^*(y^*))
		$$
		and then to refer to Theorem \ref{regproc}  to conclude that the quantity
		on the right
		is precisely the limit as $\ep\to 0$ of $\inf_{\|
			y^*\|=1}C_{\ep}(\cl\cf_{\xyb}(y^*))$, where
		the closure operation can be dropped because as we mentioned the norms (and
		therefore the
		Banach constants) of
		a convex process and its closure coincide.\endproof

		The concept of homogenization was introduced by  H\"ormander \cite{LH55}.
		The idea to apply homogenization for regularity estimation goes back to Robinson's
		\cite{SMR76}. His main result actually says  that $\sur F(\xb|\yb)
\ge C_1(\cf_{\xyb})$. In a somewhat different context homogenization techniques was applied by Lewis \cite{AL01} for estimating distance to infeasibility of so called conic systems. Full statement of Theorem \ref{moprimal} was proved also in \cite{IS08}. 	We have not discussed here some well developed problems
relating to regularity of maps with convex graphs, e.g.  stability under perturbations
of systems of convex inequalities - see e.g. \cite{CLMP09,AI13,SMR75} and references in the first two quoted papers.

		\subsection{Single-valued Lipschitz maps.}
		The collection of analytic tools that allow to compute and estimate
		regularity moduli of Lipschitz single-valued mappings contains at least two
		devices,
		not available in the general situation, which are a lot more convenient
		to work with than coderivatives. The first is the {\it scalarized
			coderivative} (associated with a subdifferential):
		$$
		\cd^*F(x)(y^*) = \sd(y^*\circ F)(x)
		$$
		and the other results from suitable local approximations of the mapping either
		by
		homogeneous set-valued mappings or by sets of linear operators.
		
		The following result is straightforward.
		\begin{proposition}\label{scalf} If $F: X\to Y$  is
			Lipschitz continuous near $x\in X$, then for every $y^*\in Y^*$
			\begin{equation}\label{6.14}
			\sd_F(y^*\circ F)(x)= D_F^*F(x)(y^*).
			\end{equation}
		\end{proposition}
		
		Things are more complicated with the Dini-Hadamard subdifferential. From now on
		we assume that all spaces are G\^ateaux smooth.

		\begin{definition}\label{dircom}{\rm
				A homogeneous set-valued mapping $\ca: X\rra Y$ is a
				{\it strict Hadamard prederivative} of $F: X\to Y$ at $\xb$ if
				$\|\ca\|_+<\infty$,
				and for any norm compact set $Q\subset X$
				\begin{equation}\label{dircomd}
				F(x+th)-F(x)\subset t\ca (h)+r(t,x)t\| h\|B_Y,\quad\forall\; h\in Q,
				\end{equation}
				where $r(t,x)=r(t,x,Q)\to 0$ when $x\to \xb,\ t\to +0$. If moreover
				the inclusion holds with $Q$  replaced by $B_X$
				then  $\ca$ is called {\it strict Fr\'echet
					prederivative} of $F$ at $\xb$. Clearly, for a Fr\'echet prederivative
				we can write $r(t,x)$ in the form $\rho (t,\| x-\xb\|)$.
			}
			
		\end{definition}

		There are some  canonical ways for constructing prederivatives. The first to mention
is the generalized Jacobian introduced by Clarke \cite{FHC76b} for mappings in the finite dimensional
case and then extended to some classes of Banach spaces by P\'ales and Zeidan \cite{PZ07,PZ07a}.
Another construction, not associated with linear operators  was intruduced in \cite{AI81}. Take an $\ep >0$  and set
		$$
		\ch_{\ep}(h):=\{\la^{-1}(F(x+\la h)- F(x)): \; x, \; x + \lambda h \in  \dom F
		\cap B(\bar{x}, \ep),\; \la >0\},\ \ h\in X.
		$$
		Then $0\in \ch_{\ep}(0)$ and for $t > 0$ we have
		$$
		\ch_{\ep}(th)= t\{(t\la)^{-1}(F(x+t\la h)- F(x)): \; x, \; x + t\lambda h \in
		\dom F \cap B(\bar{x}, \ep),\; \la >0\},
		$$
		that is $\ch_{\ep}(th)=t\ch_{\ep}(h)$.
		Thus $\ch_{\ep}$ is  positively
		homogeneous and it is an easy matter to see that   (\ref{dircomd})
		holds  with  $r(t,x)=0$.

		We say that $F: X\to Y$ is {\it directionally compact} at $\xb\in\dom F$
		if it has a (norm) compact-valued  strict
		Hadamard prederivative with closed graph. It is {\it strongly
		directionally compact} if there is a compact-valued strict Fr\'echet
		prederivative
		with closed graph.
		
		The simplest, and probably the most important example of a directionally compact
		(actually even strong directionally compact) mapping is an integral operator
		associated with a differential equation, e.g.
		$$
		x(\cdot) \mapsto F(x(\cdot))(t) = x(t) - \int _0^tf(s,x(s))ds
		$$
		with $f(t,\cdot)$ Lipschitz with summable rate.
		
		\begin{proposition}[\cite{AI97}]\label{scaldh}
			If $F: X\to Y$ is Lipschitz continuous near $x$, then
			$$
			\sd_H(y^*\circ F)(x) \subset D_H^*F(x)(y^*),\quad \forall y^*\in Y^*.
			$$
			If furthermore
			$F:X\to Y$ is directionally compact at  $x$, then
			$$
			D_H^*F(x)(y^*)= \sd_H(y^*\circ F)(x)\quad\&\quad D_G^*F(x)(y^*)= \sd_G(y^*\circ
			F)(x),\quad \forall y^*\in Y^*.
			$$
		\end{proposition}

		Combining this proposition with Theorem \ref{subcrit} we get
		\begin{theorem}\label{subcrit4}
			Let $F: X\to Y$ satisfy the Lipschitz condition in a neighborhood of $\xb$.
			If $F$ is directionally compact at all $x$ of the neighborhood, then
			$$
			\sur F(\xb)\ge\lim_{\ep\to 0} \inf\{\| x^*\|:\; x^*\in  \sd_H(y^*\circ
			F)(x),\; \| y^*\|=1,\; \| x-\xb\|<\ep\},
			$$
		\end{theorem}

		The obvious inequality
		$$
		(y^*\circ F)(x+h)-(y^*\circ F)(x)\ge \inf_{w\in\ch(x)(h)}\lan y^*,w\ran
		$$
		(where $\ch(x)$ is a strict prederivative at $x$) leads to the estimate $\sur F(\xb)\ge
		\dis\liminf_{x\to\xb}C^*(\ch(x))$  under the assumptions of the theorem. A better result can be proved with the help of the general metric regularity criteria
		if $F$ has a strict Fr\'echet prederivative at $\xb$.

		\begin{theorem}\label{t1}
			Assume that $Y$ is G\^ateaux smooth and $F:X\to Y$  satisfies the Lipschitz
			condition in a neighborhood of \
			$\xb$ and, moreover, admits at $\xb$ a strict Fr\'echet prederivative   $\ch$
			with norm compact values  	such that for any $y^*$ with $\| y^*\|=1$
			\begin{equation} \label{eqBasic1}
			\sup_{\| h\|=1} \inf_{w \in \ch(h)} \langle y^*,w\rangle \ge  \rho > 0.
			\end{equation}
			\noindent Then ${\rm sur} F(\xb) \geq \rho$.
		\end{theorem}
		
		\proof
		With no loss of generality we may assume that the norm in $Y$ is G\^ateaux
		smooth off the origin.
		Take an $\ep \in (0, \rho/3) $ and an $r>0$  such that
		\begin{equation}\label{eqBasic2}
		F(x')-F(x)\in\ch(x)+\ep\| x'-x\|,
		\end{equation}
		if $x,x'\in B(\xb,r)$. Take an $x\in\bo(\xb,r/2)$ and
		a $y\in Y$, different from $F(x)$. Let $y^*$ denote the derivative of
		$\|\cdot\|$
		at $y - F(x)$. Then
		\begin{equation} \label{eq30}
		\lim_{t \to 0} t^{-1} \big(\|y - F(x) + tw\|-\|y -F(x)\|) = \langle
		y^*,w\rangle, \quad \mbox{for every} \quad w \in Y.
		\end{equation}
		
		By \eqref{eqBasic1}, there is an $h \in S_X$ such that
		\begin{eqnarray}\label{1betar}
		\langle y^*,w\rangle > \rho - \ep, \quad \mbox{for all} \quad w \in \ch(h).
		\end{eqnarray}
		Since the set $-\ch(h)$ is  compact and the limit in \eqref{eq30} is
		uniform with respect to $w$ from any fixed compact set, we conclude  that
		for  sufficiently small $t >0$
		$$
		\|y- F(x)- tw\| - \|y- F(x)\| + \langle y^*,tw\rangle < t \varepsilon    \quad
		\mbox{for all} \quad w \in \ch(h).
		$$
		This and \eqref{1betar} imply that
		\begin{equation}\label{15bbb}
		\|y- F(x)-tw\| <\|y- F(x)\| - \langle y^*,tw\rangle   + \varepsilon t \leq  \|y-
		F(x)\| - t(\rho - 2 \varepsilon)
		\end{equation}
		for all \ $w \in \ch(h)$.
		Let $x':=x+th$.  Then $\|x' - x \| = \|th\| = t < r/2$, hence $x' \in
		B(\xb,r)$. Since $\ch$ is positively homogeneous, we have  $\ch(x' - x) =
		\ch(th) = t \ch(h)$. Thus by \eqref{eqBasic2} there is a $w \in \ch(h)$ such
		that
		\begin{eqnarray}\label{16bb}
		\|F(x')-F(x)-tw\| \le  t\ep.
		\end{eqnarray}
		Now, we are ready for the following chain of estimates
		\begin{eqnarray*}
			\|y - F(x') \| & \leq &    \big\|F(x)-F(x') + tw\big\|+ \big\|y- F(x) -
			tw\big\|\\
			&<& \ep t+ \|y- F(x)\|  -(\rho-2\varepsilon)t \quad\quad  \hbox{(by
				(\ref{16bb}) and (\ref{15bbb}))}\\
			&=& \|y- F(x)\|-(\rho- 3\varepsilon) t = \|y- F(x)\|-(\rho-
			3\varepsilon)\|x'-x\|.
		\end{eqnarray*}
		It remains to apply the criterion of Theorem \ref{secmil}.
		\endproof
		
		A slight modification of the proof allows to get the following
		\begin{theorem}\label{t2}
			Assume that   $F:X\to Y$  satisfies the Lipschitz condition in a neighborhood
			of
			$\xb$ and, moreoover, there are  a homogeneous set-valued mapping  $\ch: X\rra
			Y$
			with norm compact values and $\beta\ge 0$ 	such that (\ref{eqBasic1}) holds and
			\begin{equation} \label{eq5}
			F(x+h) - F(x)\subset \ch(h) + \beta\| x'-x\|B_Y.
			\end{equation}
			\noindent Then ${\rm sur} F(\xb) \geq \rho-\beta$.
		\end{theorem}
		
		This theorem, in turn, allows us to look at  what happens when a Lipschitz
		mapping
		is approximated by a bunch of linear operators.
		Indeed, if  $\ct$ is a collection of linear operators from $X$ to $Y$, then
		the set-valued mapping $X\ni x\longmapsto \ch(x):=\{Tx:\; T\in\ct\}$ is of
		course
		positively homogeneous. It is an easy matter to see that  $\ch$ inherits some
		properties of $\ct$: for us it is important to observe that
		when $\ct$ is (relatively) norm compact in $\cll(X,Y)$ with the norm
		$\| T\|=\sup \{\| Tx\|:\; \| x\|\le 1\}$, then so are
		the values of $\ch$, if $\ct$ is bounded, then the values of $\ch$ are
		also bounded etc.. Thus we  come to the following conclusion.
		
		\begin{theorem}\label{t3}
			Assume that for a given $\xb\in\dom F$ there is a convex subset $\ct \subset
			\cll(X,Y)$ which is norm compact in $\cll(X,Y)$ and has the following two
			properties:

			(a) there is a $\beta>0$ such that for any $x,x'$ in a neighborhood of $\xb$
			there is a $T\in \ct$
			such that
			\begin{equation}\label{eq8}
			\|F(x)-F(x') - T(x-x')\|\le \beta\| x-x'\|;
			\end{equation} 	
			\indent	 (b) there are $\rho >0$ and  $\ep >0$ such that for any $T\in T$
			\begin{equation}\label{eq9}
			\ep\rho B_Y\subset T(\ep B_X).
			\end{equation}	
			\noindent Then $\sur F(\xb)\ge \rho-\beta$.
		\end{theorem}

		Scalarization formulas first appeared in \cite{AI84} for mappings between
		finite dimensional spaces and \cite{AK85} for mappings between Fr\'echet smooth spaces,
although scalarized coderivatives were considered already in \cite{AI81,AK81}.
The very term ``coderivative" was introduced in \cite{AI81}.
The concept of prederivative was introduced in \cite{AI81}
		and a characterization of directional compactness in \cite{AI97}, see also 	
		\cite{JT95a} for an earlier result.
		
		Theorems \ref{t1} and \ref{t2} will appear in \cite{CFI}. Theorem \ref{t3}
		was proved in \cite{CF13}. An earlier result without constraints on the domain of the mapping was proved by P\'ales in \cite{ZP97}  We also refer to \cite{CFI} for  a shorter proofs of the last theorem. Note that the convexity requirement in Theorem \ref{t3} is essential (consider, for instance, $F(x)=|x|: \R\to \R$
		and $\ct$ containing two operators $T_1(x)=x$ and $T_2(x)=-x$). Because of this requirement the estimate provided by Theorem \ref{t3} is generally less precise
		than those of Theorems \ref{subcrit4} and \ref{t1} (consider for instance
		the mapping $\R^2\to\R:\; F(x_1,x_2)= |x_1|-|x_2|$), but it can be easier
		to apply in certain cases (e.g. in the finite dimensional case when we can take the generalized Jacobian as $\ct$ - see \cite{FHC76b}).

			\subsection{Polyhedral  sets and mappings}
		 This subsection contains some elementary results concerning geometry of polyhedral sets in $\R^n$ and regularity of set-valued mappings with polyhedral graphs. Deeper problems associated with variational inequalities over convex
		 polyhedral sets will be discussed in the next section.

			\begin{definition}[polyhedral sets]\label{poly}{\rm
					A {\it convex polyhedral set} (or a {\it convex polyhedron}) $Q\subset \R^n$ is an intersection of finitely many closed linear subspaces and hyperplanes,	that is
					\begin{equation}\label{polyd}
					Q=\{x\in\R^n:\; \lan x_i^*,x\ran\le\al_i,\ i=1,\ldots,k;\; \lan x_i^*,x\ran=\al_i,\; i=k+1,\ldots,m  \}
					\end{equation}
					for some nonzero $x_i^*\in\R^n$ and  $\al_i\in\R$.
					Following \cite{DR} we shall use the term {\it polyhedral set} for finite unions of convex polyhedra. 		}	
			\end{definition}
			Clearly, any polyhedral set is closed. Also: as any linear equality can be replaced by two linear inequalities, we can represent any polyhedral set by means of a system of linear inequalities only. 	Elementary geometric argument allow to reveal one of the most fundamental property of
			polyhedral sets: {\it orthogonal projection of a polyhedral set
				is a polyhedral set}. In fact a linear image of a polyhedral set
			is polyhedral (see \cite{RTR} for this and other basic properties of polyhedral sets).

			A set-valued mapping $R^n\rra \R^m$  is (convex) polyhedral if so is its graph.
			Our primary interest in this section is to study regularity properties
			of such mappings.
			
			\begin{proposition}[local tangential representation]\label{loctan}
				Let $Q\subset \R^n$ be a polyhedral  set and $\xb\in Q$. Then there is
				an $\ep >0$ such that
				$$
				Q\cap B(\xb,\ep)= \xb+ T(Q,\xb)\cap(\ep B).
				$$
			\end{proposition}

			As an immediate consequence, we conclude that {\it regularity properties of a polyhedral set-valued mapping with closed graph at a point of the graph are fully determined by
				the corresponding properties at zero of its graphical derivative at the point.}
			
			One more useful corollary concerns normal cones of a polyhedral sets.
			
			\begin{proposition}\label{locnorm}
				Let $Q\subset \R^n$ be a polyhedral set. Then for any $\xb\in Q$ there is an $\ep>0$such that $N(Q,x)\subset N(Q,\xb)$ for any $x\in Q\cap B(\xb,\ep)$.
			\end{proposition}

			Our first result is the famous Hoffmann theorem on error bounds for a system of linear
			inequalities. Set $a=(\al_1,\ldots,\al_m)\in\R^m$ and let $Q(a)$ be defined by
			(\ref{polyd}).
			
			\begin{theorem}[Hoffmann]\label{hoffmann}
				Given $x_i^*\in\R^n$. Then there is a $K>0$ such that
				the inequality
				$$
				d(x,Q(a))\le K\Big(\sum_{i=1}^k(\lan x_i^*,x\ran-\al_i)^+ +\sum_{i=k+1}^m|\lan x_i^*,x\ran-\al_i|\Big)
				$$
				holds for all $x\in\R^n$ and all $a\in\R^m$ such that $Q(a)\neq\emptyset$.
			\end{theorem}
			
			\proof We shall apply Theorem \ref{conver}. Take an $a$ and set
			$$
			f(x) = \sum_{i=1}^k (\lan x_i^*,x\ran-\al_i)^+ +\sum_{i=k+1}^m |\lan x_i^*,x\ran -\al_i| .
			$$
			Then $Q(a)= [f\le 0]$. Set
			$$
			\begin{array}{ll}
			I_1(x)=\{i\in\{1,\ldots,k\}: \; \lan x_i^*,x\ran\le\al_i \},&
			J_+(x)=\{i\in\{1,\ldots,m\}: \; \lan x_i^*,x\ran>\al_i \};\\
			I_0(x)=\{i\in\{k+1,\ldots,m\}: \; \lan x_i^*,x\ran=\al_i \},&
			J_-(x)=\{i\in\{k+1,\ldots,m\}: \; \lan x_i^*,x\ran<\al_i \}
			\end{array}
			$$
			Then
			$$
			\sd f(x) =\sum_{i\in I_1(x)}[0,1]x_i^*+\sum_{i\in I_0(x)}[-1,1]x_i^*
			+\sum_{i\in J_+(x)} x_i^*-\sum_{i\in J_-(x)} x_i^*.
			$$
			If $x\not\in Q(\al)$, then $0\not\in\sd f(x)$ and $d(0,\sd f(x))>0$.
			
			We observe now that the dependence of $\sd f(x)$ of $x$ and $a$ is fully determined
			by the decomposition of the index set ${1,\ldots,m}$. Let $\Sigma$ be the collection
			of all decompositions of the index set into four subsets $I_1,\ I_0,\ J_+,\ J_-$
			such that $I_1\subset\{1,\ldots,k\}$, $I_0,J_-\subset \{k+1,\ldots,m\}$ and
			$$
			0\not\in \sum_{i\in I_1}[0,1]x_i^*+\sum_{i\in I_0}[-1,1]x_i^*
			+\sum_{i\in J_+} x_i^*-\sum_{i\in J_-} x_i^*.
			$$
			For any $\sigma\in\Sigma$ denote by $\ga(\sigma)$ the distance from zero to the
			set in the right-hand side of the above inclusion, and let $K$ stand for the
			upper bound of $\ga(\sigma)^{-1}$ over $\sigma\in\Sigma$. Then $K<\infty$ since $\Sigma$ is a finite set. Clearly, $K$ does not depend on either $a$ or $x$. On the other hand,
			$K\sd f(x)\ge 1$. It remains to refer to Theorem \ref{conver} to conclude the proof.
			\endproof
			
			As an immediate consequence, we get

			\begin{theorem}[regularity of convex polyhedral mappings]\label{polyreg}
				Let $F: \R^n\rra \R^m$ be a polyhedral set-valued mapping. Then
				
				(a) there is a $K>0$ such that $d(y,F(\xb))\le K\| x-\xb\|$ for any $\xb\in\dom F $ and any $(x,y)\in\gr F$;
				
				(b) there is a $K>0$ (different from that in (a)) such that
				$d(x,F^{-1}(y))\le Kd(y,F(x))$ for any $x\in\dom F$ and  $y\in F(X)$.
			\end{theorem}

			\noindent and

			\begin{theorem}[global subtransversality of convex polyhedral sets]\label{subpoly}
				Any two convex polyhedral sets $Q_1$ and $Q_2$ with nonempty intersection are
				globally subtransversal: there is a $K>0$ such that
				$$
				d(x,Q_1\cap Q_2)\le K(d(x,Q_1)+ d(x,Q_2)).
				$$
			\end{theorem}
			
			To prove Theorem \ref{polyreg} we have to apply the Hoffmann estimate to the
			graph of $F$. Concerning Theorem \ref{subpoly}, it should be observed that
			global transversality does not imply transversality at any point. As a simple example, consider the half spaces $S_1=\{x:\; \lan x^*,x\ran \ge 0\}$ and
			$S_2=\{x:\; \lan x^*,x\ran \le 0\}$ with some $x^*\neq 0$. The intersection
			of the sets is $\ker x^*\neq \emptyset$. But the inclusions $x_1-x\in S_1$
			and $x_2-x\in S_2$ imply $\lan x^*,x_1\ran\ge \lan x^*,x_2\ran$, hence (see Definition
			\ref{vartrans}) $S_1$ and $S_2$ are not transversal at points of $\ker x^*$.

			The results easily extend to all (not necessarily convex) polyhedral mappings.

			\begin{theorem}[subregularity of polyhedral mappings]\label{semireg}
				Let $F: \R^n\rra \R^m$ be a semi-linear set-valued mapping with closed graph.
				Then
				
				(a) there is a $K>0$ such that
				for any  $\xb\in\dom F$  there is an $\ep >0$ such that $d(y,F(\xb))\le K\| x-\xb\|$ for all $(x,y)\in\gr F$ such that $\| x-\xb\|<\ep$;
				
				(b) there is a $K>0$ (different from that in (a)) such that
				for any $\xyb\in \gr F$  there is an  $\ep >0$ such that
				$d(x,F^{-1}(y))\le K d(y,F(x))$ if
				$\| x-\xb\|<K\ep$.
				Thus $F$ is subregular at any point of its graph.
			\end{theorem}
			
			\proof We have $F(x)=\cup_{i=1}^k F_i(x)$,
			where all $F_i$ are convex polyhedral set-valued mappings. By Theorem \ref{polyreg}  for any $i$ there is a $K_i$ such that $d(y,F_i(x))\le K_i\| x-\xb\|$ for any $\xb\in\dom F_i$
			and any $(x,y)\in\gr F_i$. Now fix some $\xb\in\dom F$, and
			let $I=\{i:\;\xb\in\dom F_i\}$. Choose an  $\ep>0$  so small that $d(x,\dom  F_i)>\ep$ if $i\not\in I$ and
			$\| x-\xb\|<\ep$. (Clearly, such an $\ep$ can be found as all $\dom F_i$ are polyhedral sets, hence closed.)
			If now $y\in F(x)$ and $\| x-\xb\|<\ep$, then $I(x,y)=\{i:\; y\in F_i(x)\}\subset I$.
			On the other hand, as we have seen, there are $K_i$ such that $y\in F_i(x)$
			implies that $d(y,F_i(\xb))\le K_i\| x-\xb\|$. Thus, if
			$y\in F(x)$ and $\| x-\xb\|<\ep$, then
			$$
			d(y,F(\xb))\le\max_{i\in I(x,y)}d(y,F_i(\xb))\le (\max_i K_i)\| x-\xb\|.
			$$
			This proves the first statement.
			
			To prove the second, we apply the first to $F^{-1}$ and find $K$ and $\ep$ such
			that $d(x,F^{-1}(\yb))\le K\| v-\yb\|$ if $v\in F(x)$ and $\| v-\yb\|<\ep$.
			If $d(\yb,F(x)) <\ep$, it follows that $d(x,F^{-1}(\yb))\le Kd(\yb,F(x))$.
			This inequality trivially holds if $d(\yb,F(x)) \ge\ep$ and $\| x-\xb\|\le K\ep$.
			\endproof
			
			The property in the second part of the theorem falls short of  metric regularity because it does not guarantee that the $\ep$ will be uniformly bounded away from zero if we slightly change $\yb$.
			The following simple  example 	illustrates the phenomenon.
			
			\begin{example}{\rm
					Let $X= Y= R,\ Y$, and let
					$$
					F_1(x)=\left\{\begin{array}{cl} \R_+,&{\rm if}\; x>0,\\ \R,&{\rm if}\; x=0,\\
					\emptyset,&{\rm if}\; x<0\end{array}\right.;\qquad
					F_2(x)=\left\{\begin{array}{cl} \R_-,&{\rm if}\; x<0,\\ \R,&{\rm if}\; x=0,\\
					\emptyset,&{\rm if}\; x>0\end{array}\right.
					$$
					and $F(x) =F_1(x)\cup F_2(x)$. Fix some $y>0$ and $x<0$. Then $F^{-1}(y)= H_+$ and  $d(x,F^{-1}(y))=|x|$, $d(y,F(x))= y$ so that for no $K$ the inequality
					$d(x,F^{-1}(y)\le Kd(y,F(x))$ holds in a neighborhood of $(0,0)$.
				}
			\end{example}
			
			\begin{corollary}[subtransversality of polyhedral  sets]\label{subsemi}
				Any two semi-linear sets $Q_1$ and $Q_2$ with nonempty intersection are
				subtransversal at any common point of the sets.
				$$
				d(x,Q_1\cap Q_2)\le K(d(x,Q_1)+ d(x,Q_2)).
				$$
			\end{corollary}

			To conclude, we mention that {\it for any polyhedral mapping $F: R^n\rra \R^n$
				the set of critical values (that is such $y\in\R^m$ such that $\sur F(x|y)=0$
				for some $x\in F^{-1}(y)$) is a polyhedral set of dimension smaller than $m$.}
			This will immediately follow from the semi-algebraic Sard theorem stated in the next subsection.

			\subsection{Semi-algebraic mappings, stratifications and the Sard theorem.}
			
			Most of the results of this subsection, including
			the Sard  theorem can be extended to a wide class
			of objects, so called {\it definable} sets, mappings and functions. We however
			confine ourselves here to semi-algebraic functions whose definition is much
			simpler (compare with the general definition of definability) and does not require
			any specific effort\footnote{It should be mentioned that recently  Barbet, Dambrine, Daniilidis, Rifford \cite{BDDR} proved a remarkable result containing extensions of the Sard theorem
				to some other important classes of non-smooth functions.}

			We shall concentrate basically on two topics: consequences of the general theory and studies associated with semi-algebraic
			geometry, mainly in connection with the Sard theorem.

			\subsubsection{Basic properties {\rm (see \cite{BCR,MC})}.}
			
			A semi-algebraic set in $\R^n$ is by definition a union of finitely many sets of
			solutions
			of a finite system of polynomial equalities and inequalities of $n$ variables:
			$$
			\{ x\in \R^n:\; P_i(x)=0,\  i=1,\ldots,k,\; P_i(x)<0,\ i=k+1,\ldots, m\}.
			$$
			As immediately follows from the definition, every algebraic set is
			semi-algebraic, every
			polyhedral set is semi-algebraic, unions and intersections of finite collections
			of semi-algebraic sets are again semi-algebraic. The main fact of the
			semi-algebraic geometry is the deep Tarski-Seidenberg theorem which roughly
			speaking says that
			a linear projection of a semi-algebraic set is a semi-algebraic set. This
			theorem determines
			stability of the class of semi-algebraic sets with respect to a broad variety of
			transformations.
			
			A mapping (no matter single or set-valued) is semi-algebraic if its graph is
			semi-algebraic.
			Here is a list of some basic properties of semi-algebraic sets and mappings:
			
			$\bullet$ \ the closure and interior of a semi-algebraic set is semi-algebraic;
			
			$\bullet$ \ Cartesian product of semialgebraic sets is semi-algebraic;
			
			$\bullet$ \ composition of semi-algebraic mappings is semi-algebraic;
			
			$\bullet$ \ image and preimage of a semi-algebraic set under a semi-algebraic
			mapping is
			semi-algebraic;
			
			$\bullet$ \ derivative of a (single-valued) semi-algebraic mapping is
			semi-algebraic;
			
			$\bullet$ \ the upper and lower bound of a finite collection of
			extended-real-valued semi-algebraic functions is semi-algebraic;
			
			$\bullet$ \ if we have a semi-algebraic function of two (vector) variables, then
			its upper
			or lower bound with respect to one of the variables on a semi-algebraic set is
			semi-algebraic;
			
			$\bullet$ \ if $F$ is a semi-algebraic set-valued mapping such that every $F(x)$ is a finite set, then  the number of elements in each $F(x)$ does not exceed certain finite $N$.	
			
			For us, in the context of variational analysis and, especially, regularity
			theory, the most important is that
			
			$\bullet$ \ subdifferential mapping of a semi-algebraic function or the
			coderivative mapping of a semi-algebraic map is semi-algebraic (no matter of
			which  subdifferential on $\R^n$: Fr\'echet, Dini-Hadamard, limiting or Clarke,
			we are talking about);
			
			$\bullet$ \ slope of a semi-algebraic function is a semi-algebraic function of
			the point;
			
			$\bullet$ \ rates of regularity of a semi-algebraic functions are also
			semi-algebraic functions of the point of the graph.

			\begin{definition}\label{stratd}{\rm
					A finite partition $(M_i)$ of a set $Q\subset\R^n$ is called} $C^r$-Whitney
				stratification of $Q$
				{\rm if each $M_i$ is a $C^r$-manifold and the following two properties are
					satisfied:
					
					(a) if $(x_k)\subset M_i$ converges to some $x$  belonging to another element
					$(M_j)$ of the partition, and the unit normal vectors $v_k\in N_{x_k}M_i$
					converge to some $v$, then $v\in N_x M_j$;
					
					(b) if $M_j\cap \cl M_i\neq\emptyset$, then $M_j\subset \cl M_i$. 		
				}
			\end{definition}
			\noindent Elements of partitions are usually called {\it strata}.
			The following remarkable fact is due to S. \L ojasievicz:
			\begin{theorem}[stratification theorem]\label{strat}
				Given a semi-algebraic set  $Q\subset \R^n$ and an $r \in \N$. Then $Q$ admits a
				Whitney stratification into semi-algebraic $C^r$-manifolds.		 
			\end{theorem}
			Of course, stratification is not unique. But it
			is easy to understand that maximal dimensions of the strata coincide for all
			Whitney stratifications. This observation justifies the following
			
			\begin{definition}
				{\rm The} dimension $\dim Q$
				{\rm of a semi-algebraic set $Q$ is the maximal dimension of
					the strata in Whitney stratifications of $Q$}.
			\end{definition}

			The most important consequence of the stratification theorem is a Sard-type
			theorem
			for semi-algebraic set-valued mappings,
			
			\begin{definition}\label{sardd} {\rm	
					Let $F: \R^n\rra \R^m$ be a set-valued mapping with semi-algebraic graph, and
					let
					$\sd$ stand either for the limiting or for  the Clarke subdifferential. A
					point
					$\yb\in \R^m$ is a {\it critical value} of $F$ if there is an $x\in\R^n$ such
					that
					$y\in F(x)$ and $0\in D^*F(x|y)(y^*)$ for some $y^*\neq 0$.  	 
				}
			\end{definition}
			
			\begin{theorem}[semi-algebraic Sard theorem]\label{sardsa}
				Critical values of a semi-algebraic set-valued mapping $F: \R^n\rra \R^m$ form a
				semi-algebraic set of dimension not exceeding $m-1$.
				
				In particular an extended-real valued semi-algebraic function can have at most
				finitely many critical values.
				
			\end{theorem}
			
			For the theory of semi-algebraic sets and mappings see \cite{BCR,YC}. The Sard theorem
			was first proved by Bolte-Daniilidis-Lewis \cite{BDL06} for real-valued functions and
			then by Ioffe \cite{AI08} for set-valued mappings (in both cases the theorems were stated
			for  more general classes of objects - semi-analytic functions in \cite{BDL06} and arbitrarily stratifiable maps in \cite{AI08}).

			\subsubsection{Transversality.}
			
			We are finally ready to extend transversality theory (not just the definition)
			beyond the smooth domain.
			To begin with, we observe that a direct extension of Proposition \ref{prethom}
			does not hold if $F$ is not smooth.

			\begin{example}\label{transfail}{\rm  Consider the function
					$$
					f(x,w) = |x|-|w|
					$$
					viewed as a mapping from $\R^2$ into $\R$.  This mapping is clearly
					semi-algebraic, even polyhedral.
					It is  easy  to verify that the mapping is regular at every point
					with the modulus of
					surjection identically equal to one (if we take the $\ell^{\infty}$ norm in
					$\R^2$). Furthermore
					$$
					Q = f^{-1}(0)=\{ (x,w):\; |x|=|w|\}
					$$
					and the restriction to $Q$ of the projection $(x,w)\to w$ is also a regular
					mapping with the modulus of
					surjection equal one. However, the partial mapping $x\to f(x,0)= |x|$ is not
					regular at zero.  }
			\end{example}

			However, the following statement is true.
			
			\begin{proposition}[\cite{AI11b}]\label{regpar}
				Let $F:  \R^m\times \R^k\rra \R^n$  be a semi-algebraic  set-valued mapping
				with locally closed graph, and let $\yb\in F(\xb,\pb)$. Assume that
				
				(a)  $F$ is regular at  $((\xb,\pb),\yb)$;
				
				(b) the set-valued mapping $ \R^m\times\R^n\rra\R^k$ associating the set
				$\{ p:\  y\in F(x,p)\}$ with any $(x,y)\in\R^n\times\R^n$ is regular at
				$((\xb,\yb),\pb)$;
				
				(c) there is a Whitney stratification $(M_i)$  of $\gr  F$    such that  the
				restriction of the projection
				$(x,p)\to p$ to the set  $S_i=\{ (x,p): (x,p,\yb)\in  M_i\}$, where $M_i$ is
				the stratum  containing $(\xb,\pb,\yb)$,  is regular at $(\xb,\pb)$.
				
				Then $F_{\pb}: x\mapsto F(x,\pb)$ is regular at $\xyb$.
			\end{proposition}
			
			It is now possible to state and prove a set-valued version of Theorem
			\ref{thomsm}.
			
			\begin{theorem}\label{thomset}
				Let the mapping $F: \R^n\times\R^k\rra \R^m$ with closed graph and a closed  set $S\subset\R^m$
				be both semi-algebraic. Denote by $F_p$ the set-valued mapping $x\mapsto
				F(x,p)$. If $F$ is transversal to $S$, then for all $p$, with possible exception
				of a semi-algebraic set of dimension smaller
				than $k$, $F_p$ is transversal to $S$.
			\end{theorem}
			
			\proof The theorem is trivial if $F(x,p)\cap S=\emptyset$ for all $(x,p)$, so we assume that
			$F(x,p)$ meets $S$ for some values of the arguments. Then $(0,0)$ is a regular
			value of
			the mapping $\Psi: \R^n\times \R^m\times\R^k\to \R^m$,
			$\Psi(x,y,p)= (F(x,p)-y)\times (S-y)$. Let $Q= \Psi^{-1}(0,0)$. This is a
			semi-algebraic set, so by Theorem \ref{sardsa} there is a semi-algebraic set
			$C_0\in\R^k$ such that $\dim C_0<k$ and every $p\in\R^k\backslash C_0$ is a
			regular value of the restriction $\pi|_Q$ of the projection $(x,y,p)\mapsto p$.
			
			Take an $r> N+m-k$, and let $(M_i)_{ i=1,\ldots r}$ be a $C^1$-Whitney
			stratification of  $\gr \Psi$ with all $M_i$ being semi-algebraic manifolds.
			Then for any $i$
			there is a semi-algebraic set $C_i\subset \R^k$ such that any $p\in
			\R^k\backslash C_i$
			is a regular value of $\pi|_{M_i}$. The union $C=\bigcup_{i=0}^r C_i$ is also a
			semi-algebraic set of dimension smaller than $k$ and, as we have just seen, for
			any $p\not\in C$ all of the assumptions of Proposition \ref{regpar} are
			satisfied for $\Psi$. Therefore
			$(0,0)$ is a regular value of $\Psi_p$. By Proposition \ref{tran1} this means
			that $F_p$ is transversal to $S$.\endproof

		\section{Some applications to analysis and optimization}
		In this section we give several examples illustrating the power of regularity
		theory
		as a working instrument for treating various problems in analysis and
		optimization. We do not try each time  to prove the result under the most
		general assumptions. The purpose is rather to demonstrate how regularity
		considerations help to understand and/or simplify the analysis of one or another
		phenomenon.  Again, it should be said that some interesting areas of application
		of metric regularity remain outside the scope of the paper. Just mention
		the role of regularity in numerical optimization (see e.g. \cite{DR,KK02,KK09})
		or connections with metric fixed point theory (e.g. \cite{AAGDO,DF11,DF12,AI11a,AI14}) or recent developments associated with
		tilt stability, quadratic growth etc. (e.g. \cite{AG08,AG14,DI15,DL13,KK02,PR98} ).

		\subsection{Subdifferential calculus}
		
		In each of the three calculus rules stated in Proposition \ref{basrul} we assume
		one function Lipschitz. One of the  reasons (especially important in the proof
		of the exact sum rule) is that Lipschitz functions have bounded
		subdifferentials. But what happens when both functions
		are not Lipschitz? For instance, what can be said about normal cone to an
		intersection of
		sets? As in the calculus of convex subdifferentials, we do need some
		qualification conditions to ensure the result.
		
		\begin{theorem}\label{sumtran}
			Let $X$ be a Banach space and $S_i,\ i=1,2$ are closed subsets of $X$. Let
			further
			$\xb\in S= S_1\cap S_2$. If $S_1$ and $S_2$ are subtransversal at $\xb$, then
			$$
			N_G(S,\xb)\subset N_G(S_1,\xb)+ N_G(S_2,\xb).
			$$
		\end{theorem}
		
		Explicitly, this theorem was first mentioned in \cite{AI00} but de facto it was
		proved
		already in \cite{AI89a} (see also \cite{IP96}, Proposition 3).
		It turns out that subtransversality is the most general of all so far
		available conditions that would guarantee the inclusion. The most popular
		subdifferential transversality condition (condition (b) of Theorem \ref{F5}) may
		be much stronger.
		
		The inclusion is among the most fundamental facts of the subdifferential
		calculus: enough
		to mention that in the majority of  publications on the subject it is used as the
		starting point for deriving all other calculus rules. Below is a sketch of the
		proof of the theorem for the finite dimensional situation.
		
		\proof We need the following elementary and/or well known facts of functions on
		and sets in $\R^n$:
		
		$\bullet$ \ $\hat N(Q,x)\cap B= \hat{\sd} d(\cdot,Q)(x)$ if $x\in Q$;
		
		$\bullet$ \ if $x^*\in\hat{\sd}d(\cdot,Q)(x)$ and $u\in Q$ is the closest to
		$x$, then $x^*\in \hat{N}(Q,u)$;
		
		$\bullet$ \ if $x\in Q$ and $f(\cdot)$ is nonnegative, equal to zero at $x$ and
		$f(u)\ge d(u,Q)$ in a neighborhood of $x$, then $\hat{\sd}
		d(\cdot,Q)(x)\subset\hat{\sd}f(x)$.
		
		Combining this with the definition of the limiting subdifferential, we conclude
		that for $Q$, $f$ and $x$ as above, $\sd d(\cdot,Q)(x)\subset \sd f(x)$ - the
		fact that is surprisingly missing from monographic publications.
		
		By the assumption there is a $K>0$ such that $d(x,S)\le K(d(x,S_1)+d(x,S_2))$,
		so applying the above to $f(x) =  K(d(x,S_1)+d(x,S_2))$ along with the exact
		calculus rule of Proposition  we conclude that
		$\sd d(\cdot,S)(\xb)\subset K(\sd(\cdot, S_1)(\xb)+ \sd(\cdot, S_1)(\xb))$ and the
		result follows.\endproof

		\subsection{Necessary conditions in constrained optimization.}
		We discuss here two ways to apply  regularity theory
		to necessary optimality conditions and then a general approach to necessary
		conditions associated with one of them. Both substantially differ from classical
		proofs
		that include linearization and separation as the major steps (see e.g.
		\cite{DM65,RVG77,IT,SMR76,SMR76c}). Verification of relevance of linearization
		is usually the central and most difficult part of the proofs. It is established
		under certain constraint qualifications which always imply and often are
		equivalent to regularity of the constraint mapping (as in case of the popular
		Mangasarian-Fromovitz and Slater qualification conditions) (see e.g. \cite{SMR76}
		where the connection with regularity was made  explicit).
		
		We refer to \cite{AK85,BM88,BM} for extensions of the classical approach
		to nondifferentiable optimization in which convex separation is replaced by
		an ``extremal principle".
		The point is however that a fuller use of regularity arguments makes the way to necessary
		conditions much shorter. To begin with we shall consider
		the problem
		\begin{equation}\label{P1}
		{\rm minimize} \; f(x),\quad {\rm  s.t.}\;  F(x)\in Q, \; x\in C
		\end{equation}
		(where $F: X\to Y$ is single-valued and $Q\subset Y$ and $C\subset X$ are closed
		sets)
		assuming for simplicity that both $X$ and $Y$ are finite dimensional although
		the
		results have been originally proved in much more general situations.

		\subsubsection{Non-covering principle.}
		
		So let $\xb\in C$ be a solution of the problem. Let $\Psi$ stand for the
		restriction to
		$C$ of the set-valued mapping $x\mapsto \{f(x)-\R_-\}\times(F(x)-Q)$ from $X$ into
		$Z=\R\times Y$. Clearly, this mapping cannot be regular near $(\xb,(f(\xb),0))
		$.
		(Indeed, if $U$ is a small neighborhood of $\xb$, then $\Psi(U)$ cannot contain
		points $(f(\xb)-\ep,0)$. It follows that the negation of any condition
		sufficient for
		regularity is a necessary  condition for $\xb$ to be a local solution in the
		problem.
		Applying Theorem \ref{F7} and Corollary \ref{CF7} we get the following result.
		
		\begin{theorem}\label{N1}
			Assume that $F: \R^n\to \R^m$ is Lipschitz in a neighborhood of $\xb$.  If $\xb$
			is a local solution of (\ref{P1}), then there is a nonzero
			pair $(\la, y^*)$ such that $\la \ge 0$, $y^*\in N(Q,\yb)$ and
			\begin{equation}\label{P2}
			0\in\sd(\la f + (y^*\circ F|_C))(\xb).
			\end{equation}
		\end{theorem}
		
		This formulation needs some comments.
		We have stated the theorem in finite dimensions for simplicity, its
		infinite dimensional version can be found e.g. in \cite{AI87}.
		Note further that a more customary
		formulation would be
		\begin{equation}\label{P3}
		0\in \sd (\la f + (y^*\circ F))(\xb) + N(C,\xb).
		\end{equation}
		\noindent This condition is usually more convenient (constraints are separated)
        but in general weaker than (\ref{P2}). It is equivalent to (\ref{P2}) if
        e.g. $C=X$ (obvious) or if both $f$ and $F$ are continuously differentiable and the constraint qualification
        \begin{equation}\label{P2b}
        0\in F'(\xb)y^*+ N_C(\xb),\quad y^*\in N_Q(F(\xb)) \;\Rightarrow\; y^*=0
        \end{equation}
        is satisfied (see e.g. \cite{RW}, Example 10.8) which means that $F|_C$ is transversal to $Q$ at $\xb$ (Proposition \ref{regpar}).
		
		Finally, we observe that the necessary condition is stated in the Lagrangian form. Again,  such condition can be substantially more
		precise than
		the "separated" condition $0\in\la\sd f(\xb) +\sd(y^*\circ F)(\xb)$ (say in the
		absence of the constraint $x\in C$) which in various forms often appears in
		literature. Both conditions are equivalent if, say $f$ is continuously differentiable.

        The ``non-covering"
		approach to necessary optimality condition was first applied  probably by Warga
		\cite{JW76} in a fairly classical setting of the standard optimal control
		problem. Warga  refers not to the Lyusternik- Graves theorem
		but to the result of Yorke \cite{JY72} which is a weakened version of the
		theorem for
		integral operators associated with ordinary differential equations. But already
		the same year the controllability - optimality dichotomy appeared as the main
		tool of
		proving necessary conditions for  nonsmooth optimal control in the papers by
		Clarke \cite{FHC76c} and Warga \cite{JW76a}.
		In the context of an abstract optimization problem a non-covering criterion
		seems to have been first applied by Dmitruk-Milyutin-Osmolowski in \cite{DMO}
		to problems with finitely many functional constraints and
		recently, to problems with mixed structure (partly smooth and partly close  to
		convex), by Avakov, Magaril-Il'yaev and
		Tikhomirov \cite{AMT13}. In the next subsection 8.3 we demonstrate the work of this
		techniques
		for an abstract relaxed optimal control problem. Theorem \ref{N1} in an infinite dimensional
        setting was obtained in \cite{AI87} with the same proof based on the non-covering criterion.

		\subsubsection{Exact penalty.}
		The immediate predecessor of the approach we are going to discuss here
		was the  idea of an
		``exact penalty" \ offered by Clarke \cite{FHC76a,FHC83}: if $f$ attains a local
		minimum on a closed set $S$ at $\xb\in S$ and satisfies the Lipschitz condition
		near $\xb$, then  $\xb$ is a point of unconstrained minimum of $g(x)=f(x)+K
		d(x,S)$
		with $K$ greater than the Lipschitz constant of $f$ near $\xb$. Clarke used a
		fairly sophisticated reduction technique to apply this idea to problems with
		functional constraints.
		The arguments however are dramatically simplified by direct invoking regularity
		considerations.
		
		Let us return to the problem (\ref{P1}), assuming as above that $F$ is
		single-valued Lipschitz $X=\R^n$, $Y=\R^m$, and set as in Theorem \ref{F7}
		$$
		\Phi(x) =\left\{\begin{array}{cl} F(x)-Q,&{\rm if}\; x\in C;\\  \emptyset,&{\rm
			otherwise.}\end{array}	\right.	
		$$
		Then our problem can be reformulated as	
		\begin{equation}\label{P2a}
		{\rm minimize} \; f(x),\quad {\rm  s.t.}\;  0\in\Phi(x).
		\end{equation}
		
		Suppose that $\Phi$ is subregular at $(\xb,0)$. This means that there is some
		$K_0>0$
		such that $d(x,\Phi^{-1}(0))\le K_0d(0,\Phi(x))$ for $x$ of a neighborhood of
		$\xb$.
		But $\Phi^{-1}(0)$ is the feasible set of our problem, so that there is
		some other $K_1>0$ such that the function
		$f(x) +K_1 d(0,\Phi(x)) $
		attains local minimum at $\xb$ or equivalently, the function
		$f(x)+K_1d(y,F(x)-Q)$ attains a local minimum at $\xb$ subject to $x\in C$.
		The last function is Lipschitz continuous near $\xb$, hence there is a $K$ such that
		\begin{equation}\label{P8}
		g(x) = f(x)+K(d(y,F(x)-Q) + d(x,C)
		\end{equation}
		attains an unconditional minimum at $\xb$.
		
		If on the other hand, $\Phi$ is nor subregular at $\xb$, Theorems \ref{F1} and
		\ref{F7} imply together that $0\in \sd(y^*\circ F)(\xb)+ N(C,\xb)$ for some
		nonzero $y^*\in N(Q,F(\xb))$.
		From here we easily get a weakened version of Theorem
		\ref{N1}
		with the Lagrangian condition replaced by its ``separated"  versions
		$$
		0\in \sd f(\xb)+ \sd(y^*\circ F)(\xb)+ N(C,\xb),\quad y^*\in N(Q,F(\xb)).
		$$
		
		This is a definite drawback, as we have already mentioned   which however is
		counterbalanced by some serious advantages.
		First we note that $g$ is defined in
		terms of the original data  which   makes it possible to study higher order
		optimality conditions using this function. This is how such a techniques was
		used for the first time in  \cite{AI79}
		in order to get necessary optimality conditions earlier obtained by
		Levitin-Milyutin-Osmolowski in  \cite{LMO}.
		
		Another advantage is that the second approach is more universal. It can work for problems for which using scalarized coderivatives is either difficult or just impossible as say, in problems involving inclusions $0\in \Phi (x)$ with general set-valued  $\Phi$. This is a typical case in optimal control of dynamic systems described by differential inclusions.   	Loewen  \cite{PDL} was the first to use this
		approach to
		prove a maximum principle in a free
		right end point problem of that sort. The analytic challenge in his proof was to find an upper estimate for the distance to the feasible set. However the next step in the development, the "optimality alternative" discussed below,		excludes even any need in such an estimate.

		\subsubsection{Optimality alternative.}
		Consider the abstract problem with $(X,d)$ being a complete metric space:
		
		\vskip 1mm
		
		\centerline{\rm minimize\quad $f(x)$, \ \ subject to \ \ $x\in Q\subset X$.}

		\begin{theorem}\label{optalt}
			Let  $\vf$ be a nonnegative lsc function on $X$ equal to zero at $\xb$. If
			$\xb\in Q$ is a local solution to the problem, then the following alternative holds
			true:
			
			\vskip 1mm
			
			\noindent$\bullet$ \ either there is a $\la >0$ such that the function $\la f+\vf$ has an unconstrained local
			minimum at $\xb$;
			
			\noindent$\bullet$ \ or there is a sequence $(x_n)\to \xb$ such that
			$\vf(x_n)<n^{-1}d(x_n,Q)$ and the function
			$x\mapsto \vf(x)+n^{-1}d(x,x_n)$  attains a local  minimum at $x_n$ for each
			$n$.
		\end{theorem}
		
		We shall speak about  {\it regular case} if the first option takes place and
		{\it singular}  or {\it non-regular case} otherwise.
		
		\vskip 1mm
		
		\proof  Indeed, either there is an $R>0$ such that $R\vf(x)\ge d(x,Q)$ for all
		$x$ of a neighborhood of $\xb$, or there is a sequence $(z_n)$ converging to
		$\xb$ and such that $n^2\vf(z_n)<d(z_n,Q)$. In the first  case (as $f$ is
		Lipschitz) we have for $x\not\in Q$ and
		$u\in Q$ close to $x$ (so that e.g. $d(x,u)< 2d(x,Q)$:
		$$
		f(x)\ge f(u)-Ld(x,u)\ge f(\xb)-2LR\vf(x),
		$$
		if $L$ is a Lipschitz constant of $f$.
		
		As $X$ is complete and $\vf$ is lower semicontinuous,  we can apply Ekeland's
		principle to $\vf$ (taking into account that $\vf(z_n)<\inf \vf +
		n^{-2}d(z_n,Q)$) and find
		$x_n$ such that $d(x_n,z_n)\le n^{-1}d(z_n,Q)$, $\vf(x_n)\le \vf(z_n)$
		and $\vf(x)+n^{-1}d(x,x_n)>\vf(x_n)$ for $x\neq x_n$. We have finally
		$$
		d(x_n,Q)\ge d(z_n,Q)-d(x_n,z_n)\ge(1-n^{-1})d(z_n,Q)\ge
		(1-n^{-1})n^2\vf(z_n)\ge n\vf(x_n)
		$$
		as claimed.\endproof
		
		Thus, a constrained problem reduces to one or a sequence of unconstrained
		minimization problems.
		Hopefully, such problems can be easier to analyze thanks to the freedom of
		choosing $\vf$ which we call {\it test function} in the sequel.
		Even before the alternative was explicitly stated  it
		was de facto used  to prove the maximum principle in various problems of optimal
		control
		\cite{GI93,AI97,RV}. Here is a brief account of how the alternative works for
		optimal control of systems governed by differential inclusions.
		
		\subsubsection{Optimal control of differential inclusion.}
		As the first example of application of the alternative we shall briefly consider the
		following
		problem of optimal control of a system governed by differential inclusion (see also the next subsection 8.3):
		minimize
		\begin{equation}\label{cost}
		\ell(x(0),x(T))
		\end{equation}
		on trajectories of the differential inclusion
		\begin{equation}\label{difinc}
		\dot x\in F(t,x),
		\end{equation}
		satisfying the end point condition
		\begin{equation}\label{end}
		(x(0),x(T))\in S.
		\end{equation}
		
		The natural space to treat the problem is $W^{1,1}$. Let $\xb(\cdot)$ be a local
		solution. For any $x(\cdot)\in W^{1,1}$ set
		$$
		\vf(x(\cdot)) = \int_0^T d(\dot x(t),F(t,x(t)))dt + d((x(0),x(T)), S).
		$$
		Clearly, $\vf$ is nonnegative and $\vf(\xb(\cdot))=0$. Thus, if $\ell$ is a
		Lipschitz function, we can apply the alternative to get necessary optimality
		condition. According to the alternative, either there is a $\la >0$ such that
		$\xb(\cdot)$ is a local minimum
		of
		$$
		\la\ell(x(0),x(T))+  d((x(0),x(T)), S) + \int_0^T d(\dot x(t),F(t,x(t)))dt,
		$$
		or there is a sequence $(x_n(\cdot))$ converging to $\xb(\cdot)$ such that every
		$x_n(\cdot)$ is not feasible in (\ref{cost})-(\ref{end}) and is a local minimum
		of the functional
		$$
		d((x(0),x(T)), S) + \int_0^T d(\dot x(t),F(t,x(t)))dt + n^{-1}\Big(
		\| x(0)-x_n(0)\| +\int_0^T\| \dot x(t)-\dot x_n(t)\|dt\Big)
		$$
		In both cases we get an (unconstrained) Bolza problem. Analysis of such problem
		needs different techniques and we refer to  \cite{AI97,RV} where necessary optimality conditions
		for the problem were obtained along these lines. A more general
		result was established a few years later by Clarke \cite{FHC05}(actually the most general for optimal control of differential inclusions so far) but a shorter proof of Clarke's theorem based on optimality alternative is now also available
		\cite{AI15}.
		
		To conclude, I wish to note that this is not the only possible application of regularity related ideas to optimal control. We can refer to \cite{RV05}
		for the discussion of the role of metric regularity in the Hamilton-Jacoby theory of optimal control.

		\subsubsection{Constraint qualification.}
		The last question we intend to briefly discuss in this subsection concerns
		constraint qualifications in optimization problems. They  often play an important
		role in proofs but their basic function is to guarantee that the multiplier
		$\la$ of the cost function is in the necessary (e.g. Lagrangian)  optimality
		conditions  is positive. The point is that constraint qualifications are often connected with regularity properties of the constraint mapping. We shall discuss just one example.

		Let us say that the problem is {\it normal} at a certain feasible point  if the
		constraint mapping is regular at the point. The {\it problem is normal} if either the feasible set is empty or the problem is normal at every feasible point.
		In the  case of the problem (\ref{P1}) the constraint mapping
		is the restriction of $F$
		to $C$, so by Theorem \ref{F7} normality is guaranteed if $F$ is transversal to
		$Q$, that is if
		$y^*\in N(Q,F(x))$ and $0\in D^*F|_C(\xb,0)(y^*)$ imply together that $y^*=0$
		which in turn imply that
		\begin{equation}\label{P5}
		0\in \sd(y^*\circ F)(x) + N(C,x),\quad\&\quad y^*\in N(Q,F(x))\; \Rightarrow
		y^*=0.
		\end{equation}
		This is the  now standard constrained qualification in nonsmooth optimization
		(see e.g. \cite{DR,KK02,BM,RW}). If $f$ and  $F$ are continuously differentiable and the sets $C$ and $Q$ are  convex,
		 (\ref{P5}) is dual to
		Robinson's constraint qualification \cite{SMR76}.

		\subsection{An abstract relaxed  optimal control problem.}
		Here we apply the optimality alternative to get necessary optimality condition in the problem
		\begin{equation}\label{glawyp}
		\text{\rm minimize}\quad f(x) \quad \text{\rm s.t.}\quad F(x,u)=0,\; u\in U.
		\end{equation}
		Here $F: X\times U\to Y$, $X$ and   $Y$ are separable Banach spaces and $U$ is a
		set.
		The problem is similar to  problems with mixed
		smooth and convex structures studied in \cite{IT,VMT82}. But contrary to
		\cite{IT,VMT82},
		here we do not assume that $F$ is continuously differentiable in $x$.
		We shall formulate the requirements on $F$ a bit later. First we need to
		introduce
		and discuss  some necessary concepts.
		
		We say that a continuous mapping $F: X\to Y$ is {\it semi-Fredholm} at $\xb$
		it has at $\xb$ a strict prederivative of the form $\ch(x)= Ax+ \| h\|Q$, where
		$A: X\to Y$ is a linear bounded operator that send $X$ onto a closed subspace
		of $Y$ of finite codimension and $Q\subset Y$ is a compact set (that can be
		assumed
		convex and symmetric).
		We say furthermore that $S\subset X$ is {\it finite-dimensionally generated} if
		$S=\Lambda^{-1}(P)$ where  $\Lambda: X\to R^n$ is a continuous linear operator
		and
		$P\subset \R^n$ is closed.

		\begin{proposition}[non-covering principle for (\ref{glawyp}) \cite{AI87,GI93}]\label{regfred}
			Let $F: X\to Y$ be semi-Fredholm at $\xb$, and let $S$ be a
			finite-dimensionally generated
			subset of $X$. Let further $F|_S$ be the restriction of $F$ to $S$, that is the
			set-valued mapping equal to $\{F(x)\}$ on $S$ and $\emptyset$ outside of $S$. If
			$F|_S$
			is not regular near $\xb$, then there is a $y^*\neq 0$ such that
			$0\in\sd_G(y^*\circ F)(\xb)+N_G(S,\xb)$. Moreover, the weak$^*$-closure of the
			set of
			such $y^*$ with norm $1$ does not contain zero\footnote{More general versions of this result can be found in many publications related to
				``point estimates" and compactness properties of subdifferentials - see e.g \cite{AI89a,JT95,JT95a,JT99,BM} }.
		\end{proposition}
	
		We intend to use this principle to prove the following theorem.
		\begin{theorem}\label{lagrange}
			Let $(\xb,\ub)$  be a solution of  (\ref{glawyp}). We assume that
			
			({\bf A$_1$}) $f$ satisfies the Lipschitz condition in a neighborhood of $\xb$;
			
			({\bf A$_2$}) for any $u\in U$ the mapping $F(\cdot,u)$  is Lipschitz in a
			neighborhood of $\xb$, and  $F(\cdot,\ub)$  is semi-Fredholm at  $\xb$;
			
			({\bf A$_3$})  $F(x,U)$ is a convex set for any $x$  of a neighborhood of
			$\xb$;
			
			({\bf A$_4$}) $S$ is finite-dimensionally generated

			Let further $\cll(\la,y^*,x,u)=\la f(x)+\lan y^*,F(x,u)\ran$ be the Lagrangian
			of the problem.
			Then there are $\la\ge 0$  and $y^*\in Y^*$ such that the following relations
			hold true:
			$$
			\begin{array}{ll}
			\la + \| y^*\|>0 & \text{(non-triviality)};\\
			0\in \sd_G\cll(\la,y^*,\cdot,\ub)(\xb) + N_G(S,\xb)& \text{(Euler-Lagrange
				inclusion)}  ;\\
			\lan y^*,F(\xb,\ub)\ran\ge\lan y^*,F(\xb,u)\ran,\quad\forall \; u\in U &
			\text{(the maximum principle)}  .
			\end{array}
			$$
		\end{theorem}
		
		\proof Given a finite collection $\cu=(u_1,\ldots,u_k)$
		of elements of $U$, we define a mapping   $\Phi_{\cu}: X\times\R^k\to Y$ by
		$$
		\Phi_{\cu}(x,\al_1,\ldots,\al_k)=
		F(x,\ub)+\sum_{i=1}^{k}\al_i(F(x,u_i)-F(x,\ub)).
		$$
		It is an easy matter to see that this mapping is also semi-Fredholm at $(\xb,0)$.
		
		Consider the problem
		$$
		\text{\rm minimize}\; f(x) \quad \text{\rm
			s.t.}\;\Phi_{\cu}(x,\al_1,\ldots,\al_k)=0,\;
		x\in S,\; \al_i\ge 0.
		\eqno{(\bf P_{\cu})}
		$$
		Then  $(\xb,0,\ldots,0)$  solves the problem (as immediately follows from ({\bf
			A$_3$})).
		Let further $\Psi: X\times \R^k \to Y$ be defined by
		$$
		\Psi(x,\al_0,\ldots,\al_k)= (f(x)+\al_0,\Phi_{\cu}(x,\al_1,\ldots,\al_k)).
		$$
		This mapping cannot be regular in a neighborhood of $(\xb,0,\ldots,0)$ because
		no point $(f(\xb)-\ep, 0,\ldots,0)$ can be a value of $\Psi$ at $x\in S$ close
		to
		$\xb$ and $\al$ close to zero. It is an easy matter to verify that $\Psi$
		is also
		semi-Fredholm at $(\xb,0,\ldots,0)$ and we can apply Proposition \ref{regfred}.
		
		Set
		$\hat S= S\times \R_-^{k+1},\; \hat{\cll}(\la,y^*,x,\al_0,\ldots,\al_k)=
		\la(f(x)+\al_0)+ \lan y^*,\Psi(x,\al_0,\ldots,\al_k)\ran$. By the proposition
		there are
		multipliers $(\la,y^*)\neq 0$  such that
		$$
		0\in\sd_G\hat{\cll}(\la,y^*,\cdot)(\xb,0,\ldots,0)+ N_G(\hat
		S,(\xb,0,\ldots,0)).
		$$
		We have (using the standard rules of subdifferential calculus )
		$$
		\begin{array}{l}
		N_G(\hat S,(\xb,0,\ldots,0))= N_G(\xb,S)\times\R_-^{k+1};\\
		\sd_G\hat{\cll}(\la,y^*,\cdot)(\xb,0,\ldots,0)\subset
		\sd_G\cll(\la,y^*,\cdot,\ub)(\xb) \\   \qquad\qquad\qquad\qquad\qquad
		+(\la,\lan y^*,F(\xb,u_1)-F(\xb,\ub)\ran,\ldots,\lan y^*,
		F(\xb,u_i)-F(\xb,\ub)\ran).
		\end{array}
		$$
		
		It follows that there are $\xi_i\le 0,\; i=0,\ldots,k$   such that
		$$
		\begin{array}{l}
		0\in \sd_G \cll(\la,y^*,\cdot,\ub)(\xb)+N_G(S,\xb);\\
		\la=-\xi_0\ge 0;\\
		\lan y^*,F(\xb,u_i)-F(\xb,\ub)\ran = \xi_i\ge 0,\quad i=1,\ldots,k.
		\end{array}
		$$

		The relations remain obviously valid if we replace $\la,y^*$ by $r\la,ry^*$ with
		some positive $r$. Thus for any finite collection $(u_1,\ldots,u_k)\subset U$ we
		can find
		a pair of multipliers $(\la,y^*)$ satisfying the three above mentioned
		relations
		and the normalization condition $\la + \| y^*\|=1$. Let $\Omega(u_1,\ldots,u_k)$
		be
		the weak$^*$-closure of all such pairs. Then  $\Omega(u_1,\ldots,u_k)$ is
		weak$^*$-compact and by Proposition \ref{regfred} does not contain zero. It
		remains to notice
		that  the increase of the set $(u_1,\ldots,u_k)$  may result only in
		decrease of $\Omega(u_1,\ldots,u_k)$ and therefore there is a nonzero pair
		$\la,y^*$ common to all
		sets $\Omega(u_1,\ldots,u_k)$. \endproof

		\subsection{Genericity in tame optimization.}
		Here by ``tame optimization"
		we mean optimization problems with semi-algebraic data. We  consider
		the same class of problems as in  (\ref{P1}). This time however we are interested in the effects of perturbations and shall work with a family of problems depending on a parameter $p$:
		\begin{equation}\label{P3a}
		{\rm minimize}\quad f(x,p),\quad {\rm s.t.}\quad F(x,p)\in Q,\; x\in C.
		\end{equation}
		\noindent Here $x$ is an argument in the problem and $p$ is a parameter.
		So subdifferentials and derivatives that will appear below are always with respect to $x$ alone. If $p$ is fixed,  then we denote the corresponding problem by $\cp_p$.

		Before we continue, we have to mention that for a semi-algebraic set $S\subset \R^n$ the properties
		
		$\bullet$ \ $S$ is a set of first Baire category in $\R^n$;
		
		$\bullet$ \ $S$ has $n$-dimensional Lebesgue measure zero;
		
		$\bullet$ \ $\dim S< n$
		
		\noindent are equivalent. Thus, when we  deal with semi-algebraic objects e.g. in $\R^k$, the
		word ``generic" means "up to a semi-algebraic set of dimension smaller than $k$."

		We shall assume that $p$ is taken from an open set $P\subset \R^k$ and, as before, $x\in\R^n$ and $F$ takes values in $\R^m$. Our main assumption is
		that
		\vskip 1mm
		
		\centerline{\it the restriction $F|_C(x,p)$ of $F$ to $C$ is transversal to $Q$.}
		
		\vskip 1mm
		
		\noindent This is definitely the case when	$k=m$ and $F(x,p) = F(x)-p$. As to $F$ itself, we assume that it is continuous with respect to $(x,p)$ and locally
		Lipschitz in $x$. The sets $C$ and $Q$ as usual are assumed closed.

		\begin{theorem}[generic normality]\label{F8}
			Under the stated assumptions for a generic $p\in P$, the mapping $F|_C(\cdot,p)$
			is transversal to $Q$.  Thus for a generic $p$ the problem $\cp_p$ is normal.	
		\end{theorem}
		
		\proof The first statement is immediate from Theorem \ref{thomset}, while the second from the comments following the statement of Theorem \ref{F7}.\endproof
		
		Let us call a point $x$ feasible in $\cp_p$  a {\it critical point} of the problem if the non-degenerate Lagrangian necessary condition of 8.2.1
		$$
		0\in\sd(f + (y^*\circ F|_C))(x,p),\quad y^*\in N(Q,F(x,p))
		$$
		is satisfied. In this case the value of $f$ at $x$ is called a {\it critical value}
		of $\cp_p$.
		
		\begin{theorem}[generic finiteness of critical values]\label{F9}
			If under the stated assumptions, $\cp_p$ is normal, then the problem may have only finitely many critical values. Thus there is an integer $N$ such that
			for a generic $p$ the number of critical values in the problem does not exceed $N$.
		\end{theorem}

		\proof Consider the function
		$$
		\cll_p(x,y,y^*)= f(x,p) +\lan y^*,F|_C(x,p)-y\ran + i_Q(y).
		$$
		As follows from the standard calculus rules,
		$$
		\sd\cll_p(x,y,y^*) = \sd(f+ y^*\circ F|_C)(x,p)\times
		(N(Q,y)-y^*)\times\{F(x,p)-y\}.
		$$
		
		Thus, $(x,y,y^*)$ is a critical point of $\cll_p$ if and only
		if $F(x,p)=y$, $0\in N(Q,y)-y^*$, that is $y\in Q$ and $y^*\in
		N(Q,y)$, \ and $0\in \sd(f+y^*\circ F|_C)(x,p)$. In other words,
		$(x,y,y^*)$ is a  critical point of $\cll_p$ if and only if $x$
		is a feasible point in ({\bf P}), $y=F(x,p)$ and the necessary
		optimality condition is satisfied at $x$ with $y^*$ being the
		Lagrange multiplier. We also see that in this case
		$\cll_p(x,y,y^*)=f(x,p)$. In other words, critical values of the
		problem are precisely  critical values of $\cll$.
		
		By the  Sard theorem
		$\cll_p$ may have  at most  finitely many
		critical values, whence the theorem.\endproof
		
		The last result we are going to present here has been so far proved only
		under some additional assumptions  on elements of the problem. We shall explain it for the classical case, although semi-algebraic nature of the data remains
		crucial.

		\begin{theorem}[generic finiteness of critical points]\label{F10}
			Assume that $p=(q,y)$ with $q\in\R^n$ and $y\in\R^m$, $f(x,p)= f(x)-\lan q,x\ran$, $F(x,p)=F(x)-y$ with $f(x)$ and $F(x)$ both continuously differentiable
			Assume further that  the sets
			$C$ and $Q$ are closed and convex. Then there is an integer $N$ such that for a generic $p$ the number of pairs $(x,y^*)$, such that $x$ is a critical point in $\cp_p$ and $y^*$ a corresponding Lagrange multiplier does not exceed $N$.
		\end{theorem}
		
		The theorem follows from the two results below that contain valuable information about geometry of subdifferential mappings of semi-algebraic functions.
		
		\begin{proposition}[dimension of the subdifferential graph \cite{DL12}]\label{dime}
			The dimension of the graph of the subdifferential (no matter which, Fr\'echet, limiting or Clarke) mapping of a semi-algebraic function on $R^n$ is $n$.	
		\end{proposition}
		
		\begin{proposition}[finiteness of preimage \cite{AI11b,DL12}]\label{fine}
			Let $F:\R^n\rra \R^n$ be a semi-algebraic set-valued mapping such that
			$\dim(\gr F)\le n$. If $y$ is a regular value of $F$, then $F^{-1}(y)$ contains at most
			finitely many elements. Thus, there is an integer $N$ such that for a   generic $y$ the number of elements in $F^{-1}(y)$ cannot exceed $N$.	
		\end{proposition}
		
		To see how the propositions lead to the proof of the theorem, we note first
		that $D^*F|_C(x)(y^*)=F'(x)y^*+N_C(x)$ if $x\in C$, $F$ is smooth and $C$ convex.
		By Theorem \ref{tran1} $F|_C$ is transversal to $Q$ if and only if
		\begin{equation}\label{P10}
		x\in C,\; F(x)\in Q+y,\; 0\in F'(x)y^*+N_C(x),\; y^*\in N(Q,F(x)-y)\; \Rightarrow y^*=0,
		\end{equation}
		and by Theorem \ref{thomset} this holds for a generic $y$.
		
		Consider the function
		$$
		g(x,y)= f(x) + i_C(x)+i_Q(F(x)-y)
		$$
		By Proposition \ref{dime} the dimension of the graph of its subdifferential
		is $n+m$. Then so is the dimension of the graph of the mapping
		$$
		\Gamma (x,y^*)= \{(q,y):\; (q,y^*)\in\sd g(x,y)  \}.
		$$
		Now by the Sard theorem generic $(q,y)$ is a regular value of $\Gamma$
		so (Proposition \ref{fine})  for a generic $(q,y)$ there are finitely many
		$(x,y^*)$ such that $(q,y)\in \Gamma(x,y^*)$. Finally, if for such $(q,y)$
		the qualification condition (\ref{P10}) is satisfied, then
		$$
		\sd g(x,y) = \{(q,y^*):\; f'(x)+ (y^*\circ F(\cdot))'(x) + N(C,x),\;
		y^*\in N(Q,F(x)-y)\}
		$$
		(even if $Q$ is not convex - see again Exercise 10.8 in  \cite{RW}) which in particular means that $x$ is a critical point of $\cp_p$ and $y^*$ is a
		Lagrange multiplier in the problem.

		\subsection{Method of alternating projection.}
		This is one of the most popular methods to solve feasibility problem due to its simplicity and efficiency. The feasibility problem in its simplest form consists in finding a common point of two sets, say $Q$ and $S$. The recipe offered by
		the method of alternating projection is the following: starting with a certain
		$x_0$, we choose for $k=0,1,\ldots$
		$$
		x_{2k+1}\in\pi_Q(x_{2k}), \quad x_{2k+2}\in\pi_S(x_{2k+1}),
		$$
		where $\pi_Q(x)$ is the collection of points of $Q$ closest to $x$ etc.. 	 
		
		Von Neumann was the first to show in mid-30s (see \cite{JVN50}) that in case of two subspaces
		the method converges to a certain point in the intersection of two closed subspaces in a Hilbert space (depending of course on the starting point). Later in the 60s Bregman
		\cite{LB65} and Gubin-Polyak-Raik \cite{GPR67} applied it to convex subsets
		in $\R^n$. In particular it was shown in \cite{GPR67} that the convergence is linear if relative interiors of the sets meet. Later Bauschke and Borwein
		\cite{BB96} proved linear convergence if the sets are  subtransversal at any common point.
		
		But in computational practice the method was successfully applied even for nonconvex sets. The first explanation was given by Lewis, Luke and Malik \cite{LLM09}: if at a certain point $\xb$ in the intersection the sets are transversal and at least one of the sets
		is not ``too non-convex" in a certain sense (super-regular in the terminology of the authors) then linear convergence of alternating projections to a certain point common to the sets (not necessarily $\xb$) if the starting point is sufficiently close to $\xb$. And very recently it was shown by Druzviatskyj,  Ioffe and Lewis \cite{DIL15} that transversality alone guarantees linear convergence.
		In fact linear convergence was proved in \cite{DIL15} under a substantially weaker condition of
		``intrinsic transversality" of the sets, but we believe that geometric essence of the phenomenon	 is captured by the transversality $\Rightarrow$ linear convergence
		implication. The question whether linear convergence is guaranteed by
		subtransversality, as in the convex case, remains open (see \cite{HL12}).

		Here is a short proof of linear convergence under the transversality assumption. Set
		$$
		\vf (x,y)=i_Q(x)+i_S(y) + \| x-y\|.
		$$

		We claim that if $Q$ and $S$ are transversal at $\xb\in Q\cap S$, then
		there are $\kappa>0$ and $\del>0$ such that for any $x\in Q,\ y\in S$
		close to $\xb$
		$$
		\max\{|\nabla\vf(\cdot,y)|(x),|\nabla \vf (x,\cdot)|(y)\}\ge \kappa.
		$$
	    To this end, we first note that by Theorem \ref{F5}
	    $$
	    \theta = \sup\{\lan u,v\ran:\; u\in N(Q,\xb),\; v\in -N(S,\xb),\; \|u\|=\| v\|=1\}<1.
	    $$
	
	Fix a certain $\kappa\in(0,1)$ and assume that there are sequences $(x_n)\subset Q,\ (y_n)\subset S$,
		$x_n\neq y_n$, converging to $\xb$ and such that
		$$
		|\nabla\vf(\cdot,y_n)|(x_n)< \kappa,\quad |\nabla \vf (x_n,\cdot)|(y_n)<\kappa,
		$$
		that is the functions
		\begin{equation}\label{min}
		x\mapsto \vf (x,y_n)+ \kappa\| x-x_n\|\quad{\rm and}\quad y\mapsto \vf(x_n,y)+\kappa\|y-y_n\|
		\end{equation}
		attain local minima respectively at $x_n$ and $y_n$. This means that
		\begin{equation}\label{3}
		0\in w_n^*+\frac{x_n-y_n}{\|x_n-y_n\|}+\kappa B;\qquad 0\in z_n^*+\frac{y_n-x_n}{\|x_n-y_n\|}+\kappa B
		\end{equation}
		for some $w_n^*\in N(Q,x_n)$ and $z_n^*\in N(S,y_n)$. Thus, for any limit point  $(w^*,z^*)$
		of $(w_n^*,z_n^*)$, we have
		$$
		w^*= e+a,\qquad z^* = -e+b,
		$$
		where $\| e\|=1,\; \| a\|\le\kappa,\; \| b\|\le \kappa$. Consequently
		$$
		\theta\ge\frac{\lan e+a,e+b\ran}{\|e+a\|\|e+b\|}\ge \frac{(1-\kappa)^2}{(1+\kappa)^2}
		$$
		and we get
		\begin{equation}\label{ineq}
		\kappa\ge\frac{1-\sqrt{\theta}}{1+\sqrt{\theta}},
		\end{equation}
		This proves the claim.
		
		Then $\pi_Q(y)={\rm argmin}~\vf(\cdot,y)$ and the method of alternating projections can be written as follows:
		$$
		x_{n+1}\in {\rm argmin}~\vf(x_n,\cdot);\quad x_{n+2}\in{\rm argmin}~\vf(\cdot,x_{n+1}).
		$$
		We obviously have $|\nabla\vf(x_n,\cdot)|(x_{n+1})|=0$. For a given $x$ (not necessarily in Q), consider the function $\psi_x(y)=i_S(y)+\|x-y\|$. For any $c\in (0,1)$  condition
$|\nabla\psi_x|(x_{n+1})\le c$ obviously holds if
$$
\lan x-x_{n+1},x_n-x_{n+1}\ran\ge \sqrt {1- c^2}\| x-x_{n+1}\|\| x_n-x_{n+1}\|.
$$
 Take a $c<\kappa$, and let $K_c$ be the collection of $c$ satisfying the above inequality.
 This is an ice-cream cone with vertex at $x_{n+1}$. If $x\in Q\cap K_c$, then
 $\nabla\vf(\cdot,x_n+1)(x)\ge \kappa>c$. On the other hand, as is easy to see, the distance from
 $x_n$ to the boundary of $K_c$ is precisely $cr$, where $r=\| x_n-x_n+1\|$.
Applying Basic lemma for error bounds (Lemma \ref{baslem}), we conclude that there is an $x\in Q$ with $\vf(x,x_{n+1})\le \vf(x_n,x_{n+1}) - c\kappa\| x_{n+1}-x_n\|$. It follows that
		$$
		\|x_{n+2}-x_{n+1}\|= \vf(x_{n+2},x_{n+1})\le (1-c^2)\| x_{n+1}-x_{n+1}\|
		$$
		which is linear convergence of $(x_n)$.

		\subsection{Generalized equations.}
		By a generalized equation we mean the relation
		$$
		0\in f(x)+F(x),
		$$
		where $f$ is a single-valued  and $F: X\rra Y$ a set-valued mapping. Variational
		inequalities and necessary optimality conditions in constraint optimization
		with
		smooth cost and constraint functions are typical examples. The problem discussed
		in the theorem below is what happens with the set of solutions of the generalized
		equation if the single-valued term is slightly perturbed.
		
		\begin{theorem}[implicit function for generalized equations]\label{imgeq}
			Let $X$, $Y$ be metric spaces, and let $Z$ be a normed space.  Consider the
			generalized
			equation
			\begin{equation}\label{3.315}
			0\in f(x,p)+F(x),
			\end{equation}
			where $f: X\times P\to Z$ and $F: X\rra Z$.
			Let $(\xb,\pb)$ be a solution to the equation. Set $\zb=-f(\xb,\pb)$
   and	suppose that the following two properties hold:
			
			(a) Either $X$ or the graph of $F$ is complete in the product metric
			and  $F$ is regular near $(\xb,\zb)$ with $\sur F(\xb|\zb)>r$;
			
			(b) there is a $\rho>0$ such that $f$ is continuous on
			$\bo (\xb,\rho)\times \bo (\pb,\rho)$ and $ f(\cdot,p)$
			satisfies on $\bo (\xb,\rho)$ the Lipschitz condition with constant $\ell<r$
			for all
			$p\in\bo (\pb,\rho)$.
			
			Let $S(p)$ stand for the solution mapping of (\ref{3.315}). Then
			$$
			d(x, S(p'))\le (r-\ell)^{-1}\| f(x,p)-f(x,p')\|.
			$$
			if $x\in S(p)$ is close to $\xb$ and $p, p'$ are sufficiently close to $\pb$.
			Thus, if $f(x,\cdot)$ satisfies the Lipschitz condition with constant $\al$ on
			a neighborhood
			of $\pb$ for all $x\in\bo(\xb,\rho)$, then $S(\cdot)$  has the Aubin property
			near $(\pb,\xb)$ with
			$\lip S(\pb|\xb)\le\al(r-\ell)^{-1}$.
			
			Finally, if  in addition $F$ is strongly regular near $(\xb,\zb)$,
			then
			$S(\cdot)$ has a Lipschitz localization $s(\cdot)$ at $\xyb$ with Lipschitz
			constant not greater than
			$\al(r-\ell)^{-1}$, so that
			$$
			d(s(p), s(p'))\le (r-\ell)^{-1}\| f(s(p),p)-f(s(p),p')\|\le
			\al(r-\ell)^{-1}d(p,p').
			$$
		\end{theorem}
		
		Note that in view of Theorem \ref{critgen} condition (a) is equivalent to the
		assumption that
		there are $r>0$ and $\xi >0$ such that $|\nabla_{\xi}\vf_z|(x,v)>r$ (where
		$\vf_z(x,v)=d(z,v)+i_{\gr F}(x,v)$) if e.g.
		$d(z,\pb)<\rho,\;  \| z\|<\rho$ and $z\neq v\in F(x)$.
		
		\proof Set $G(x,p)=f(x,p)+F(x)$ and let $H(p,z)= (G(\cdot,p))^{-1}(z)$, so that
		$S(p) = H(p,0)$.
		As the Lipschitz constants of functions $f(\cdot,p)$ are bounded by the same
		$\ell$ for all
		$p\in \bo(\pb,\rho)$,  it follows from Theorem \ref{milt} that there is a
		$\del>0$ such that
		for every $p\in \bo (\pb,\rho)$  the inequality $d(x,H(p,z))\le
		(r-\ell)^{-1}d(z,G(x,p))$
		holds if $d(x,\xb)<\del$ and $\|z-z(p)\|<\del$, where $z(p)=f(\xb,p)-f(\xb,\pb)\in G(\xb,p)$. As
		$f$
		is continuous, we can choose $\la>0$ such that $\| z(p)\|<\del$ for
		$p\in\bo(\pb,\la)$.
		For such $p$ we have $0\in\bo(z(p),\del)$ and therefore if $d(\pb,p')<\la$, we
		get, taking into
		account that $0\in f(x,p)+F(x)$ by the assumption,
		$$
		\begin{array}{lcl}
		d(x,S(p'))&\le& (r-\ell)^{-1}d(0,G(x,p'))= (r-\ell)^{-1}d(0,f(x,p')+F(x))\\
		&=& (r-\ell)^{-1}d(-f(x,p'),F(x))\le(r-\ell)^{-1}\|f(x,p')-f(x,p)\|
		\end{array}
		$$
		
		This proves the first part of the theorem. The second now follows from Theorem
		\ref{mimpl}.
		\endproof
		
		The concept of generalized equation was introduced by Robinson in \cite{SMR79}.
		The theorem proved in \cite{SMR79,SMR80} corresponded to $f$ continuously
		differentiable in $x$ and $F$ being either a maximal monotone operator or
		$F(x)=N(C,x)$, where $C$ is a closed convex set.
		We refer to \cite{DR} for further results and bibliographic comments on
		generalized equations which is one of the central objects of interest in the
		monograph.
		
		An earlier version of part (a) of the theorem with a less precise estimate
		can be found in \cite{KK02} (Theorem 4.9). Part (b) of the theorem relating to
		 strong regularity
		is  the basic statement of  Theorem 5F.4 of \cite{DR} (generalizing the
		earlier results of Robinson in \cite{SMR80,SMR91};  see also \cite{ALD95} for
		an earlier  result). Our proof however is different:
		here the theorem appears as a direct consequence of Milyutin's perturbation
		theorem.
		Note that in most of the related results in \cite{DR} it is assumed (following
		\cite{SMR91}) that  there exists a ``strict estimator $h(x)$ for $f$ of modulus
		$\ell$"  such that $\sur(F+h)(x|\yb+h(\xb))\ge r$. This is a fairly convenient
		device for practical purpose but it adds no generality to the result as the case
		with $h$ reduces to the setting of the theorem if we replace $F+h$ by $F$ and
		$f-h$ by $f$.

		\subsection{Variational inequalities over polyhedral sets.}
		
		{\it Variational inequality} is a relation of the form
		\begin{equation}\label{vareq}
		0\in \vf(x)+N(C,x),
		\end{equation}
		where $\vf:\R^n\to\R^n$ is a single-valued mapping and $C\subset \R^n$ is a convex set. If $C$ is a cone, it is equivalent to
		$$
		x\in K,\quad F(x)\in K^{\circ},\quad \lan x,F(x)\ran =0.
		$$
		The problem of finding such an $x$ is known as a {\it complementarity problem}
		(see e.g. \cite{FP}). Problems of this kind typically appear in nonlinear
		programming in connection with necessary optimality conditions.

	    Consider for instance the problem
	    \begin{equation}\label{nesop}
	    {\rm minimize}\quad f_0(x)\quad{\rm s.t}\quad f_i(x)\le 0,\ i=1,\ldots,k,\
	    f_i(x)\le 0,\ i=k+1,\ldots,m.
	    \end{equation}
	    with $f_0,\ldots,f_m$ twice continuously differentiable.
	    If $\xb$ is a solution of the problem, then (assuming that the problem is normal
	    and setting $f = (f_1,\ldots,f_m)$)
	    there is a $\yb\in \R^m$ such that
	    $$
	    \nabla f_0(\xb)+\lan \yb,\nabla f(\xb)\ran =0.
	    $$
	    Setting
	    $$
	    \vf(x,y)=\left(\begin{array}{c} \nabla f_0(\xb)+\lan \yb,\nabla f(\xb)\ran,\\ f(x)\end{array}\right);\qquad C=\R^n\times \R^m_+,
	    $$
	    we see that $\xyb$ solves (\ref{vareq}) (with $x$ replaced by $(x,y)$).
	
	    Consider the set valued mapping $\Psi(x)= \vf(x)+N(C,x)$ associated with
	    (\ref{vareq}) assuming that $C$ is a convex polyhedral set. What can be said
	    about regularity of such mapping near a certain $\xyb\in\gr\Phi$? Applying
	    Milyutin's perturbation theorem (Theorem \ref{milt}) and Theorem \ref{strlip}
	    and taking into account that the Lipschitz constant of $h\to \vf(x+h)-\phi'(x)h$
	    at zero is zero,	  we immediately get
	
	    \begin{proposition}\label{lvar}
	    Let $\yb\in\Psi(\xb)$ for some $\xb\in C$. Set $A=\vf'(\xb)$ and
	    $\hat{\Psi}(x) = Ax +N(C-\xb,x)$. Then $\Psi$ is (strongly) regular
	    near $\xyb$ if and only if $\hat{\Psi}$ is (strongly) regular near
	    $(0,0)$ and $\sur\Psi(\xb|\yb)=\sur\hat{\Psi}(0|0)$.
	    \end{proposition}
	
	    In other words, the regularity properties of $\Psi$ are the same as
	    of its ``linearization" $\hat{\Psi}$. Therefore in what follows we can deal only
	    with the {\it linear variational inequality}
	    \begin{equation}\label{lvareq}
	    0\in Ax+N(C,x)
	    \end{equation}
	    and the associated mapping
	    $$
	    \Phi(x) = Ax+N(C,x).
	    $$
	
	    The key role in our analysis is played by the concept of a {\it face}
	    of a polyhedral set $C$ which is any closed subset $F$ of $C$ such that
	    any segment $\Delta\subset C$ containing a point $x\in F$ in its interior
	    lies in $F$.
	    A face of $C$ {\it proper} if it is different from $C$.
	    We refer to \cite{RTR} for all necessary information about faces. The following
	    facts are important for our discussion:
	
	    $\bullet$ the set $\cf_C$ of all faces of $C$ is finite;
	
	    $\bullet$ $F\in\cf_C$ if and only if there is a $y\in\R^n$ such that $F=\{x\in C:\; \lan y,x\ran\ge\lan y,u\ran,\;\forall\; u\in C  \}$;

	    $\bullet$ if $F,F'\in\cf_C$ and  $F\cap {\rm ri}~\!F'\neq\emptyset$, then $F'\subset F$; a proper face of $C$ lies in the relative boundary of $C$;

	    $\bullet$ if $F\in\cf_C$ and $x_1,\ x_2$ belong to the relative interior of $F$,
	    then $T(C,x_1)= T(C,x_2)$ and $N(C,x_1)= N(C,x_2)$.
	
	    The last property allows to speak about the tangent and normal cones to
	    $C$ at $F$ which we shall denote by $T(C,F)$ and $N(C,F)$.
	    It is an easy matter to see that
	    \begin{equation}\label{dimpol}
	    \dim F+\dim N(C,F)= n;\qquad \dim(F+N(C,F))=n.
	    \end{equation}
	
	    For any $x\in C$  denote by $F_{\min}(x)$ the minimal element of $\cf_C$
	    containing $x$. The   is straightforward
	    \begin{equation}\label{ri}
	    	x\in F\in \cf_C, \;\&\; F=F_{\min}(x)\; \Leftrightarrow\; x\in{\rm ri}~\! F.
	    \end{equation}

	    \begin{proposition}\label{nonsing}
	    If $\Phi$ is regular near $(x,y)$ and  $F= F_{\min}(x)$, then\\
	    $$\dim (A(F)+ N(C,F))=n.$$
	    In particular, $A$ is one-to-one on $F$.	
	    \end{proposition}
	
	   \proof If $\dim F=0$, then $x$ is an extreme point of $C$ in which case
	   $T(C,x)$ is a convex cone containing no lines and its polar therefore has
	   nonempty interior. On the other hand, if $x\in\intr C$, then $N(C,u)=\{0\}$
	   for all $u$ of a neighborhood of $x$ and $\Phi(u)=Au$ for such $u$. So by regularity $A$ is an isomorphism.

	   Thus in the sequel we may assume that the dimensions of both $F$ and
	   $N(C,F)$ are positive.	
	   By changing $(x,y)$ slightly, we can guarantee that $y$ belong to the
	   relative interior of $N(C,F)$.
	   Let $\ep>0$ be so small that the distances from $x$ and $y$ to the relative boundaries of $F$ and $N(C,F)$ are greater than $\ep$.
	   Then any $(u,v)$ such that $u\in C$, $v\in N(C,u)$, $\| u-x\|<\ep,\ \| v-y\|<\ep$ must belong to
	   $F\times N(C,F)$. This means that
	   $\Phi(B(x,\ep))\cap B(y,\ep)\subset A(F)+N(C,F)$ and the result follows from
	   (\ref{dimpol}). Indeed, the dimension equality is immediate from the last inclusion. On the other hand, if $A$ is not one-to one on $F$, then
	   $\dim A(F)<\dim F$ and by (\ref{dimpol}) $\dim A(F)+\dim N(C,F)<n$.
	   \endproof

	   Let $C\subset \R^n$ be a convex polyhedron, and let $F$ be a proper face of $C$.
	   Let $L$ be the linear subspace spanned by $F$ and $M$ the linear subspace
	   spanned by $N(C,F)$. These subspaces are complementary by (\ref{dimpol}) and orthogonal.
	   By Proposition \ref{nonsing} $A(L)$ and $M$ are also complementary subspaces
	   if $\Phi$ is regular near any point of the graph.

	   Let $\pi_M$ be the projection onto $M$ parallel to $A(L)$, so that $\pi_M(A(F))=0$. Set
	   $K_M=(T(C,F))\cap M$, and let $A_M$ be the restriction of $\pi_M\circ A$
	   to $M$.
	   Then $K_M$ is a convex polyhedral cone in $M$ and its polar $K_M^{\circ}$ (in $M$) coincides with $N(C,F)$.
	
	   \begin{definition}\label{fact}
	   	{\rm  The set-valued mapping $\Phi_M(x) = A_Mx + N(K_M,x)$ viewed as a mapping from	$M$ into $M$  will be called {\it factorization of $\Phi$ along $F$}. }
	   \end{definition}
	   Observe that the graph of a factorization mapping is a union of convex polyhedral
	   cones.
	   \begin{proposition}\label{regfact}
	   	If $\Phi$ is regular near $(\xb,A\xb)$ for some $\xb\in C$,
	   	then the factorization of $\Phi$ along $F=F_{\min}(\xb)$ is
	   	globally regular on $\R^n$.
	   	\end{proposition}
	   	
	   	\proof  Set $K_1=T(C,F)=T(C,\xb)$ and consider the mapping $\Phi_1(x)= Ax+N(K_1,x)$. By Proposition \ref{loctan},
	   	$\Phi_1(x)=\Phi (\xb+x)- A\xb$ for $x$ close to zero. Therefore $\Phi_1$ is regular near $(0,0)$, hence globally regular by Proposition \ref{loctan}. Observe that $K_1=K_M+L$ and  $K_1^{\circ}= N(K,F)$ and consequently $N(K_1,x)\subset N(K,\xb)=N(K,F)$ for any $x\in K_1$.
	   	
	   	As $\Phi_1$ is globally regular, there is a $\rho>0$ such that
	   	$d(x,\Phi_1^{-1}(z))\le \rho d(z,\Phi_1(x))$ for all $x,z\in\R^n$. Take now
	   	$x,z\in M$. We have (taking into account that $N(K_M,x)=N(K_1,x+\xi)$ for any $\xi\in L$ and $A_Mx = A(x+\xi)$ for some $\xi\in L$)
	   	$$
	   	\begin{array}{lcl}
	   	d(z,\Phi_M(x))&=&\inf\{\|z-A_Mx-y\|:\; y\in N(K_M,x)    \}\\
	   	&\ge& \inf\{\|z-A(x+\xi)-y\|:\;\xi\in L,\;  y\in N(K_1,x+\xi)    \}\\
	   	&=& \dis\inf_{\xi\in L}d(z,\Phi_1(x+\xi))=d(z,\Phi_1(w))
	   	\end{array}
	   	$$
	   	for some $w\in x+L$. On the other hand, there is  a $w'\in\R^n$ such that
	   	$z\in\Phi_1(w')$ and $\| w-w'\|=d(w,\Phi_1^{-1}(z))$. Let $x'$ be the orthogonal
	   	projection of $w'$ to $M$. We have $z=Aw'+y$ for some $y\in N(K_1,w')\subset M$.
	   	Therefore $Aw'\in M$ and moreover $A_Mx'=Aw'$.  The latter is a consequence of the following simple observation:
	   	\begin{equation}\label{observe}
	   	v=Aw\in M,\quad x\in M,\; x\perp(w-x)\; \Rightarrow\; A_Mx= v.
	   	\end{equation}
	   	Indeed, $z=w-x\in L$, hence $Ax= Aw+Az= v+Az$ and, as $v\in M$ and $Az\in A(L)$
	   	we have $\pi_M(Ax)= v +\pi_M(Az)=v$.

	   	It follows,  as $N(K_M,x')= N(K_1,w')$), that $z\in\Phi_M(x')$ and
	   $$
	   	d(x,\Phi_M^{-1}(z))\le \|x-x'\|\le\| w-w'\|=d(w,\Phi_1^{-1}(z))\le\rho d(z,\Phi_1(w))\le d(x,\Phi_M(x)),
	   	$$
	   	that is $\Phi_M$ is regular on $M$ (with the  rate of metric regularity not greater  than $\rho$).
	   	\endproof

	   The following theorem is the key observation that paves way for proofs of the main result.
	
	   \begin{theorem}\label{geom}
	   	Let $C=K$ be a convex polyhedral cone. If $\Phi$ is regular near $(0,0)$ (hence
	   	globally regular by Proposition \ref{globa}), then $A(K)\cap K^{\circ} =\{0\}$.
	   \end{theorem}
	
	   \proof The result is trivial if $n=1$. Assume that it holds for $n=m-1$, and
	   let $m=n$.
	   Note that the inclusion $A(K)\subset K^{\circ}$ can hold only if $K=\{0\}$.
Indeed, if the inclusion is valid, then $\Phi(x)\in A(K)+ K^{\circ}=K^{\circ}$ for any $x\in K$, so by regularity
$K^{\circ}$ must coincide with the whole of $\R^n$ and hence $K=\{0\}$. Thus if there is a nonzero
$u\in A(K)\cap K^{\circ}$, we can harmlessly assume that $u$ is a boundary point of $K^{\circ}$ and there is a nonzero $w\in N(K^{\circ},u)$. Then $w\in K$
and $u\in N(K,w)$. Let $F=F_{\min}(w)$ so that $u\in N(K,F)$. Let as before, $L$ be the linear subspace spanned by $F$ and $M$ the linear subspace
spanned by $N(K,F)$. These subspaces are complementary by (\ref{dimpol}) and orthogonal.
By Proposition \ref{nonsing} $A(L)$ and $M$ are also complementary subspaces.
Clearly, $u$ does not belong either to $L$ or to $A(L)$, the latter because
otherwise the dimension of $A(F)+ N(K,F)$ would be strictly smaller than $n$.
	
Consider the factorization $\Phi_M$ of $\Phi$ along $F$. 	
Then $u\in K_M^{\circ}$ by definition. But as follows from (\ref{observe}) $u$ also belongs to $A_M(K_M)$. As $\Phi_M$ is regular by Proposition \ref{regfact} and $\dim M<m$,
the existence of such a $u$ contradicts to the induction hypothesis.

\endproof

We are ready to state and proof the main result of the subsection.
\begin{theorem}[regularity implies strong regularity]\label{main}
	Let $C$ be a polyhedral set and  $\Phi(x) = Ax+N(C,x)$. If $\Phi$ is
	globally  regular then the
	inverse mapping $\Phi^{-1}$ is single-valued and Lipschitz on $\R^n$. Thus,  global regularity
	of $\Phi$ implies global strong regularity.
\end{theorem}

In other words, the solution map of $y\in \Phi(x)$ is everywhere single-valued and Lipschitz.

\proof
 We only need to show that $\Phi^{-1}$ is single-valued: the Lipschitz property
 will then automatically follow from regularity. The theorem is trivially valid if $n=1$.
 Suppose it is true for $n\le m-1$ and consider the case $n=m$.
 We have to show that, given a convex polyhedron $C\in\R^m$ and a linear operator $A$ in $\R^m$ such that $\Phi(x)= Ax+ N(C,x)$ is globally regular on $\R^n$, the equality
 $Ax+y=Au+z$ for some $x,u\in C$, $y\in N(C,x)$, $z\in N(C,u)$ can hold only if $x=u$ and $y=z$.

  \noindent {\bf Step 1}. To begin with we observe that
  the equality $Au=Ax+y$ for some $u,x\in C$ and $y\in N(C,x)$ may hold only if $u=x$. Indeed, $u-x\in T(C,x)$. The same argument as in the proof of Proposition \ref{regfact} shows that $\Phi_1(w)= Aw+ N(T(C,x),w)$ is also globally regular and therefore by Theorem \ref{geom}	$A(T(C,x))\cap N(C,x)=\{0\}$.
  It follows that  $A(u-x)= y=0$. But regularity of
  $\Phi_1$ implies (by Proposition \ref{nonsing}) that $A$ is one-to one on $T(C,x)$,
  hence $u=x$.

 \noindent {\bf Step 2}.
 Assume now that for some $x,u\in C,\; u\neq x$, the equality
 $Ax+y=Au+z$, or $A(u-x)=y-z$,  holds with $y\in N(C,x)$, $z\in N(C,u)$. We first show that this is impossible if $x\in F_{\min}(u)$. If under this condition $x\in{\rm ri}~\!C$, then $u$ is also in ${\rm ri}~\!C$ which means that $N(C,x)=N(C,u)$ coincides with the orthogonal
 complement $E$ to the subspace spanned by $C-C$. We have $y-z\in E$ and
 $u-x\in C-C$. By Proposition \ref{nonsing} $A(u-x)=y-z=0$ and the second part of
 the proposition implies that $u=x$.

 Let now $F=F_{\min}(x)$ be a proper face of $C$.  Then $F\subset F_{\min}(u)$ and therefore $z\in N(C,F)$.
 Denote as before by  $L$  the subspace spanned by $F$ and by $M$ the subspace spanned by $N(C,F)$, and let $\Phi_M$ be the factorization of $\Phi$ along $F$.   Set
 $v=A(u-x)=y-z$. Then $v\in M$ as both $y$ and $z$ are in $N(C,F)$. Let $w$ be the orthogonal
 projection of $u-x$ onto $M$. Then by (\ref{observe}) $Aw=v$ and therefore $A_Mw=v$.

 Thus (recall that $y,z\in M$)
 $$
 A_Mw+z= (\pi_M\circ A)(u-x) +z=\pi_M(A(u-x)+z)= \pi_My=y.
 $$
 On the other hand, it is clear that $y\in N(K_M,0)$ and $z\in N(K_M,w)$. Indeed,
 $z\in N(T(C,x),u-x)$ (since $\lan z,v-x\ran \le\lan z,u-x\ran$ for all $v\in C$ on the one hand and, as we have seen,  $z\in N(C,x)$, on the other) and therefore
 $z\in N(K_M,w)$ as $z\in M$ and $w-(u-x)\in L$. As $\dim M<m$, we conclude by the induction  hypotheses that $w=0$, hence $u-x\in L$. But $A(u-x)=y-z\in M$ and
 a reference to proposition \ref{nonsing} again proves that  $u=x$.

 \noindent{\bf Step 3}.
 It remains to consider the case when neither  $x$ nor $u$ belongs to the minimal face of the other.
 Let $\kappa$ be the modulus of metric regularity of $\Phi$ or any bigger number.
 Choose $\ep>0$ so small that the ball of radius $(1+\kappa)\ep$ around $x$ does not meet any face
 $F\in\cf_C$  not containing $x$.
 This means that $x\in F_{min}(w)$ whenever $w\in C$ and $\| w- x\|\le (1+\kappa)\ep$.
 Let further $N$ be an integer
 big enough to guarantee that $\del= N^{-1}\| y\|<\ep$. Regularity of $\Phi$ allows to construct recursively
 a finite sequence of pairs $(u_k,z_k), \ k=0,1,\ldots,m$ such that
 $$
 (u_0,z_0)=(u,z),\quad z_k\in F_{max}(u_k),\quad u_k+z_k= x+(1-m^{-1}k)y,\quad \|u_k-u_{k-1}\|\le \kappa\del.
 $$
 Then $u_N+z_N=  x$. As follows from the result obtained at the first step of the proof,
 this means that $u_N=x$. This in turn implies, as $u_0\neq x$,  that for a certain $k$ we have $u_k\neq x,\; \| u_k-x\|\le \kappa\del< \kappa\ep$. By the choice of $\ep$ this implies that
 $x\in F_{min}(u_k)$. But in this case the result obtained at the second step
 excludes the possibility of the equality $u_k+z_k= x+(1-m^{-1}k)y$
 unless $u_k=x$. So we again get a contradiction that completes the proof.\endproof

The material presented in this subsection is a part of my recent paper \cite{AI15a}
which contains also a proof (based on a similar ideas) of another principal result concerning uniqueness and Lipschitz behavior of solutions to variational inequalities over polyhedral sets due to Robinson \cite{SMR92}. Theorem \ref{main} was first stated
by Dontchev-Rockafellar \cite{DR96} with a comment that it follows from a comparison
of the mentioned Robinson's result and another theorem (proved by Eaves and Rothblum \cite{ER90}) containing an openness criterion for piecewise affine mappings. The
given proof seems to give the first self-contained and reasonably short justification for the result. We refer to \cite{DR,FP} for further details.


		\subsection{Differential inclusions -- existence of solutions.}
		Here we consider the Cauchy problem for differential inclusions:
		\begin{equation}\label{di1}
		\dot x\in F(t,x),\quad x(0)=x_0,
		\end{equation}
		where $F:\R\times\R^n\rra\R^n$. We assume that
		
		$\bullet$ \ $F$ is defined on some $\Delta\times U$ (that is
		$F(t,x)\neq\emptyset$ for all $x\in U$ and almost all  $t\in\Delta$), where
		$\Delta = [0,T]$ and $U$ is an open subset of $\R^n$ containing  $x_0$;
		
		$\bullet$ \ the graph of $F(t,\cdot)$ is closed for almost every $t\in \Delta$;
		
		$\bullet$ \ $F$ is measurable in $t$ in the sense that the function
		$t\mapsto d((x,y),\gr F(t,\cdot))$ is measurable for all pairs
		$(x,y)\in\R^n\times\R^n$.

		By a solution of (\ref{di1}) on $[0,\tau]\subset [0,\Delta]$ we mean any
		absolutely continuous $x(t)$ defined on  $[0,\tau]$ and such that
		$\dot x(t)\in F(t,x(t))$ almost everywhere on $[0,\tau]$.
		
		\begin{theorem}\label{exdif}
			Assume that  there is a summable $k(t)$ such that
			\begin{equation}\label{di2}
			h(F(t,x),F(t,x'))\le k(t)\| x-x'\|,\quad\forall\; x,x'\in U,\;\text{a.e.
				on}\quad [0,1].
			\end{equation}
			Let further $x_0(\cdot)$ be an absolutely continuous function on $[0,T]$
			with values in $U$ such that $x_0(0)=x_0$ and $\rho(t)= d(\dot
			x_0(t),F(t,x_0(t)))$ is a summable function.
			
			Then there is a solution of (\ref{di1}) defined on some $[0,\tau]$, $\tau>0$.
			Specifically, set $r=d(x_0,\R^n\backslash U)$, and let $\tau\in
			(0,T]$
			be so small that
			\begin{equation}\label{di3}
			1> k_{\tau}= \int_0^{\tau} k(t)dt;\quad(1-k_{\tau})r>
			\xi_{\tau}=\int_0^{\tau}d(\dot x_0(t),F(t,x_0(t)))dt.
			\end{equation}
			Then for any $\ep>0$ there is a solution $x(\cdot)$ of (\ref{di1}) defined on
			$[0,\tau]$
			and satisfying
			\begin{equation}\label{di4}
			\int_0^{\tau}\| \dot x(t)-\dot x_0(t)\|\le \frac{1+\ep}{1-k_{\tau}}\xi_{\tau}.
			\end{equation}
		\end{theorem}
		\noindent Recall that $h(P,Q)$ is the Hausdorff distance between $P$ and $Q$.
		
		\proof We may set  $x_0(t)\equiv0$ (replacing if necessary $F(t,x)$ by
		$F(t,x_0(t)+x)-\dot x(t)$ and $U$ by $r\bo$). Let $X=W_0^{1,1}[0,\tau]$ stand
		for the space of $\R^n$-valued absolutely continuous functions on $[0,\tau]$
		equal
		to zero at zero with the norm
		$$
		\| x(\cdot)\|_{\tau}=\int_0^{\tau} \| \dot x(t)\|dt,
		$$
		and let $I$ denote the identity map in $X$. Let finally $\cf$ be the set-valued
		mapping from $X$ into itself that associates with every $x(\cdot)$ the
		collection
		of absolutely continuous functions $y(\cdot)$ such that $y(0)=0$ and $\dot
		y(t)\in F(t,x(t))$ a.e..
		We have to prove the existence of an
		$x(\cdot)\in X$ satisfying (\ref{di4}) and
		\begin{equation}\label{di5}
		0\in (I-\cf)(x(\cdot))
		\end{equation}
		
		Note first that the graph of $\cf$ is closed, that is whenever $x_n(\cdot)\to
		x(\cdot)$, $y_n(\cdot)\in\cf(x_n(\cdot))$ and
		$y_n(\cdot)$ norm converge to $y(\cdot)$, then $y(\cdot)\in \cf(x(\cdot))$.
		Let $\cu$ be the open ball of radius $r$ around zero in
		$X$. Thus $x(t)\in U$ for any $t\in [0,\tau]$ whenever $x(\cdot)\in \cu$ and
		therefore by
		(\ref{di2}) $\cf$ is Lipschitz on $\cu$ with $\lip\cf(\cu)\le  k_{\tau}$. On the
		other hand, $I$ is Milyutin regular on $\cu$ with $\sur_mI(\cu)=1$. By Theorem
		\ref{milt1}
		\begin{equation}\label{di6}
		\sur_m(I-\cf)(\cu)\ge 1-k_{\tau}.
		\end{equation}
		In particular $B(y(\cdot),(1-k_{\tau})\rho)\subset (I-\cf)(\rho B)$ for any
		$y(\cdot)\in (I-\cf)(0)$ if $\rho< r$.
		Take a $y(\cdot)\in X$ such that $ \dot y(t)\in F(t,0)$ and
		$\|\dot y(t)\|=d(0,F(t,0))$ a.e.. Then $\| y(\cdot)\|_{\tau}=
		\xi_{\tau}<(1-k_{\tau})r$ by (\ref{di3}). Thus $0\in
		B(y(\cdot),(1-k_{\tau})\rho)$ for some $\rho<r$ and therefore there is an
		$x(\cdot)$ with $\| x(\cdot)\|_{\tau}<\rho$,\   $0\in(I-\cf)(x(\cdot))$.
		\endproof

		The theorem is close to the original result  of Filippov \cite{AF67}. Versions
		of this
		results and its applications can be found in many subsequent publications, see
		e.g
		\cite{AC,AF}. Typical proofs of existence results for differential inclusions
		use  either some iteration procedures or selection theorems to reduce the
		problem to
		existence of solutions of differential equations.  Observe that our proof
		appeals to
		non-local regularity theory.

	\end{document}